\documentclass[a4paper,12 pt,twoside, reqno]{amsart}

\usepackage[utf8]{inputenc}
\usepackage[english]{babel}
\usepackage{times}
\usepackage{amsmath}
\usepackage{amsfonts}
\usepackage{amssymb}
\usepackage{amsmath}
\usepackage{amsthm}
\usepackage{graphicx}
\usepackage{array}
\usepackage{color}
\usepackage{mathrsfs}
\usepackage{hyperref}
\usepackage{esint} %%% for \fint
\usepackage{tikz}
\usepackage{ upgreek }
\usepackage{enumitem}

\definecolor{darkgreen}{rgb}{0.1,0.7,0.1}
\definecolor{darkred}{rgb}{0.7,0.1,0.1}

\newcommand{\bin}[2]{\left (
\begin{array} {c}
#1 \\
#2
\end{array}
\right )}
\newcommand{\mat}[4]{\left( \begin{array}{cc}
#1 & #2  \\
#3 & #4  \\
\end{array} \right)}
\newtheorem{theorem}{Theorem}
\newtheorem{lemma}{Lemma}[section]
\newtheorem{proposition}[lemma]{Proposition}
\newtheorem{corollary}[theorem]{Corollary}
\newtheorem{conjecture}[lemma]{Conjecture}
\newtheorem{remark}[lemma]{Remark}

\newcommand{\eps}{\varepsilon}

\newcommand\symb[2][\bf]{{\mathchoice{\hbox{#1#2}}{\hbox{#1#2}}%
        {\hbox{\scriptsize#1#2}}{\hbox{\tiny#1#2}}}}

\def\R{{\symb R}}

\def\Z{{\symb Z}}
\def\C{{\symb C}}
\def\Q{{\symb Q}}
\def\P{{\symb P}}

\def\un{\mathbf{1}}

\renewcommand{\P}{\mathbb{P}}
\newcommand{\E}{\mathbb{E}}

\newcommand{\cA}{\mathcal{A}}
\newcommand{\cB}{\mathcal{B}}
\newcommand{\cC}{\mathcal{C}}
\newcommand{\cD}{\mathcal{D}}

\newcommand{\cF}{\mathcal{F}}

\newcommand{\cM}{\mathcal{M}}
\newcommand{\cN}{\mathcal{N}}

\newcommand{\cQ}{\mathcal{Q}}

\newcommand{\cS}{\mathcal{S}}
\newcommand{\cT}{\mathcal{T}}

\newcommand{\cV}{\mathcal{V}}

\newcommand{\ccM}{\mathscr{M}}

\newcommand{\gl}{\lambda}

\def\tSch{\mathtt{CS}}

\hypersetup{
citecolor=blue,
colorlinks=true, %colorise les liens
breaklinks=true, %permet le retour a la ligne dans les liens trop longs
urlcolor= blue, %couleur des hyperliens
linkcolor= black, %couleur des liens internes
bookmarksopen=true, %si les signets Acrobat sont crees,
pdftitle={S infinity operator}, %informations apparaissant dans
}

\begin{document}

\title[The $1d$ Stochastic Dirac Operator]{The $1d$ Stochastic Dirac Operator}

\author{Laure Dumaz}
\address{CNRS \& Department of Mathematics and Applications, \'Ecole Normale Sup\'erieure (Paris), 45 rue d’Ulm, 75005 Paris, France}
\email{laure.dumaz@ens.fr}

\author{Cyril Labb\'e}
\address{Université Paris Cité and Sorbonne Université, CNRS, Laboratoire de Probabilités, Statistique et Modélisation, F-75013 Paris, France and Institut Universitaire de France (IUF)}
\email{clabbe@lpsm.paris}

\vspace{2mm}

\date{\today}

\maketitle

\begin{abstract}
	
	We introduce a random Dirac operator on $L^2(\R,\R^2)$ that we call $1d$ Stochastic Dirac Operator and which is indexed by two parameters $\kappa_0,\kappa$. In two companion articles, we show that this operator arises generically as the small disorder limit of one-dimensional random Schr\"odinger operators, thus revealing a universal behavior. In the present paper, we provide a detailed construction of this singular operator, and investigate its spectral properties. In particular, we show that its spectrum is either absolutely continuous when $\kappa = 0$, or pure point when $\kappa>0$. In the latter case, we show that the shape of the eigenfunctions is given by the exponential of a two sided Brownian motion with drift, in connection with a result of Rifkind and Vir\'ag~\cite{RV}. Finally, we establish a delicate convergence result of the sequence of eigenvalues / eigenfunctions from the finite-volume operator to the infinite-volume one.

\medskip

\end{abstract}

\setcounter{tocdepth}{1}
\tableofcontents

\section{Introduction}

The object of this article is to introduce a family of random Dirac operators with white noise matrix potential that we call the $1d$ Stochastic Dirac Operators. 
Those random operators $\cS_{\infty,\kappa_0,\kappa}$, indexed by $\kappa_0,\kappa\in [0,\infty)$, act on functions $f:\R\to\R^2$ in the following (informal) way:
\begin{align}\label{Eq:Sinfty}
	\cS_{\infty,\kappa_0,\kappa} f(t) = J \partial_t f(t) + \dot{\cV}(t) f(t) \,,\quad t\in \R\;,
\end{align}
with
\begin{align}\label{Eq:potential}
	J = \mat{0}{-1}{1}{0}\;,\quad \dot{\cV}(t) = \begin{pmatrix} \kappa_0\dot{B}_0(t)+\kappa\dot{B}_1(t) & \kappa\dot{B}_2(t)\\ \kappa\dot{B}_2(t) &   \kappa_0\dot{B}_0(t)- \kappa\dot{B}_1(t)\end{pmatrix}\,,\quad t\in \R\;,
\end{align}
where $\dot{B}_0$, $\dot{B}_1$ and $\dot{B}_2$ are independent white noises on $\R$.\\

Our main motivation for introducing and studying this family of operators comes from its universality behavior revealed in our companion articles~\cite{DL_UnivCont,DL_UnivDisc}. Therein we prove that $\cS_{\infty,\kappa_0,\kappa}$ arises as the small disorder limit of a large class of $1$-d random Schr\"odinger operators. More precisely, we start from the operator $H := -\Delta + \sigma V$ where $\sigma > 0$ is a deterministic parameter, $\Delta$ is the continuous $1$-d Laplacian and $V$ is a random potential on $\R$. Under mild assumptions on $V$, we show in~\cite{DL_UnivCont} that for any energy $E>0$, the operator $H-E$, after some appropriate transformation, converges as $\sigma \downarrow 0$ to $\cS_{\infty,\kappa_0,\kappa}$, where $\kappa_0,\kappa$ admit explicit expressions in terms of  the covariance of the potential $V$ and of the energy $E$. Let us point out that all possible values $\kappa_0,\kappa \ge 0$ may arise in the limit, according to the law of $V$ one starts from. In~\cite{DL_UnivDisc}, we prove a similar convergence result starting from $1$-d random Schr\"odinger operators on $\Z$. These works are inspired by the study of Kritchevski, Valk\'o and Vir\'ag in~\cite{KVV} of critical discrete one-dimensional random Schr\"odinger operators.\\

Generic one-dimensional Dirac operators take the form
$$ J \partial_t + \sum_{i \in \{0,1,3\}} V_i\sigma_i \;,$$
where the $\sigma_i$ are the Pauli matrices
\begin{align*}
	\sigma_0 := I,\quad \sigma_1 := \mat{1}{0}{0}{-1}, \quad \sigma_3 := \mat{0}{1}{1}{0}.
\end{align*}
In the physics literature~\cite{Thaller} $V_0$ stands for the electric potential, $V_1$ the scalar potential and $V_3$ the anomalous magnetic potential. We refer to~\cite{Weidmann, LevitanSargsjan1991, Zalczer2023} for the construction and general properties of generic Dirac operators.\\
The random Dirac operator obtained by taking $V_1$ a white noise and $V_0=V_3=0$ is intimately related to the Supersymmetric Schr\"odinger Hamiltonian, it was the subject of an important literature in physics~\cite{Bouchaud, Bocquet1999, ComtetTexier, Texier, Texier_2010}.\\
The $1$d Stochastic Dirac Operators are obtained by taking
$$ V_0 = \kappa_0 \dot{B}_0\;,\quad V_1 = \kappa \dot{B}_1 \;,\quad V_3 = \kappa \dot{B}_2\;.$$
It is possible to consider the whole class of random Dirac operators with white noise potentials where each $V_i$ is of the form $\kappa_i \dot{B}_i + \mu_i$. We concentrate on the $1$d Stochastic Dirac Operators because they arise as universal small disorder limits of $1$d random Schr\"odinger operators. From a mathematical point of view, this family is the most tractable to establish exact spectral statistics.

Valk\'o and Vir\'ag in~\cite{VVoperator} introduced random operators, in the framework of canonical systems, whose spectra are given by the limits of (appropriately rescaled) eigenvalues of large random matrices, such as Sine$_\beta$ and Airy$_\beta$. Let us emphasize that these spectra are discrete. The Dirac operators considered here are strongly linked to those introduced in~\cite{VVoperator}: they arise as random conjugations of appropriate canonical systems. These random conjugations bring an additional information on the eigenfunctions. Furthermore, contrary to~\cite{VVoperator}, the spectra of our operators are dense or absolutely continuous.

\subsection*{Overview of the results}

The main contributions of the present article are:\begin{enumerate}[label=(\roman*)]
	\item a rigorous construction of the operators via stochastic flows and canonical systems, along with the convergence of the finite-volume to infinite-volume operator in the strong resolvent sense, see Proposition \ref{Prop:Construction}.
	\item the identification of the a.s.~spectrum as being the whole line $\R$ and of the density of states as being the constant $1/\pi$, whatever $\kappa_0,\kappa$ are, see Proposition~\ref{Th:IDS}
	\item the determination of the nature of the spectrum. If $\kappa>0$, we establish an Anderson localization-type result: the spectrum is pure point with exponentially localized eigenfunctions; while if $\kappa=0$, we show that the spectrum is purely absolutely continuous, see Theorem \ref{Th:Localization} and Proposition \ref{Th:Spectrumk=0}.
	\item the convergence, in the localized regime $\kappa > 0$, of an extension of the spectral measure of the finite-volume operator to the infinite-volume one, see Theorem \ref{Th:CVrandommeas}.
	\item the convergence, in the localized regime $\kappa > 0$, of the sequence of eigenvalues / eigenfunctions of the finite-volume operator to the infinite-volume one, see Corollary \ref{Cor:CVeigenfunctions}.
	\item an explicit expression, in the localized regime $\kappa > 0$, of the intensity measure of the point process of eigenvalues / eigenfunctions of the infinite-volume operator, see Theorem \ref{Th:Intensity}.
\end{enumerate}

Let us make a few comments on these results. The defining expression~\eqref{Eq:Sinfty} is formal and the rigorous construction (i) of the operator is not straightforward. Indeed, white noise is distribution-valued and one needs to lift an indeterminacy on the meaning of the product $\dot{\cV}(t) f(t)$ in~\eqref{Eq:Sinfty}. The operator that we consider arises if one understands the defining expression~\eqref{Eq:Sinfty} in the Stratonovich sense, or equivalently, if one takes as a defining expression
\begin{align}\label{Eq:SinftyIto}
	J \partial_t f(t) + \dot{\cV}(t) f(t) + \Big(\frac{\kappa_0^2}{2} - \kappa^2\Big) J f(t)\,,\quad t\in \R\;,
\end{align}
understood in the It\^o sense. The additional term is nothing but the It\^o-Stratonovich correction term, see Lemma~\ref{Lemma:Strato} for more details. In the sequel, we will always deal with the Stratonovich formulation of the operator.\\
The lack of regularity of the potential induces the following unusual property: the domains of our operators do \emph{not} contain any non-zero smooth function.

Remarkably the parameter $\kappa_0$ has (essentially) no influence on the spectral properties of the operators. On the other hand, our result (iii) reveals a localization dichotomy according to whether $\kappa$ is positive or null.

In (iv) we introduce a random point measure whose atoms are the eigenpairs of the (finite-volume or infinite-volume) operator, which appears as a natural extension of the spectral measure. Then, we prove its convergence from finite to infinite-volume. This is the main achievement of this article, and the last two results rest on it.

The finite to infinite-volume convergence (v), in the localized regime, of the point process of eigenvalues / eigenfunctions is delicate. Note that the countable collection of the eigenvalues of the infinite-volume operator forms a (random) dense subset of $\R$: there is no natural ordering of them and their matching with their finite-volume counterparts is far from being straightforward.

The expression (vi) shows that the eigenfunctions follow a shape that first arose in the work~\cite{RV} of Rifkind and Vir\'ag on critical discrete one-dimensional Schr\"odinger operators. Therein they conjectured its universality among critical one-dimensional operators in the Poisson limit, this was confirmed in the white noise case~\cite{DL_Critical}.

\section{Main results}

\subsection*{Construction of the operators}

We will define the operators $\cS_{\infty,\kappa_0,\kappa}$, along with their finite-volume counterparts $\cS_{\tau,\kappa_0,\kappa}$, indexed by $\tau > 0$ and $\kappa_0,\kappa\ge 0$. The latter take the same defining expression~\eqref{Eq:Sinfty} except that they act on functions $f:(-\tau/2,\tau/2)\to\R^2$ satisfying some prescribed separated boundary conditions: $f(\pm \tau/2)$ is parallel to $b_\pm \in \R^2 \backslash \{0\}$.

Our main tool, and one of our original contributions in this article, is a stochastic flow of $2\times 2$ complex random matrices $Y_z(s,t)$, $z\in \C$, $s,t \in \R$ which solve the Stratonovich SDE
\begin{align}\label{Eq:MatrixSDE}
	J dY_z(s,t) + \circ d\cV(t) Y_z(s,t) = z Y_z(s,t) dt\;,\quad t \in [s,\infty)\;,
%	J dY_z(s,t) + d\cV(t) Y_z(s,t) + (\frac{\kappa_0^2}{2} - \kappa^2) J Y_z(s,t) dt = z Y_z(s,t) dt\;,\quad t \in [s,\infty)\;,
\end{align}
starting from $Y_z(s,s) = I$ (and similarly when $t \le s$). Here
$$ \cV(t) := \begin{pmatrix} \kappa_0B_0(t) + \kappa B_1(t) & \kappa B_2(t)\\ \kappa B_2(t) & \kappa_0 B_0(t) - \kappa B_1(t)\end{pmatrix}\;,$$
where the $B_i$'s are two-sided Brownian motions obtained by integrating the $\dot{B}_i$'s.
By stochastic flow, we notably mean that we deal with a collection of solutions of the above SDEs which is \emph{continuous} in $s,t,z$ and satisfies the \emph{cocycle property} $$Y_z(t,u)Y_z(s,t) = Y_z(s,u)$$ for all $s,t,u$ and $z$, with probability one. We refer to Subsection~\ref{Subsec:StochasticFlow} for more details.\\

Our construction of $\cS_{\tau,\kappa_0,\kappa}$ and $\cS_{\infty,\kappa_0,\kappa}$ as self-adjoint operators on $L^2((-\tau/2,\tau/2),\R^2)$ and $L^2(\R,\R^2)$ respectively, will rest on a mapping between Dirac operators and canonical systems, that relies crucially on the processes $Y_z$ and which will be carefully presented in Section \ref{Sec:Construction}.\\

While the mapping between Dirac operators and canonical systems was used previously in the literature, see~\cite{Remling,VVoperator,DL_Critical}, here we apply this mapping \emph{simultaneously} for all $\tau \in (0,+\infty]$, that is for an uncountable collection of operators and this requires to construct a \emph{coherent} collection of random variables $Y_z(s,t)$, hence the introduction of the stochastic flow. This framework allows to construct on a single probability space all finite-volume operators (for all b.c.~at once) along with the infinite-volume operator.\\

Although our construction relies on the aforementioned mapping to canonical systems, the formal expression~\eqref{Eq:Sinfty} finds a justification a posteriori: the spectra of $\cS_{\tau,\kappa_0,\kappa}$ and $\cS_{\infty,\kappa_0,\kappa}$ are given by solutions of the stochastic flow that either satisfy the chosen boundary conditions (in the finite volume case) or satisfy some integrability conditions (in the infinite volume case). We refer to Section \ref{Sec:Construction} for precise statements.\\

The following result relates the finite and infinite-volume operators:

\begin{proposition}[Finite to infinite-volume operator convergence]\label{Prop:Construction}
	Fix $\kappa_0,\kappa\ge 0$. Almost surely, for any choice of boundary conditions $b_\pm \in \R^2 \backslash \{0\}$ the self-adjoint operator $\cS_{\tau,\kappa_0,\kappa}$ converges as $\tau\to\infty$ in the strong resolvent sense towards $\cS_{\infty,\kappa_0,\kappa}$.
\end{proposition}

We refer to Remark \ref{Rem:highenergy} for details on the relation with our previous work~\cite{DL_Critical}.

\subsection*{Symmetry properties and density of states}

Before addressing the spectral properties of the infinite-volume operator $\cS_{\infty,\kappa_0,\kappa}$, let us mention two symmetry properties it satisfies:\begin{enumerate}
	\item Invariance in law under shifts in space, see Proposition \ref{Prop:SymmetrySpatialShift}.
	\item Invariance in law under shifts in energy (up to a unitary transformation), see Proposition \ref{Prop:SymmetrySpectralShift}.
\end{enumerate}

The second property indicates that the density of states should be constant. To put this on a firm ground, we define
\begin{equation}\label{Eq:NtauStates}
	N_{\tau,\kappa_0,\kappa}(a,b) := \#\{\mbox{eigenvalues of }\cS_{\tau,\kappa_0,\kappa}\mbox{ lying in }(a,b]\}\;,\quad a \le b\;.
\end{equation}

\begin{proposition}[Density of states]\label{Th:IDS}
	Fix $\kappa_0,\kappa\ge 0$ and some boundary conditions $b_\pm \in \R^2 \backslash \{0\}$. For all $a<b$, the following convergence holds in $L^1(\Omega)$
	$$ \frac{N_{\tau,\kappa_0,\kappa}(a,b)}{\tau} \to N(a,b) := \frac{b-a}{\pi}\;,\quad \tau \uparrow \infty\;.$$
\end{proposition}
\begin{remark}
	This convergence also holds almost surely. However the proof then requires much more effort and we preferred not to present the details.
\end{remark}

This implies that the density of states, $n(E)$, defined through $N(a,b) = \int_{(a,b]} n(E)dE$, is constant equal to $1/\pi$. This density does not depend on $\kappa_0,\kappa$ and, thus, remains the same in the absence of noise ($\kappa_0=\kappa=0$).

\subsection*{The localized regime}

We turn to the investigation of the nature of the spectrum of $\cS_{\infty,\kappa_0,\kappa}$. It turns out that the parameter $\kappa$ plays a crucial role. We start by investigating the regime where $\kappa>0$.

\begin{theorem}[Exponential localization]\label{Th:Localization}
	Take $\kappa_0\ge 0$ and $\kappa>0$. Almost surely, the spectrum of $\cS_{\infty,\kappa_0,\kappa}$ is equal to $\R$ and is pure point, and each eigenvalue is of multiplicity one. In other words, almost surely $\cS_{\infty,\kappa_0,\kappa}$ admits a dense sequence $(\lambda_{j,\infty})_{j\in \Z}$ of distinct eigenvalues and the associated normalized eigenfunctions $(\Psi_{j,\infty})_{j\in \Z}$ form an orthonormal basis of $L^2(\R,\R^2)$. In addition, each eigenfunction decays exponentially fast at rate $\kappa^2$: almost surely for every $j \in \Z$
	\begin{equation}\label{Eq:ExpoLoc}
		\lim_{t\to \pm\infty} -\frac1{|t|} \ln |\Psi_{j,\infty}(t)| = \kappa^2\;,
	\end{equation}
	where $|\Psi_{j,\infty}(t)|$ denotes the Euclidean norm of the vector $\Psi_{j,\infty}(t) \in \R^2$.
\end{theorem}

We now establish a convergence result of the eigenpairs of the finite-volume operator $\cS_{\tau,\kappa_0,\kappa}$ endowed with fixed, but arbitrary, boundary conditions $b_\pm$, to those of the infinite-volume operator $\cS_{\infty,\kappa_0,\kappa}$. Since the limiting eigenvalues are a.s.~dense in $\R$, there is no canonical ordering of the limiting eigenpairs and this makes the convergence particularly delicate.

We first introduce convenient function spaces to which our eigenfunctions belong. We let $\mathcal{F}$ be the separable Banach space of all continuous functions $f:\R\to\R^2$ satisfying
\begin{equation}\label{Eq:DefF}
	\| f\|_{\mathcal{F}} := \sup_{t\in\R} | f(t)| e^{\frac{\kappa^2 |t|}{2}} < \infty\;,\quad \mbox{ and }\lim_{t\to\pm\infty} f(t) e^{\frac{\kappa^2 |t|}{2}}=0\;.
\end{equation}
By Theorem \ref{Th:Localization}, a.s.~each eigenfunction $\Psi_{j,\infty}$ belongs to $\cF$ (note that the exponential rate $\kappa^2/2$ is arbitrary: any value in $(0,\kappa^2)$ would do, as it would fit with~\eqref{Eq:ExpoLoc}). The finite-volume eigenfunctions are only defined on $[-\tau/2,\tau/2]$; we (arbitrarily) extend them into continuous functions on the whole line $\R$ by setting their values to zero outside $[-\tau/2-1,\tau/2+1]$ and by linear interpolation in between. These extensions belong to $\cF$.\\

The normalized eigenfunctions are unique up to a change of sign. We could decide of a convention to select one of the two versions of each eigenfunction: however, this would cause issues in scaling limits later on. Instead, we prefer to work with the equivalence classes associated to the relation
$$ f \sim g \Longleftrightarrow f = \pm g\;.$$
We let $\tilde{\cF}$ be the quotient of $\cF$ by this equivalence relation, endowed with the metric
$$ \tilde{d}(\tilde{f}, \tilde{g}) = \min( \| f-g \|_\cF, \|f+g\|_\cF)\;.$$
It is straightforward to check that $\tilde{\cF}$ is a Polish space. From now on, we view the eigenfunctions as equivalence classes in $\tilde{\cF}$.\\

To establish the aforementioned convergence, our main tools are the following (weighted) point measures:
$$ \mu_{\tau,\kappa_0,\kappa} := \sum_{j\in\Z}  |\Psi_{j,\tau}(0)|^2 \delta_{(\lambda_{j,\tau},\Psi_{j,\tau})}\;,\quad \mu_{\infty,\kappa_0,\kappa} := \sum_{j\in\Z}  |\Psi_{j,\infty}(0)|^2 \delta_{(\lambda_{j,\infty},\Psi_{j,\infty})}\;.$$
These measures are viewed as elements of the space $\ccM$ of all non-negative measures $m$ on $\R\times \tilde{\cF}$ that satisfy $m([a,b]\times \tilde{\cF}) < \infty$ for all $a<b$, endowed with the topology of vague convergence for the first coordinate, and weak convergence for the second coordinate.

The weights in front of the Dirac masses might seem unnatural at first sight. However, if we removed them then the infinite-volume measure would put infinite mass on $[a,b]\times \tilde{\cF}$ for any $a<b$ a.s. (recall that the spectrum is dense a.s.). The weights somehow ``select'' the eigenvalues whose eigenfunctions are localized close to $t=0$. Since the operator is invariant in law under spatial shifts, the choice $t=0$ is arbitrary and has no influence on the law of the measure. We are then able to prove the following convergence result.

\begin{theorem}[Convergence of the weighted eigenpair measures]\label{Th:CVrandommeas}
	Take $\kappa_0\ge 0$ and $\kappa >0$ and fix some boundary conditions $b_\pm \in \R^2 \backslash \{0\}$. The measure $\mu_{\tau,\kappa_0,\kappa} $ converges to $\mu_{\infty,\kappa_0,\kappa} $ in probability in $\ccM$ as $\tau\to\infty$. Moreover, we have
	\begin{align*}
		\E\Big[ \Big|\int F d\mu_{\tau,\kappa_0,\kappa} - \int F d\mu_{\infty,\kappa_0,\kappa} \Big|^p\Big] \to_{\tau \to \infty} 0
	\end{align*}
	for any $p \geq 1$ and any $F \,:\, \R \times \tilde{\cF} \to \R$ continuous, bounded and compactly supported in the first coordinate.
\end{theorem}

\medskip

From this result, we are able to deduce a very strong convergence of the eigenvalues / eigenfunctions, provided one enumerates them appropriately (recall that no natural ordering exist on the eigenpairs since the limiting spectrum is dense). We order the eigenpairs $(\lambda_{j,\infty},\Psi_{j,\infty})_j$ in the decreasing order of their weights $w(\lambda_{j,\infty},\Psi_{j,\infty})$ where
\begin{align*}
	w(\gl,\Psi) := \frac{|\Psi(0)|^2}{1+ |\gl|^2}\;,
\end{align*}
we refer to Subsection \ref{Subsec:CVeigenfunction} for more details.

\begin{corollary}[From finite to infinite-volume eigenfunctions]\label{Cor:CVeigenfunctions}
	Take $\kappa_0\ge 0$ and $\kappa >0$ and fix some boundary conditions $b_\pm \in \R^2 \backslash \{0\}$. One can enumerate the eigenvalues and eigenfunctions of $\cS_{\tau,\kappa_0,\kappa}$ and $\cS_{\infty,\kappa_0,\kappa}$
	$$ (\lambda_{m,\tau},\Psi_{m,\tau})_{m\ge 1}\;,\quad (\lambda_{m,\infty},\Psi_{m,\infty})_{m\ge 1}\;,$$
	in such a way that for any fixed $M \geq 1$, $(\lambda_{m,\tau},\Psi_{m,\tau})_{1 \leq m\leq M}$ converges in probability to $(\lambda_{m,\infty},\Psi_{m,\infty})_{1 \leq m\leq M}$ as $\tau\to\infty$.
\end{corollary}

Let us mention that the ordering of the finite-volume eigenpairs depends on the limiting eigenpairs, we refer to the proof of the corollary for more details. Although the convergences in the last two results look natural, we are not aware of any similar result in the literature on Schr\"odinger nor Dirac operators. Indeed, the finite to infinite-volume convergence is usually done at the level of locally unfolded eigenvalues (keeping their spacing of order $1$) where one expects to have a convergence towards a Poisson point process for random Schr\"odinger operators (in any dimension) in the localized regime, see \cite{KillipNakano,GerminetKlopp2014} for precise assumptions. Such a local limit at the level of the point process of eigenvalues for our operators follows from a straightforward adaptation of the arguments of the papers~\cite{Nakano2015, AllezDumaz}.

\smallskip

We now compute explicitly the intensity measure of $\mu_{\infty,\kappa_0,\kappa}$, that is, the expectation of $\mu_{\infty,\kappa_0,\kappa}$. We will see that we obtain a remarkably simple expression thanks to the special noise structure of the operators.
To that end, we need to introduce a random variable $(S,\gamma^{(\lambda)})$ whose law depends on a parameter $\lambda \in\R$. Given two independent two-sided Brownian motions $W_1$ and $W_2$ and an independent uniform random variable $\gamma_0$ on $[0,2\pi]$, we set
$$ S(t) :=  \exp\Big(\kappa W_1(t) - \kappa^2|t|\Big)\;,\quad t\in\R\;,$$
and
$$ \gamma^{(\lambda)}(t) := \gamma_0 + \sqrt{\kappa_0^2+\kappa^2} \,W_2(t) + \lambda t\;,\quad t\in\R\;.$$
Note that $S$ only depends on $W_1$, while $\gamma^{(\lambda)}$ only depends on $W_2$ and $\gamma_0$.

\begin{theorem}[Intensity of the eigenvalue/eigenfunction pairs]\label{Th:Intensity}
	Take $\kappa_0\ge 0$ and $\kappa >0$. For any non-negative measurable map $G:\R\times \tilde{\cF} \to \R$, we have
	\begin{equation}\label{Eq:Intensity}
		\E\Big[\sum_{j\in \Z} | \Psi_{j,\infty}(0)|^2 G\Big(\lambda_{j,\infty},\Psi_{j,\infty}\Big)\Big]= \int \E\bigg[G\bigg(\gl, \frac{S(\cdot)}{\| S\|_{L^2(\R)}} \begin{pmatrix} \sin \gamma^{(\lambda)} (\cdot) \\\cos \gamma^{(\lambda)} (\cdot) \end{pmatrix}\bigg)\bigg] \frac{d\gl}{\pi}\;.
	\end{equation}
	In particular, for any non-negative measurable map $h:\R\to\R$ we have
	\begin{equation}\label{Eq:IntensityLambda}
		\E\Big[\sum_{j\in \Z} | \Psi_{j,\infty}(0)|^2 h(\lambda_{j,\infty})\Big]= \int h(\lambda) n(\lambda)d\lambda\;,
	\end{equation}
	where $n(\gl) = (1/\pi)$ is the density of states defined above.
\end{theorem}
This identifies the law of a typical eigenfunction picked from the spectral measure: its shape is given by the normalized exponential two-sided Brownian motion with negative drift $-\kappa^2$ on both sides, while its phase is an independent Brownian motion with linear drift $\gl t$.
A similar result in finite volume is stated in Proposition \ref{Prop:GMP2}, and actually covers the whole range $\kappa_0,\kappa\ge 0$.

As a corollary, we obtain an explicit formula for the density-density correlator, a quantity introduced by Anderson in~\cite{Anderson58}. In the physics literature, a non-vanishing infinite-volume limit of the density-density correlator is used as a criterion for localization.
\begin{corollary}[Density-density correlator]\label{Cor:correlator}
	Assume that $\kappa>0$. For $E,t,t'\in\R$, and for all Borel set $A \subset \R$, define the measure
\begin{align*}
	m_{E,t,t'}(A)
	:=
	\E\Big[
	\sum_{j \in \Z}
	|\Psi_{j,\infty}(t)|^2
	|\Psi_{j,\infty}(t')|^2\un_{E-\lambda_{j,\infty} \in A}
	\Big].
\end{align*}
	Then
\begin{align*}
	m_{E,t,t'}(A)
	=
	\frac1\pi 	\E\bigg[
	\frac{S(t'-t)^2}{\|S\|_{L^2(\R)}^2}
	\bigg] \, |A|\,.
\end{align*}
where $|A|$ is the Lebesgue measure of $A$.
\end{corollary}
Note that the density-density correlator is independent of the energy $E$ and depends on $t, t'$ only through $t'-t$, which was expected due to the spectral and spatial stationarity of the operator. This explicit expression should match with the formula of Kolokolov~\cite{Kolokolov1993} obtained in terms of special functions, see Conjecture \ref{Conj:Kolokolov} for more details.

\subsection*{The delocalized regime}

Let us now cover the regime where $\kappa =0$. Although the density of states is the same as in the regime $\kappa >0$, the nature of the spectrum is drastically different.

\begin{proposition}\label{Th:Spectrumk=0}
	Take $\kappa_0\ge 0$ and $\kappa=0$. Almost surely, the spectrum of $\cS_{\infty,\kappa_0,0}$ is equal to $\R$ and is purely absolutely continuous.
\end{proposition}

The spectral measures admit explicit expressions, see Subsection \ref{Subsec:k=0}. It can be checked that the random measure $\mu_{\tau,\kappa_0,0}$ converges in probability as $\tau\uparrow\infty$ to the trivial limit:
$$\mu_{\infty,\kappa_0,0}(d\lambda,d\Psi) := (1/\pi) d\lambda \otimes \delta_0(\Psi)\;.$$

\subsection*{Dependence on $\kappa_0$ and $\kappa$}

The noise $\kappa_0\dot{B}_0$ has played almost no role so far: if one sets $\kappa_0$ to zero, then only the variance of the process $\gamma^{(\lambda)}$ is modified, while the density of states, the exponential rate of decay of the eigenfunctions and the shape process $S$ are left unchanged. On the other hand, our results reveal a dichotomy, in terms of spectral property, between $\kappa > 0$ and $\kappa=0$.
In Subsection \ref{Subsec:kappa}, we will see that:\begin{enumerate}
	\item (Removing $\kappa_0 \dot{B}_0$) There exists an explicit, random unitary transform that maps $\cS_{\infty,\kappa_0,\kappa}$ to $\cS_{\infty,0,\kappa}$,
	\item (Normalizing $\kappa$) Provided $\kappa > 0$, there exists an explicit, deterministic unitary transform that maps $\cS_{\infty,\kappa_0,\kappa}$ to $\kappa^2 \cS_{\infty,\kappa_0/\kappa,1}$.
\end{enumerate}
Put together, these results entail that there are only two fundamentally distinct operators among this family: $\cS_{\infty,0,0}$ and $\cS_{\infty,0,1}$. We could have presented all our results only for these two particular operators: for completeness and readability, we preferred to keep the dependence on $\kappa_0$ and $\kappa$ explicit everywhere.

\subsection*{Organization of the paper}

The rest of the paper is organized as follows. Section \ref{Sec:Construction} gives the rigorous construction of the finite and infinite-volume stochastic Dirac operators through stochastic flows and canonical systems, and proves the strong resolvent convergence of the finite volume operator towards the infinite-volume one. Section \ref{Sec:DOSGMP} proves the density of states formula, and establishes a GMP formula for the finite-volume operator. Section \ref{Sec:Spectral} studies the spectral representation and the spectral type of the infinite-volume operator, proving localization for $\kappa > 0$ and absolute continuity for $\kappa = 0$. Section \ref{Sec:CVeigenfunctions} proves the convergence of the weighted eigenvalue/eigenfunction point measure, derives the convergence of the pairs of eigenvalues/eigenfunctions under a suitable enumeration, and computes the intensity measure and density-density correlator. Appendix \ref{Appendix} collects technical results used in the previous sections.

\subsection*{Notations}
Given a matrix $M$, we denote by $M^\intercal$ its transpose. We introduce the basis vectors
\begin{equation}\label{Eq:e1e2}
	e_1 = \begin{pmatrix} 1 \\ 0 \end{pmatrix}\;,\quad e_2 = \begin{pmatrix} 0 \\ 1 \end{pmatrix}\;,
\end{equation}
and we define the matrix
$$ J = \begin{pmatrix} 0 & -1 \\ 1 & 0 \end{pmatrix}\;.$$
We observe that $J^{-1} = J^\intercal = - J$ and for any invertible matrix $A$
\begin{equation}\label{Eq:JMJ}
	JAJ^{-1} = \det(A) (A^\intercal)^{-1}\;.
\end{equation}
Given two vectors $u,v \in \C^2$, we define their Wronskian as
\begin{equation}\label{Eq:Wronskien}
	W(u,v) = u_1v_2 - u_2v_1 = - u^\intercal J v\;.
\end{equation}
We also write $u\parallel v$ if there exists $q\in \C$ such that $u = q v$.\\
Given a locally integrable function $R$ valued in the set of $2\times 2$ real symmetric positive definite matrices and some interval $E \subset \R$, we let $L^2_R(E)$ be the (quotient) space of all measurable functions $f: E\to \R^2$ such that
$$ \int_E f(s)^\intercal R(s) f(s) ds < \infty\;.$$
As any $L^2$ space, this is a Hilbert space for the inner product $(f,g) \mapsto \int_E f(s)^\intercal R(s) g(s) ds$.\\
\textbf{From now on, we will omit the subscripts $\kappa_0,\kappa$ in $\cS_{\infty,\kappa_0,\kappa}$, $\cS_{\tau,\kappa_0,\kappa}$, $\mu_{\infty,\kappa_0,\kappa}$, $\mu_{\tau,\kappa_0,\kappa}$ in order to alleviate the notation.}

\subsection*{Acknowledgements} Our work was partially supported by the ANR project RANDOP ANR-24-CE40-3377. C.L.~also acknowledges support from the Institut Universitaire de France and from the ANR project Smooth ANR-22-CE40-0017. We would like to thank Yan Fyodorov for suggesting to investigate the density-density correlator in Corollary \ref{Cor:correlator} and for pointing out the article by Kolokolov~\cite{Kolokolov1993}.

\subsection*{AI statement} After the completion of this work, we used ChatGPT Pro and Claude Opus 5.5 to proofread the article. Besides typographic errors and bibliographic suggestions, they identified an issue with our original convention on the sign of the eigenfunctions: this led us to introduce the quotient space $\tilde{\cF}$.

\section{Construction of the operator}\label{Sec:Construction}

In this section, except when explicitly mentioned, we will work implicitly with generic values of $\kappa_0,\kappa \ge 0$.

\subsection{The SDEs}

As mentioned  in the introduction, the rigorous construction of the operators $\cS_\tau$ and $\cS_\infty$ relies on the family of SDE associated with the formal eigenproblem $\cS_\infty f = z f$. We thus need to study the collection of solutions $y_z(t) \in\C^2, t \ge 0$ to the family of Stratonovich SDE parametrized by $z\in \C$
\begin{align}\label{Eq:FwdSDE}
	   J dy_z(t) + \circ d\cV(t) y_z(t) = z y_z(t) dt\;,
%	J dy_z(t) + d\cV(t) y_z(t) +(\frac{\kappa_0^2}{2}-\kappa^2) J y_z(t) dt = z y_z(t) dt\;,
\end{align}
starting from $y_z(0) = x \in \C^2$, along with translated versions thereof where the initial condition is taken at some $t_0\in\R$ rather than $0$. Here we recall that
$$ \cV(t) := \begin{pmatrix} \kappa_0B_0(t) + \kappa B_1(t) & \kappa B_2(t)\\ \kappa B_2(t) &\kappa_0B_0(t)-\kappa B_1(t)\end{pmatrix}\;,\quad t\in\R\;.$$

By~\eqref{Eq:FwdSDE}, we mean
\begin{equation*}
	\left\{
	\begin{array}{rcl}
	-y_z^{(2)}(t) + y_z^{(2)}(0) + \int_0^t \sum_{i=1}^2 y_z^{(i)}(s) \circ d\cV^{(1,i)}(s) &= & z \int_0^t y_z^{(1)}(s) ds \;,\\[0.8em]
		y_z^{(1)}(t) - y_z^{(1)}(0) + \int_0^t \sum_{i=1}^2 y_z^{(i)}(s) \circ d\cV^{(2,i)}(s)&=& z \int_0^t y_z^{(2)}(s) ds \;,
	\end{array}
	\right.
\end{equation*}
where $\circ$ stands for Stratonovich integral.

\smallskip

The next result identifies the It\^o-Stratonovich correction for this SDE:
\begin{lemma}\label{Lemma:Strato}
	The It\^o-Stratonovich correction for the SDE \eqref{Eq:FwdSDE} is equal to $(\frac{\kappa_0^2}{2}-\kappa^2)J y_z(t)$. If we rewrite the equation in the It\^o sense then
$$	J dy_z(t) + d\cV(t)  y_z(t) + \Big(\frac{\kappa_0^2}{2}-\kappa^2\Big) J y_z(t) dt = z y_z(t) dt\;,\quad t \ge 0\;.$$
\end{lemma}
\begin{proof}
	It is well known that, given a real-valued Brownian motion $B$ and a semimartingale $X$ taking values in $\C$, the It\^o-Stratonovich correction is given by
	$$ \int_0^t X(s) \circ dB(s) = \int_0^t X(s) dB(s) + \frac12  \langle X,B\rangle_t\;.$$
	Now if $X$ takes values in $\C^2$ and satisfies $dX(t) = \sum_{i=0}^2 A_i(t)\circ dB_i(t) + C(t) dt$ for some $\C^2$-valued processes $A_i$, $C$ and some real-valued Brownian motions $B_i$, then for any $i\in\{0,1,2\}$
	$$ \int_0^t X(s) \circ dB_i(s) = \int_0^t X(s) dB_i(s) + \frac12  \int_0^t A_i(s) ds\;.$$
Recall the Pauli matrices from the introduction. We obtain
\begin{align*}
	J dy_z &= -\kappa_0 y_z \circ dB_0-\kappa \sigma_1 y_z \circ dB_1 -\kappa \sigma_3 y_z \circ dB_2 + z y_z dt \\
	&=-\kappa_0 y_z dB_0 - \frac{\kappa_0^2}{2} J y_z dt -\kappa \sigma_1 y_z  dB_1 - \frac{\kappa^2}{2} \sigma_1 J \sigma_1 y_z dt \\
	&\,\quad-\kappa \sigma_3 y_z  dB_2 - \frac{\kappa^2}{2} \sigma_3 J \sigma_3 y_z dt + z y_z dt \\
	&=-\kappa_0 y_z dB_0  -\kappa \sigma_1 y_z  dB_1 -\kappa \sigma_3 y_z  dB_2 -( \frac{\kappa_0^2}{2} - \kappa^2) J y_z dt  + z y_z dt \,,
\end{align*}
which gives the result.
\end{proof}

Second, when $z=\lambda \in \R$, it is convenient to write $y_\lambda$ in polar coordinates
$$ y_\lambda(t) = e^{\rho_\lambda(t)} \begin{pmatrix} \sin \theta_\lambda(t) \\ \cos \theta_\lambda(t) \end{pmatrix}\;,\quad t\ge 0\;,$$
where $(\theta_\lambda, \rho_\lambda)$ solve the Stratonovich SDE
\begin{equation}\label{Eq:FwdPolarSDE}\begin{split}
	d\theta_\lambda(t) &= \lambda dt - \kappa_0 \circ dB_0(t) + \kappa\,\Re(e^{2i \theta_\lambda(t)} \circ d\cB(t))\;,\quad t\ge 0\;,\\
	d\rho_\lambda(t) &= \kappa\,\Im(e^{2i \theta_\lambda(t)} \circ d\cB(t))\;,\quad t\ge 0\;,
\end{split}
\end{equation}
for some initial condition $(\theta_\lambda(0) , \rho_\lambda(0))$ compatible with $y_\lambda(0) = x$ and where
\begin{equation}\label{Eq:ComplexBM}
	\dot{\cB}(t) := \dot{B}_1(t) + i \dot{B}_2(t)\;,\quad t\in \R\;,
\end{equation}
is a (non-normalized) complex white noise. These SDE can be written in the It\^o sense (note that it does not change the equation for the phase):
\begin{equation}\label{Eq:FwdPolarSDEIto}\begin{split}
		d\theta_\lambda(t) &= \lambda dt - \kappa_0 dB_0(t) + \kappa\,\Re(e^{2i \theta_\lambda(t)} d\cB(t))\;,\quad t\ge 0\;,\\
		d\rho_\lambda(t) &= \kappa^2 dt + \kappa\,\Im(e^{2i \theta_\lambda(t)} d\cB(t))\;,\quad t\ge 0\;.
	\end{split}
\end{equation}

Let us point out that adding $a$ to $\rho_\lambda(0)$ has the trivial effect of adding $a$ to $\rho_\lambda(t)$ as well: this property is a direct consequence of the fact that \eqref{Eq:FwdSDE} is linear in $y_\lambda$.

\subsection{Stochastic flows}\label{Subsec:StochasticFlow}

It is a priori unclear how to associate a genuine self-adjoint operator to the expression \eqref{Eq:Sinfty}. If we restrict ourselves to constructing the operator on a segment $(-\tau/2,\tau/2)$ given some b.c.~at the two endpoints, our starting point will be to require that the set of eigenvalues / eigenfunctions is given by those $(\lambda,y_\lambda)$ that are such that: (a) $y_\lambda(-\tau/2)$ satisfies the left b.c., (b) $y_\lambda$ satisfies \eqref{Eq:FwdSDE} on the segment $(-\tau/2,\tau/2)$ instead of $[0,\infty)$, and (c) $y_\lambda(\tau/2)$ satisfies the right b.c.\\

To implement this idea simultaneously for all $\tau > 0$, one immediately faces the task of defining \emph{all} solutions $y_\lambda$, $\lambda \in \R$, starting from all possible initial conditions at any time. Stochastic flows appear to be the right object to fulfill this task: the next result introduces a stochastic flow associated to the SDEs \eqref{Eq:FwdSDE} and collects some important properties for the sequel.

\begin{proposition}\label{Prop:Flow}
	There exists a collection of random variables $Y_z(s,t)$ indexed by $z\in \C$, $(s,t)\in \R^2$, taking values in the set of $2\times 2$ complex matrices and satisfying the following properties:\begin{enumerate}
		\item (Continuity) The following map is continuous :
			\begin{align*}
			\C\times \R^2 \times \C^2 &\to \C^2\\
			(z,(s,t),x)&\mapsto Y_z(s,t)x\;,
		\end{align*}
		\item (Cocycle) For any $s,t,u \in \R$ and any $z\in \C$
		$$Y_z(t,u) Y_z(s,t) = Y_z(s,u)\;,\quad\mbox{ and }\quad Y_z(s,s) = I\;.$$
		\item (SDE) For any given $z\in\C,s\in\R,x\in\C^2$, if we set $y(t) := Y_z(s,t)x$ then almost surely $y$ solves
		\begin{align}\label{Eq:FwdSDEbis}
			   J dy(t) + \circ d\cV(t) y(t) = z y(t)dt\;,\quad t \ge s\;,
%			J dy(t) + d\cV(t) y(t) + (\frac{\kappa_0^2}{2}-\kappa^2) J y(t)dt = z y(t)dt\;,\quad t \ge s\;,
		\end{align}
		while if we set $y(t) := Y_z(-s,-t)x$ and $\hat{\cV}(t) := \cV(-t)$, then almost surely
		\begin{align}\label{Eq:BwdSDEbis}
			J dy(t) + \circ d\hat{\cV}(t) y(t) = -z y(t)dt\;,\quad t \ge s\;.
%			J dy(t) + d\hat{\cV}(t) y(t) +(\frac{\kappa_0^2}{2}-\kappa^2) J y(t)dt = -z y(t)dt\;,\quad t \ge s\;.
		\end{align}
	\end{enumerate}
	Moreover, the flow is unique up to a modification on a set of probability $0$.
\end{proposition}
The first two properties are satisfied for all $\omega \in\Omega$, this is one of the crucial aspects of stochastic flows. On the other hand, since the uniqueness of solutions to SDE is only true up to $\P$-null set, the third property can only be satisfied almost surely.\\
The cocycle property implicitly asserts that for any $s<t$, $Y_z(s,t)$ is invertible with inverse $Y_z(t,s)$.\\
\begin{proof}[Proof of Proposition \ref{Prop:Flow}]
	For any $s\in \R$, $z\in\C$ and $x\in\C^2$, we let $(Y_z(s,t)x,t\ge s)$ and $(Y_z(s,t)x,t\le s)$ be the a.s.~unique strong solutions of the forward and backward SDE presented in the statement. By~\cite[Lemma 4.5.6]{Kunita}, for all $p\ge 2$ and all $R>0$, there exists $C_p(R) > 0$ such that for all $z \in \C$ such that $|z|\le R$, all $x,x' \in \C^2$ such that $|x|,|x'| \le R$ and all $-R\le s \le t\le R$, $-R\le s' \le t' \le R$
	\begin{equation}\label{Eq:MomentSDE}
		\E[|Y_z(s,t)x - Y_{z}(s',t')x'|^p] \le  C_p(R)\big(|s-s'|^{\frac{p}{2}} + |t-t'|^{\frac{p}{2}} + |x-x'|^p\big)\;.
	\end{equation}
	Furthermore a direct estimation based on Gr\"onwall's Lemma shows that, up to increasing the constant $C_p(R)$, for all $z,z'\in\C$ such that $|z|,|z'| \le R$, all $x\in\C^2$ such that $|x|\le R$ and all $-R\le s\le t \le R$
	\begin{equation}\label{Eq:MomentSDE2}
		\E[|Y_z(s,t)x - Y_{z'}(s,t)x|^p] \le  C_p(R)|z-z'|^p\;.
	\end{equation}
	The same bounds hold for the backward SDE. By Kolmogorov Continuity Theorem~\cite[Theorem I.2.1]{RevuzYor}, there exists a modification, still denoted $Y$, such that $(z,(s,t),x)\mapsto Y_z(s,t)x$ is continuous. Since $t\mapsto Y_z(s,t)x$ was already continuous, the modification still satisfies the forward and backward SDE a.s., namely (3).\\
	For any fixed $z\in\C$, the arguments presented in~\cite[Chapter 4]{Kunita} ensure the existence of a continuous modification of $Y_z$ such that (2) and (3) are satisfied. Since the process was already continuous, $Y_z$ and its modification coincide on some event $\Omega_z$ of probability one.\\
	We now work on $\cap_{z\in D} \Omega_z$, where the intersection is taken over a countable dense subset $D$ of $\C$.  On this event, (1) is satisfied, and (2) is satisfied for all $z \in D$, and by density of $D$ in $\C$, and the continuity spelled out in (1), we deduce that (2) is satisfied for all $z\in\C$.\\
	There is one property that needs to be justified: the linearity in $x$ of the solution, which is implicit in the requirement that $Y_z(s,t)$ takes values in the set of $2\times 2$ complex matrices. Let us provide the arguments. It is immediate to check at the level of the SDE that, for fixed $z,s,t,x,y,a,b$, the identity $Y_z(s,t)(ax + by) = aY_z(s,t)x + bY_z(s,t)y$ holds almost surely. Thus, this property remains true simultaneously for a dense countable collection of such parameters, and one can conclude by the continuity property (1).\\
	Finally, uniqueness follows from pathwise uniqueness for the forward/backward SDEs. Indeed, if $\widetilde Y$ is another flow satisfying the same properties and driven by the same Brownian motions, for every fixed $z \in \C$, $s \in \R$, 
	\begin{align*}
		\P( \forall t \in \R, \; \forall i \in \{1,2\},\; \widetilde Y_z(s,t)e_i = Y_z(s,t) e_i) = 1\,,
	\end{align*}
where $e_i$ are the two canonical basis vectors of $\C^2$. Taking the parameters $(z,s)$ in dense sets and using continuity, we obtain
	\begin{align*}
			\P( \forall t \in \R, \; \forall z \in \C,\; \forall s \in \R, \; \forall x \in \C^2,\; \widetilde Y_z(s,t)x = Y_z(s,t)x) = 1\,,
	\end{align*}
	which proves uniqueness up to a null set.
\end{proof}

We refer to Subsection \ref{Sec:Polar} for a similar result on the polar coordinates. Let us conclude with a useful observation on the solutions of the stochastic flow.

\begin{lemma}\label{Lemma:Wronskian}
	Almost surely, for any $z\in \C$ and any $s,t\in \R$, we have $Y_z(s,t)^\intercal J Y_z(s,t) = J$. As a consequence, almost surely for any $x,y \in \C^2$ and any $s\in\R$, the Wronskian~\eqref{Eq:Wronskien} of the vectors $Y_z(s,t)x$ and $Y_z(s,t)y$ is constant (as a function of $t$) equal to the Wronskian of $x$ and $y$.
\end{lemma}
\begin{proof}
	Using the continuity in $z$ of the stochastic flow, it suffices to prove that for any given $z\in\C$, almost surely the process $Y_z(s,t)^\intercal J Y_z(s,t)$ is constant in $t$, and thus matches with its values at $t=s$ which is $J$. Fix $z\in \C$. We are going to show that the differential in $t$ vanishes for $t\ge s$ (the case $t\le s$ is similar). Since the SDE are written in the Stratonovich sense, we obtain
	\begin{align*}
		d\big( Y_z(s,t)^\intercal J Y_z(s,t)\big) &= (dY_z(s,t))^\intercal J Y_z(s,t) + Y_z(s,t)^\intercal J d Y_z(s,t)\\
		&=\big( -J(zdt- \circ d\cV(t)) Y_z(s,t) \big)^\intercal J Y_z(s,t)\\
		&+ Y_z(s,t)^\intercal (zdt- \circ d\cV(t)) Y_z(s,t)\;.
	\end{align*}
	A simple computation, based on the symmetry of the matrix $\cV(t)$ and the identity $J^\intercal = J^{-1} = - J$, shows that the sum of these two terms vanishes.
\end{proof}

\subsection{Canonical system}\label{Subsec:CanSys}

We start by recalling a few general, deterministic facts from the theory of canonical systems, that can be found in~\cite{Weidmann, Remling, Romanov}. We give ourselves a locally integrable function $(K(t),t\in\R)$ taking values in the set of $2\times 2$ real, non-negative invertible matrices.\\
For any $z\in\C$, let $(U_z(s,t), s,t \in \R)$ be the $2\times 2$ matrix-valued solution of
\begin{equation}\label{Eq:RJconstant}
	J \partial_t U_z(s,t) = z K(t)U_z(s,t)\;,\quad t\in\R\;,\mbox{ and }\quad U_z(s,s) = I\;.
\end{equation}
We now recall how the theory of canonical systems allows to associate a self-adjoint operator (in finite and infinite volume) to the expression $K(t)^{-1} J\partial_t$.

\medskip

\emph{Operator in finite volume.} For any $\tau>0$, and given two ``boundary conditions'' $x_-,x_+ \in \R^2\backslash \{(0,0)^\intercal\}$, let $\cD(\cT_\tau)$ be the set of all functions $f\in L^2_{K}((-\tau/2,\tau/2),\R^2)$ such that:
\begin{enumerate}
	\item $f$ is absolutely continuous,
	\item There exists $g\in L^2_{K}((-\tau/2,\tau/2),\R^2)$ s.t.~$J \partial_t f = K g$,
	\item $f(-\tau/2) \parallel x_-$ and $f(\tau/2) \parallel x_+$.
\end{enumerate}
The set $\cD(\cT_\tau)$ is dense in $L^2_{K}((-\tau/2,\tau/2),\R^2)$ and the operator $\cT_\tau f =g$ (where $g$ is prescribed by (2)) is self-adjoint on $\cD(\cT_\tau)$. This operator has a discrete spectrum and the following holds: $\lambda$ is an eigenvalue of $\cT_\tau$ if and only if the function $v_\lambda(t) := U_\lambda(-\tau/2,t) x_-$ satisfies $v_\lambda(\tau/2) \parallel x_+$. When this last condition holds, $v_\lambda$ is an eigenfunction associated with $\lambda$.\\

\emph{Operator in infinite volume.} We say that the operator is limit point at $+\infty$ (resp.~$-\infty$) if the matrix-valued function $K$ is not integrable at $+\infty$ (resp.~$-\infty$). Equivalently, the operator is limit point at $+\infty$ (resp.~$-\infty$) if there exists $z\in\C$, $s\in\R$ and $x\in \R^2$ such that $U_z(s,t)x$ does not lie in $L^2_K([s,+\infty),\R^2)$ (resp.~$L^2_K((-\infty,s],\R^2)$), see~\cite[Th.5.3]{Weidmann}.\\
In this article, we will be only concerned with operators that are limit point at both infinities, so let us concentrate on this case. Let $\cD(\cT_\infty)$ be the set of all functions $f\in L^2_{K}(\R,\R^2)$ such that:
\begin{enumerate}
	\item $f$ is absolutely continuous,
	\item There exists $g\in L^2_{K}(\R,\R^2)$ s.t.~$J \partial_t f = K g$.
\end{enumerate}
The set $\cD(\cT_\infty)$ is dense in $L^2_{K}(\R,\R^2)$ and the operator $\cT_\infty f =g$ (where $g$ is prescribed by (2)) is self-adjoint on $\cD(\cT_\infty)$.\\

Let us apply the above facts to our setting. We set $K(t) := Y_0(0,t)^\intercal \, Y_0(0,t)$, $t\in\R$. By Proposition \ref{Prop:Flow}, we already know that this matrix is invertible at all times. Since the determinant of $Y_0(0,t)$ is nothing but the Wronskian of its two column vectors, Lemma \ref{Lemma:Wronskian} implies that $Y_0(0,t)$ has determinant $1$ for all $t \in \R$.\\

At this point, let us notice that the collection of ODEs \eqref{Eq:RJconstant} are the counterpart of the solution of SDEs \eqref{Eq:FwdSDE}. However, contrary to the latter, there is no stochastic integral involved in the ODE so that its construction is carried out in a deterministic manner and it is straightforward to check that $U_z(s,t)$ is continuous in $s$, $t$ and $z$. Actually the two sets of equations are in correspondence:

\begin{lemma}\label{Lemma:Flowvz}
	Almost surely, for all $z\in \C$ and all $t\in\R$
	$$ Y_z(0,t) = Y_0(0,t)U_z(0,t)\;.$$
\end{lemma}
\begin{remark}
	Given the process $(Y_0(0,t),t\in\R)$ and the flow $(U_z,z\in\C)$, it is possible to \emph{define}
	$$Y_z(s,t) := Y_0(0,t)U_z(0,t) (Y_0(0,s)U_z(0,s))^{-1}\;.$$
	Since the construction of $(Y_0(0,t),t\in\R)$ goes through standard SDE arguments, this might appear as a simple construction of the stochastic flow. However, we are not aware of any direct argument ensuring that the SDE property in Proposition \ref{Prop:Flow} is fulfilled by this object.
\end{remark}
\begin{proof}
	Set $Z(t) := Y_0(0,t) U_z(0,t)$. Then for $t\ge 0$
	\begin{align*}
		d (J Z(t)) &= -\circ d\cV(t) Z(t) + J Y_0(0,t) J^{-1} z K(t) U_z(0,t)dt\\
		&= -\circ d\cV(t) Z(t) + z Z(t)dt\;,
	\end{align*}
	where the second equality follows from \eqref{Eq:JMJ}. We thus deduce that for any given $z\in\C$, $(Y_z(0,t),t\ge 0)$ and $(Y_0(0,t) U_z(0,t),t\ge 0)$ satisfy the same SDE, and thus coincide a.s. From the continuity in $z$ of the stochastic flow and of the set of solutions of \eqref{Eq:RJconstant}, we deduce that the two processes coincide simultaneously for all $z$, almost surely. The arguments backward in time are virtually the same.
\end{proof}

We are now in position to apply, in a deterministic way, the theory of canonical systems recalled above. We start with the operator in finite-volume. We give ourselves some boundary conditions $b_\pm$ that are vectors in $\R^2$ with unit Euclidean norm. For any $\tau>0$, consider the boundary conditions $x_\pm := Y_0(0,\pm \tau/2)^{-1} b_\pm$, and let $\cT_\tau$ be the self-adjoint operator introduced by the arguments above with domain $\cD(\cT_\tau)$.

We now turn to the operator in infinite-volume. We first check that the operator is limit point at both infinities.

\begin{lemma}\label{Lemma:LimitPoint}
	Fix a non-zero $x\in\R^2$. Almost surely, $Y_0(0,t)x$ is square integrable neither towards $+\infty$ nor towards $-\infty$. As a consequence we are almost surely in the limit point case at both infinities.
\end{lemma}
\begin{proof}
	By linearity we can assume that $x$ has unit Euclidean norm. Pick some $\theta \in [0,2\pi)$ such that $x = (\sin \theta, \cos \theta)^\intercal$. Given \eqref{Eq:BwdSDEbis}, we observe a symmetry in law between the solutions on $[0,\infty)$ and on $(-\infty,0]$, and we thus concentrate on the behaviour at $+\infty$. Then $\| Y_0(0,t)x \|_{L^2([0,+\infty),\R^2)}^2 = \int_{\R_+} \exp(2 \rho_0(t)) dt$ where $\rho_0$ satisfies \eqref{Eq:FwdPolarSDE} starting from $\rho_0(0) = 0$ and $\theta_0(0) = \theta$. The process $(\rho_0(t),t\ge 0)$ has the same law as the process $(\kappa^2 t + \kappa B(t), t\ge 0)$ where $B$ is a standard Brownian motion. If $\kappa=0$ it is constant equal to $0$, while if $\kappa>0$ it goes to $+\infty$ a.s. Consequently in both cases almost surely $\| Y_0(0,t)x \|_{L^2([0,+\infty),\R^2)}^2 = +\infty$.\\
	Let us now deduce that we are in the limit point case at $+\infty$. Recall that $K(t) = Y_0(0,t)^\intercal Y_0(0,t)$ so that
	$$ \int_0^\infty x^\intercal K(t) x \,dt = \int_0^\infty (Y_0(0,t) x)^\intercal (Y_0(0,t) x) dt = \| Y_0(0,t)x \|_{L^2([0,+\infty),\R^2)}^2\;,$$
	so that $K$ is not integrable at $+\infty$. The same holds at $-\infty$.
\end{proof}
Using the elements of theory of canonical systems in infinite volume recalled above, we can define the self-adjoint operator $\cT_\infty$ on the domain $\cD(\cT_\infty)$.

\subsection{The Dirac operator}

Let us observe that multiplication by $Y_0(0,\cdot)$ defines a unitary map from $L^2_K(E,\R^2)$ onto $L^2_{I}(E,\R^2)$, for $E = (-\tau/2,\tau/2)$ or $E = \R$. Indeed, for $v \in L^2_K(E, \R^2)$, we have
\begin{align*}
	\int_E |Y_0(0,t) v(t)|^2 dt = \int_E v(t)^\intercal K(t) v(t) dt \,.
\end{align*}

For any $\tau > 0$, we define $\cS_\tau := Y_0(0,\cdot) \cT_\tau Y_0(0,\cdot)^{-1}$ which is self-adjoint on the domain
$$ \cD(\cS_\tau) := \Big\{f \in L^2_{I}([-\tau/2,\tau/2], \R^2):  Y_0(0,\cdot)^{-1} f \in \cD(\cT_\tau)\Big\}\;.$$
Note that any $f\in \cD(\cS_\tau)$ satisfies the prescribed b.c.
$$ f(-\tau/2) \parallel b_-\;,\quad f(\tau/2) \parallel b_+\;.$$

Similarly, we set $\cS_\infty := Y_0(0,\cdot) \cT_\infty Y_0(0,\cdot)^{-1}$ which is self-adjoint on the domain
$$ \cD(\cS_\infty) := \Big\{f \in L^2_{I}(\R,\R^2):  Y_0(0,\cdot)^{-1} f \in \cD(\cT_\infty)\Big\}\;.$$

The following proposition shows that the operator $\cS_\tau$ fulfils the objective that we gave ourselves: its spectrum is dictated by those solutions of \eqref{Eq:FwdSDE} that satisfy the b.c. (recall that the theory of canonical systems ensure that the spectrum is discrete, and thus pure point). 

\begin{proposition}\label{Prop:Spec}
	There exists an event of probability one on which the following property holds. For any $\tau > 0$, any $b_\pm$ and any $\lambda \in \R$: $\lambda$ is an eigenvalue of $\cS_\tau$ if and only if one of the following equivalent conditions holds:
	\begin{enumerate}
		\item $Y_\lambda(-\tau/2,\tau/2)b_- \parallel b_+$,
		\item $Y_\lambda(\tau/2,-\tau/2)b_+ \parallel b_-$.
	\end{enumerate}
\end{proposition}

Note that this proposition implies that the spectrum is simple.
\begin{proof}
	We work on an event of probability one on which the statement of Lemma \ref{Lemma:Flowvz} holds. Fix $\tau > 0$ and $b_\pm$. For all $\lambda \in \R$: $\lambda$ is an eigenvalue of $\cS_\tau$ if and only if $\lambda$ is an eigenvalue of $\cT_\tau$ if and only if $v_\lambda(t) := U_\lambda(-\tau/2,t) x_-$ with $x_- := Y_0(0,-\tau/2)^{-1}b_-$ satisfies $v_\lambda(\tau/2) \parallel x_+$ with $x_+:= Y_0(0,\tau/2)^{-1}b_+$. Using the invertibility of $Y_0(0,\tau/2)$, the latter condition is equivalent to: $Y_0(0,\tau/2)v_\lambda(\tau/2) \parallel b_+$. Using the content of Lemma \ref{Lemma:Flowvz} we get 
	\begin{align*}
		Y_0(0,\tau/2)v_\lambda(\tau/2) &= Y_0(0,\tau/2)U_\lambda(0,\tau/2)U_\lambda(-\tau/2,0)Y_0(0,-\tau/2)^{-1}b_-\\
		&= Y_\lambda(-\tau/2,\tau/2)b_-\;.
	\end{align*}
	Consequently the latter condition is equivalent to the first item of the statement. The second item is equivalent to the first upon applying the invertible matrix $Y_{\lambda}(\tau/2,-\tau/2)$.
\end{proof}

In infinite volume, nothing ensures in general that the spectrum is pure point. However, the following proposition shows that the pure point component admits a simple characterization in terms of the stochastic flow.

\begin{proposition}\label{Prop:CaracEigenvalueTinfty}
	Almost surely, $\lambda$ is an eigenvalue of $\cS_\infty$ if and only if there exists a non-zero $x\in \R^2$ such that $t\mapsto Y_\lambda(0,t)x$ lies in $L^2_{I}(\R)$.
\end{proposition}
\begin{proof}
	We work on an event of probability one on which the statement of Lemma \ref{Lemma:Flowvz} holds. Then: $\lambda$ is an eigenvalue of $\cS_\infty$ if and only if it is an eigenvalue of $\cT_\infty$. Given the explicit domain of $\cT_\infty$, we deduce that $\lambda$ is an eigenvalue of $\cT_\infty$ if and only if there exists $x\in \R^2$ such that $v_\lambda(t) := U_\lambda(0,t)x$ lies in $L^2_K(\R,\R^2)$. By Lemma \ref{Lemma:Flowvz}, $Y_0(0,t) v_\lambda(t) = Y_\lambda(0,t) x$ so that $v_\lambda$ lies in $L^2_K(\R,\R^2)$ if and only if $Y_\lambda(0,t) x$ lies in $L^2_{I}(\R)$.
\end{proof}

\subsection{Dependence on $\kappa_0$ and $\kappa$}\label{Subsec:kappa}

So far, we have worked with generic values of the parameters $\kappa_0,\kappa\ge 0$ and we have kept implicit the dependence of the operators with respect to these two parameters. In this subsection, we show how to ``remove'' $\kappa_0$ and how to ``rescale'' $\kappa$. We first establish a technical lemma.

\begin{lemma}\label{Lemma:CalculSto}
	Let $(\theta(t),t\ge 0)$ be a continuous $\R$-valued semimartingale and set
	$$Q(t) = \mat{\cos \theta(t)}{-\sin \theta(t)}{\sin \theta(t)}{\cos \theta(t)} \;,\quad t\ge 0\;.$$
	Then the matrix-valued process $N(t) = Q(t) Y_z(0,t)$, $t\ge 0$ satisfies (in the sense of Stratonovich calculus)
	\begin{align}
		JdN(t) &=\Big(- Q(t)  (\circ d\cV(t) + \circ d\theta(t) \, I) Q(t)^{-1} + z \, I dt\Big) N \, . \notag
	\end{align}	
\end{lemma}
\begin{proof}
	We have
	$$ JdN(t) = JdQ(t) Y_z(0,t) + JQ(t)J^{-1} J dY_z(0,t)\;,\quad t\ge 0\;.$$
	Since $dQ(t) = J (\circ d\theta(t)) Q(t)$ and $QJ = JQ$ (both are rotation matrices), a straightforward computation leads to the asserted result.
\end{proof}

\begin{proposition}[Removing $\kappa_0$]\label{Prop:Variousk}
	Take $\kappa_0>0$ and $\kappa \ge 0$. Almost surely, the map $\cQ \;:\; f\mapsto Q f$ with
	\begin{equation}\label{Eq:Q}
		Q(t) :=  \begin{pmatrix} \cos(- \kappa_0 B_0(t)) & -\sin (-\kappa_0 B_0(t))\\ \sin (- \kappa_0 B_0(t)) & \cos(- \kappa_0 B_0(t))\end{pmatrix}\;,\quad t\in\R\;,
	\end{equation}
	is unitary on $L^2(\R,\R^2)$ and we have the identity in law $$\cQ \; \cS_{\infty,\kappa_0,\kappa} \; \cQ^{-1} \overset{(d)}{=} \cS_{\infty,0,\kappa} \,.$$
\end{proposition}
\begin{proof}
	Fix $\kappa_0 > 0$ and $\kappa\ge 0$. Let $Y_z(s,t)$ be the stochastic flow associated with the noise $(\kappa_0\dot{B}_0,\kappa \dot{B}_1, \kappa \dot{B}_2)$ and denote by $\cS_{\infty,\kappa_0,\kappa}$ the associated operator. Consider the process $\theta(t) := -\kappa_0 B_0(t)$, $t\ge 0$ so that the matrix $Q(t)$ defined in~\eqref{Eq:Q} matches with that defined in Lemma \ref{Lemma:CalculSto}. Recall the Pauli matrices from the introduction. We set $N(t) := Q(t) Y_z(0,t)$, $t\ge 0$ and deduce from that lemma that
	\begin{align*}
		JdN(t) &= \Big(- Q(t) \sigma_1 Q(t)^{-1} \kappa  \circ dB_1(t) - Q(t) \sigma_3 Q(t)^{-1} \kappa  \circ dB_2(t) + z \, I dt\Big) N(t) \\
		&= \Big(-\kappa  \sigma_1 \circ d \tilde B_1(t) -  \kappa \sigma_3 \circ d\tilde B_2(t) + z \, I dt\Big) N(t) 
	\end{align*}
	where
		\begin{equation}\label{Eq:BBtilde}\begin{split}
					d\tilde{B}_1(t) &:= \cos(-2\kappa_0B_0(t)) \circ d{B}_1(t) - \sin(-2\kappa_0B_0(t)) \circ d{B}_2(t)\;,\\
					d\tilde{B}_2(t) &:= \cos(-2\kappa_0B_0(t)) \circ d{B}_2(t) + \sin(-2\kappa_0B_0(t)) \circ d{B}_1(t)\;.
				\end{split}
		\end{equation}
	Since $B_0$ is independent from $(B_1,B_2)$, these stochastic differentials can also be understood in the It\^o sense. As $\langle\tilde{B}_1,\tilde{B}_1 \rangle_t = \langle \tilde{B}_2, \tilde{B}_2 \rangle_t = t$ and $\langle \tilde{B}_1,\tilde{B}_2\rangle_t = 0$, the processes $(\tilde{B}_1, \tilde{B}_2)$ are two independent Brownian motions.\\
	Denote by $\tilde{Y}_z(s,t)$ the stochastic flow associated to the noises $(0,\kappa \dot{\tilde{B}}_1, \kappa \dot{\tilde{B}}_2)$.	The associated operator $\tilde{\cS}_{\infty,0,\kappa}$ is given by $\tilde{Y}_0(0,\cdot) \tilde{\cT}_\infty \tilde{Y}_0(0,\cdot)^{-1}$ where $\tilde{\cT}_\infty$ is the canonical system associated to $\tilde{K}(t) := \tilde{Y}_0(0,t)^\intercal \tilde{Y}_0(0,t)$.\\
	The computation above shows that for all $z\in\C$, the processes $\tilde{Y}_z(0,\cdot)$ and $Q(\cdot) Y_z(0,\cdot)$ follow the same SDE starting from the same initial condition at $t=0$. From the continuity properties of the stochastic flows, it implies that almost surely for all $z\in\C$ and all $t\in\R$, $\tilde{Y}_z(0,t) = Q(t) Y_z(0,t)$. Since at any time $t\in\R$, $Q(t)$ is a rotation matrix, we deduce that almost surely for all $t\in\R$
	$$ \tilde{K}(t) = \tilde{Y}_0(0,t)^\intercal \tilde{Y}_0(0,t) = {Y}_0(0,t)^\intercal {Y}_0(0,t) = K(t)\;,$$
	so that $\tilde{\cT}_\infty = \cT_\infty$, and therefore
	\begin{align*}
			\tilde{\cS}_{\infty,0,\kappa} &= \tilde{Y}_0(0,\cdot) \tilde{\cT}_\infty \tilde{Y}_0(0,\cdot)^{-1}= Q(\cdot) Y_0(0,\cdot) \cT_\infty Y_0(0,\cdot)^{-1} Q(\cdot)^{-1}= Q(\cdot) \cS_{\infty,\kappa_0,\kappa} Q(\cdot)^{-1}\;.
		\end{align*}
	from which we deduce Proposition \ref{Prop:Variousk}.
\end{proof}

\begin{proposition}[Rescaling $\kappa$]
	Assume that $\kappa > 0$ and consider the unitary map $U_\kappa: f \mapsto \kappa^{-1} f(\cdot/\kappa^2)$. Then
	\begin{align*}
		\kappa^{-2} U_\kappa	\cS_{\infty,\kappa_0,\kappa} U_\kappa^{-1} \stackrel{(d)}{=}  \cS_{\infty,\kappa_0/\kappa,1} \,.
	\end{align*}
\end{proposition}
The proof is relatively straightforward so we do not present it.

\begin{remark}[Relation with the critical operator $ \tSch_\tau$ in~\cite{DL_Critical}]\label{Rem:highenergy}
	In a previous work~\cite{DL_Critical}, we constructed the operator
	$$ \tSch_\tau := 2J \partial_t + \sqrt{\tau} \begin{pmatrix} \dot{B}_0(t) + \frac1{\sqrt 2}\dot{B}_1(t) & \frac1{\sqrt 2}\dot{B}_2(t)\\ \frac1{\sqrt 2}\dot{B}_2(t) &\dot{B}_0(t)-\frac1{\sqrt 2} \dot{B}_1(t)\end{pmatrix}\;,\quad t\in (0,1)\;,$$
	which coincides in law with $\cS_{\tau/4,1,1/\sqrt{2}}$ up to an affine transformation. More precisely we can construct a coupling such that $(\mu,\Phi_\mu)$ is an eigenpair of $\tSch_\tau$ if and only if $(\lambda,\Psi_\lambda)$ is an eigenpair of $\cS_{\tau/4,1,1/\sqrt{2}}$ with
	$$ \mu = \frac{\tau \lambda}{2}\;,\quad \Phi_\mu(t) := \sqrt{\frac{\tau}{4}} \Psi_\lambda(\frac{\tau}{4}t - \frac{\tau}{8})\;.$$
\end{remark}

\subsection{Weyl's theory and the resolvent operators}\label{Subsec:Weyl}

We now present Weyl's theory for our Dirac operator $\cS_\infty$. The material is well-known for canonical systems, see for instance~\cite[Sec.~3.4 to 3.7]{Remling}, and is readily transferred to Dirac operators using the explicit mapping between the two objects. In this subsection, we work with generic values $\kappa_0,\kappa\ge 0$.

For convenience, we single out the so-called Neumann and Dirichlet solutions of the stochastic flow (recall~\eqref{Eq:e1e2}):
\begin{equation}\label{Eq:DirichletNeumann}
	y_z^N(t) := Y_z(0,t)e_1\;,\quad y_z^D(t) := Y_z(0,t)e_2\;,\quad t\in\R\;,
\end{equation}
for any $z\in\C$. Fix $\tau > 0$ and $z\in \C\backslash\R$. For any $\beta \in [0,\pi)$, there exists a unique $m^{+,\tau}_z(\beta) \in\C$ such that the random vector
\begin{equation*}
	y_z^N(\tau/2) + m^{+,\tau}_z(\beta) y_z^D(\tau/2)\;,
\end{equation*}
is parallel to $\begin{pmatrix} \sin \beta \\ \cos\beta \end{pmatrix}$. Furthermore, $\beta \mapsto m^{+,\tau}_z(\beta)$ forms a circle in $\C\backslash\R$, whose radius is given by
$$ \frac1{2 |\Im(z)| \int_0^{\tau/2} |  y_z^D(t)|^2 dt}\;.$$
This expression readily shows that the radius decreases with $\tau$. Actually, the corresponding disks are embedded into each others. Consequently as $\tau \to +\infty$, these disks converge either to a non-trivial disk or to a point. In the case where these disks converge to a point, we denote this point $m^+_z$. Actually, the following dichotomy holds: in the limit circle case, all disks (indexed by $z\in\C\backslash\R$) converge to non-trivial disks ; in the limit point case, all disks converge to points.

Towards $-\infty$, we introduce the analogous objects: let $m^{-,\tau}_z(\beta) \in\C$ be such that
$$ y^N_z(-\tau/2) - m^{-,\tau}_z(\beta) y^D_z(-\tau/2)\;,$$
is parallel to $\begin{pmatrix} \sin \beta \\ \cos\beta \end{pmatrix}$. It enjoys the same properties as $m^{+,\tau}_z(\beta)$.\\

By Lemma \ref{Lemma:LimitPoint}, we know that our operator is limit point at both infinities, almost surely. We then denote the points that arise in the limit by $m^\pm_z := \lim_{\tau \to \infty} m^{\pm,\tau}_z(\beta)$. Weyl's theory then ensures the following fact: for any $z\in\C\backslash\R$, $t\mapsto Y_z(0,t)x$ is square-integrable at $+\infty$ (resp.~$-\infty$) if and only if $x$ is parallel to $e_1+m^+_ze_2$ (resp.~$e_1-m^-_z e_2$).

Recall that our finite-volume operators are endowed with some b.c.~$b_\pm$. We let $\beta_\pm \in [0,\pi)$ be such that $b_\pm$ is parallel to $\begin{pmatrix} \sin \beta_\pm \\ \cos\beta_\pm \end{pmatrix}$. The following proposition provides an explicit expression for the kernels of the resolvent operators, see~\cite[Theorem 3.3]{Remling}.

\begin{proposition}\label{Prop:Kernel}
	Fix $\tau > 0$ and $z \in \C\backslash\R$. Set
	$$ y^{+,\tau}_z(t) := y_z^N(t) + m^{+,\tau}_z(\beta_+) y_z^D(t)\;,\quad y^{-,\tau}_z(t) := y_z^N(t) - m^{-,\tau}_z(\beta_-) y_z^D(t)\;,\quad t\in \R\;.$$
	Then almost surely, for any function $g\in L^2_{I}(-\tau/2,\tau/2)$ we have
	$$ (\cS_\tau - z)^{-1}g(t) = \int_{s \in [-\tau/2,\tau/2]} (\cS_\tau - z)^{-1}(s,t) g(s) ds\;,$$
	where for any $s,t \in [-\tau/2,\tau/2]$
	$$ (\cS_\tau - z)^{-1}(s,t) = \frac1{W(y^{+,\tau}_z,y^{-,\tau}_z)} \Big( y^{+,\tau}_z(t) y^{-,\tau}_z(s)^\intercal \mathbf{1}_{s\le t} + y^{-,\tau}_z(t) y^{+,\tau}_z(s)^\intercal \mathbf{1}_{s> t} \Big)\;.$$
	Similarly if we set
	$$ y^+_z(t) := y_z^N(t) + m^{+}_z y_z^D(t)\;,\quad y^-_z(t) := y_z^N(t) - m^{-}_z y_z^D(t)\;,\quad t\in \R\;,$$
	then almost surely, for any function $g\in L^2_{I}(\R)$ we have
	$$ (\cS_\infty - z)^{-1}g(t) = \int_{s \in \R} (\cS_\infty - z)^{-1}(s,t) g(s) ds\;,$$
	where for any $s,t\in\R$
	$$ (\cS_\infty - z)^{-1}(s,t) = \frac1{W(y^{+}_z,y^{-}_z)} \Big( y^{+}_z(t) y^{-}_z(s)^\intercal \mathbf{1}_{s\le t} + y^{-}_z(t) y^{+}_z(s)^\intercal \mathbf{1}_{s> t} \Big)\;.$$
\end{proposition}

From these identities, we deduce that each operator $(\cS_\tau-z)^{-1}$ is Hilbert-Schmidt and thus compact, so that $\cS_\tau$ has discrete spectrum.\\

We now have all the ingredients at hand to prove the first result presented in the introduction.
\begin{proof}[Proof of Proposition \ref{Prop:Construction}]
	Since our operator is limit point at both infinities, we have $m^{\pm,\tau}_z(\beta_\pm) \to m^{\pm}_z$ as $\tau\to\infty$, and consequently $y^{\pm,\tau}_z$ converges to $y^\pm_z$ locally uniformly as $\tau\to\infty$. From the explicit expressions of the kernels of the Green's functions, it is easy to check that for any functions $g,h\in L^2_{I}$ with compact supports
	$$ \langle h, (\cS_\tau - z)^{-1}g\rangle \to \langle h, (\cS_\infty - z)^{-1}g\rangle\;,\quad \tau \uparrow \infty\;.$$
	Since the operator norms of $(\cS_\tau-z)^{-1}, (\cS_\infty - z)^{-1}$ are bounded (by $|\Im(z)|^{-1}$), this convergence remains true without the compactness of the supports of $g,h$. This implies that $\cS_\tau$ converges to $\cS_\infty$ in the weak resolvent sense (using the terminology of Weidmann~\cite{WeidmannCV}, meaning that $(\cS_\tau-z)^{-1}P_\tau$ converges weakly to $(\cS_\infty-z)^{-1}$ if $P_\tau$ is the orthogonal projection of $L^2(\R,\R^2)$ onto the space $L^2((-\tau/2,\tau/2),\R^2)$ mapping $f$ to $\un_{[-\tau/2,\tau/2]} f$, see e.g. Section 6.6 of  \cite{Teschl}). By~\cite[Lemma 6.37]{Teschl}, this implies strong resolvent convergence.
\end{proof}
	We now state and prove the symmetry properties mentioned in the introduction.
\begin{proposition}[Spectral invariance]\label{Prop:SymmetrySpectralShift}
		Fix $\kappa_0,\kappa\ge 0$. For any $E\in\R$, let $\cQ_E$ be the unitary map $f \mapsto Q_E f$ where
		$$Q_E := \mat{\cos (E \; \cdot)}{-\sin(E \; \cdot)}{\sin(E \; \cdot)}{\cos(E \; \cdot)}\, .$$
		The operator $\cQ_E \big( \cS_{\infty,\kappa_0,\kappa} - E \big) \cQ_E^{-1}$
		has the same law as the operator $\cS_{\infty,\kappa_0,\kappa}$.
\end{proposition}
\begin{proof}
	Fix $E\in\R$ and set
	\begin{equation*}
		Q_E(t) :=  \begin{pmatrix} \cos(Et) & -\sin (Et)\\ \sin (Et) & \cos(Et)\end{pmatrix}\;,\quad t\in\R\;.
	\end{equation*}
	Let $\tilde{Y}$ be the stochastic flow associated to the noises $(\kappa_0 \dot{B}_0,\kappa \dot{\tilde{B}}_1, \kappa \dot{\tilde{B}}_2)$ where
	\begin{equation*}\begin{split}
			d\tilde{B}_1(t) &:= \cos(2 Et) d{B}_1(t) - \sin(2 Et) d{B}_2(t)\;,\\
			d\tilde{B}_2(t) &:= \cos(2 Et) d{B}_2(t) + \sin(2 Et) d{B}_1(t)\;.
		\end{split}
	\end{equation*}
	(note that $\tilde{B}_1$, $\tilde {B}_2$ are two independent Brownian motions).
	Using Lemma \ref{Lemma:CalculSto} with $\theta(t) := Et$ and again the continuity properties of the stochastic flows, we deduce that $\tilde{Y}_{z}(0,t) = Q_E(t) Y_{z+E}(0,t)$ for all $t\in\R$ and all $z\in\C$.\\
	The definition of $\cS_\infty$ and Lemma \ref{Lemma:Flowvz} which gives $Y_E(0,\cdot) = Y_0(0,\cdot)U_E(0,\cdot)$ yield
	\begin{align*}
		\cS_\infty-E &= Y_E(0,\cdot) U_E^{-1}(0,\cdot)(\cT_\infty - E) U_E(0,\cdot) Y_E(0,\cdot)^{-1}\\
		&= Q_E^{-1}(\cdot) \tilde{Y}_0(0,\cdot) U_E^{-1}(0,\cdot)(\cT_\infty - E) U_E(0,\cdot) \tilde{Y}_0(0,\cdot)^{-1} Q_E(\cdot)\;.
	\end{align*}
	We claim that $U_E^{-1}(0,\cdot)(\cT_\infty - E) U_E(0,\cdot)$ coincides with the operator $\tilde{\cT}_\infty$ associated to
	$$\tilde{K}(t) := Y_E(0,t)^\intercal Y_E(0,t)\;.$$
	If the claim holds, then, using that $$Y_E(0,t)^\intercal Y_E(0,t) = (Q_E^{-1}(t)\tilde{Y}_{0}(0,t))^\intercal Q_E^{-1}(t)\tilde{Y}_{0}(0,t) = \tilde{Y}_{0}(0,t)^\intercal \tilde{Y}_{0}(0,t)\,,$$ 
	and the fact that $\tilde{Y}$ has the same law as $Y$, we deduce that $\tilde{Y}_0(0,\cdot) \tilde{\cT}_\infty \tilde{Y}_0(0,\cdot)^{-1}$ has the same law as $\cS_{\infty}$, and consequently $Q_E(\cS_\infty-E)Q_E^{-1}$ is equal in law to $\cS_{\infty}$.\\	
	We turn to the proof of the claim. For simplicity, we set $U_E(\cdot) := U_E(0,\cdot)$ and we recall that $Y_E(0,\cdot) = Y_0(0,\cdot) U_E(\cdot)$. We note that $J \partial_t U_E  = E K U_E$ and $U_E^\intercal J U_E = J$ (since $\det U_E = 1$).
	Take $\tilde{f} \in \cD(\tilde{\cT}_\infty)$. It is straightforward that $f:= U_E \tilde{f}$ belongs to $L^2_K$. Furthermore
	\begin{align*}
		J \partial_t f &= J(\partial_t U_E) \tilde f + J U_E \partial_t \tilde f = E KU_E \tilde{f} + (U_E^\intercal)^{-1} \tilde K (\tilde K)^{-1} J \partial_t \tilde f\\
		&= E K f + (U_E^\intercal)^{-1} \tilde K \tilde{\cT}_\infty \tilde{f}\;.
	\end{align*}
	Since $ U_E^\intercal K U_E  = \tilde K$, we deduce that
	$$ K^{-1} J \partial_t f = E f + U_E \tilde{\cT}_\infty \tilde{f}\;.$$
	The r.h.s.~belongs to $L^2_K$. Consequently $f\in\cD(\cT_\infty)$ and $(\cT_\infty - E) f = U_E \tilde{\cT}_\infty \tilde{f}$. The same arguments show that the same identity holds if we start from $f\in\cD(\cT_\infty)$ and set $\tilde{f} := U_E^{-1} f$.
\end{proof}
\begin{proposition}[Spatial invariance]\label{Prop:SymmetrySpatialShift}
	Fix $\kappa_0,\kappa\ge 0$. For any $t_0\in\R$, the operator $T_{t_0}\cS_{\infty,\kappa_0,\kappa} T_{t_0}^{-1}$ has the same law as $\cS_{\infty,\kappa_0,\kappa}$, where $T_{t_0}$ is the shift operator $T_{t_0} f(t) = f(t+t_0)$, $t\in\R$.
\end{proposition}
\begin{proof}
	We fix $t_0\in\R$. Set $\tilde{Y}_z(s,t) := Y_z(s+t_0,t+t_0)$ for any $z\in\C$, and any $s,t\in\R$. Clearly $\tilde{Y}$ has the same law as $Y$, and is associated to the noises $(\kappa_0 \dot{B}_0(\cdot+t_0), \kappa \dot{B}_1(\cdot+t_0), \kappa \dot{B}_2(\cdot+t_0))$.\\
	Then we observe that
	\begin{align*}
		T_{t_0}\cS_\infty T_{t_0}^{-1}  &= Y_0(t_0,\cdot + t_0) Y_0(0,t_0) T_{t_0} \cT_\infty T_{t_0}^{-1} Y_0(0,t_0)^{-1} Y_0(t_0,\cdot + t_0)^{-1}\\
		&= \tilde{Y}_0(0,\cdot) \tilde{Y}_0(0,-t_0)^{-1}  T_{t_0} \cT_\infty T_{t_0}^{-1} \tilde{Y}_0(0,-t_0) \tilde{Y}_0(0,\cdot)^{-1}\;.
	\end{align*}
	Thus it suffices to prove that $\tilde{Y}_0(0,-t_0)^{-1} T_{t_0} \cT_\infty T_{t_0}^{-1} \tilde{Y}_0(0,-t_0)$ coincides with the operator $\tilde{\cT}_\infty$ associated to $\tilde{K}(t) := \tilde{Y}_0(0,t)^\intercal \tilde{Y}_0(0,t)$ which implies $T_{t_0} \cS_\infty T_{t_0}^{-1} = \tilde{\cS}_\infty$. The remaining arguments are very similar to those provided for the shift in energy so we leave the details to the reader.
\end{proof}

\section{Density of states and GMP formula}\label{Sec:DOSGMP}

In this section, we prove Proposition \ref{Th:IDS} on the density of states, and we compute the intensity measure of $\mu_{\tau}$ for any finite $\tau$. In both cases, the proof relies on the SDEs in polar coordinates. Recall the complex white noise introduced in~\eqref{Eq:ComplexBM}.

\subsection{Polar coordinates}\label{Sec:Polar}

For any given $s\in\R$, we consider the SDEs
\begin{equation}\label{Eq:SDEPolars}
	\begin{split}
	d\theta_\lambda(t) &= \lambda dt  - \kappa_0 dB_0(t) + \kappa\,\Re(e^{2i \theta_\lambda(t)} \circ d\cB(t))\;,\quad t\ge s\;,\\
	d\rho_\lambda(t) &= \kappa \,\Im(e^{2i \theta_\lambda(t)} \circ d\cB(t))\;,\quad t\ge s\;.
\end{split}
\end{equation}
starting from $\theta_\lambda(s) = \theta_0$ and $\rho_\lambda(s) = \rho_0$. In the sequel, we will rely on the following technical result.

\begin{proposition}\label{Prop:FlowThetaRho}
	Fix $s\in\R$. There exists a collection $\big(\theta_\lambda(\theta_0,s,t), \rho_\lambda(\theta_0,\rho_0,s,t); t\in [s,\infty), \lambda\in\R, \theta_0\in\R, \rho_0 \in \R\big)$ of r.v.~such that:\begin{enumerate}
		\item (Continuity) The map $(\lambda,t,\theta_0,\rho_0)\mapsto (\theta_\lambda(\theta_0,s,t), \rho_\lambda(\theta_0,\rho_0,s,t))$ is continuous from $\R\times [s,\infty)\times\R^2$ into $\R^2$,
		\item (Symmetry) For any $k\in\Z$ and $\rho_0'\in\R$
		\begin{align*}
			\theta_\lambda(\theta_0+k\pi,s,t) &= \theta_\lambda(\theta_0,s,t) + k\pi \;,\\
			\rho_\lambda(\theta_0+k\pi,\rho_0+\rho_0',s,t) &= \rho_\lambda(\theta_0,\rho_0,s,t) + \rho_0'\;.
		\end{align*}
		\item (Monotonicity) For any $t \ge s$, $\theta_0 \mapsto \theta_\lambda(\theta_0,s,t)$ is continuous increasing.
		\item (Differentiability) The map $\lambda \mapsto \theta_\lambda(\theta_0,s,t)$ is a $\cC^1$ diffeomorphism for any $s<t$, and if we set $X(t):=\partial_\lambda \theta_\lambda(\theta_0,s,t)$ then a.s.
		\begin{equation}\label{Eq:SDEPartial}
			d X(t) = dt - 2 \kappa  X(t) \Im(e^{2i \theta_\lambda(\theta_0,s,t)} \circ d\cB(t))\;,\quad t\in [s,\infty)\;.
		\end{equation}
		\item (SDE) The process $(\theta_\lambda(\theta_0,s,t), \rho_\lambda(\theta_0,\rho_0,s,t); t\in [s,\infty))$ is an almost sure solution of~\eqref{Eq:SDEPolars}.
		\end{enumerate}
\end{proposition}
\begin{remark}
	We could have built a stochastic flow associated to the polar coordinates as in Proposition \ref{Prop:Flow}. With such a stochastic flow, it would have been possible to prove an a.s.~convergence in Proposition \ref{Th:IDS}, but for simplicity we refrain from doing so.
\end{remark}
\begin{proof}
	For convenience, we add to the SDE \eqref{Eq:SDEPolars} a third component which is constant equal to $\lambda$. From~\cite[Chap.V, Theorems 37, 40 and 44]{Protter2005}, there exists a process $(\lambda,\theta_\lambda(\theta_0,s,t), \rho_\lambda(\theta_0,\rho_0,s,t); \lambda\in\R, \theta_0\in\R, \rho_0 \in \R, t\in [s,\infty))$ such that properties (Continuity) and (SDE) are true, and furthermore, such that almost surely for all $t\ge s$
	$$(\lambda,\theta_0,\rho_0) \mapsto \big(\lambda,\theta_\lambda(\theta_0,s,t),\rho_\lambda(\theta_0,\rho_0,s,t)\big)\;,$$
	 is injective and continuously differentiable, with non-negative partial derivatives. In particular, the derivative in $\lambda$ is given by~\eqref{Eq:SDEPartial}. The properties (Monotonicity) and (Differentiability) follow. Finally, if we fix $k$ and $\rho_0'$, then the two identities of the (Symmetry) property are satisfied a.s.~by the solutions of the SDE~\eqref{Eq:SDEPolars}. Taking a union over a countable and dense collection of parameters $k,\rho'_0$, this remains true, and we conclude by continuity.
\end{proof}

We can now spell out a characterization of the spectrum of $\cS_{\tau}$ in terms of the SDEs in polar coordinates. Recall that $\cS_{\tau}$ is endowed with some b.c.~$b_\pm$, and that $\beta_\pm \in [0,\pi)$ are chosen in such a way that $b_\pm$ is parallel to $\begin{pmatrix} \sin \beta_\pm \\ \cos\beta_\pm \end{pmatrix}$.

\begin{proposition}\label{Prop:SpecAngle}
	Fix $\tau > 0$. There exists an event of probability one on which the following property holds. For any $\lambda \in \R$: $\lambda$ is an eigenvalue of $\cS_\tau$ if and only if $\theta_\lambda(\beta_-,-\tau/2,\tau/2)$ lies in $\beta_+ + \pi\Z$.
\end{proposition}
\begin{remark}
	For this property to hold with probability one for \emph{all} $\tau>0$ simultaneously, one would need to construct a stochastic flow for the polar coordinates.
\end{remark}
\begin{proof}
	Fix $s \in \R$. The proof relies on an almost sure relationship between the stochastic flow $Y$ and the process constructed in Proposition \ref{Prop:FlowThetaRho}. Set for any $\theta,\rho\in \R$
	$$ g(\theta,\rho) :=  e^\rho (\sin \theta, \cos \theta)^\intercal\;.$$
	We claim that almost surely for all $x\in\R^2$, all $\lambda \in \R$, and all $\theta_0,\rho_0$ such that $x = g(\theta_0,\rho_0)$
	$$ Y_\lambda(s,t)x = g\big((\theta_\lambda(\theta_0,s,t),\rho_\lambda(\theta_0,\rho_0,s,t)\big)\;.$$
	If we fix $\lambda,x,\theta_0,\rho_0$ then stochastic calculus ensures that $t\mapsto g\big(\theta_\lambda(\theta_0,s,t),\rho_\lambda(\theta_0,\rho_0,s,t)\big)$ satisfies the same SDE as $t\mapsto Y_\lambda(s,t)x$. By strong uniqueness of the solution of this SDE, we deduce that almost surely the two processes coincide. This remains true on a dense countable collection of parameters $\lambda,x,\theta_0,\rho_0$. The continuity of the stochastic flow $Y$, of the process $(\theta_\lambda(\theta_0,s,t), \rho_\lambda(\theta_0,\rho_0,s,t); t\in [s,\infty), \lambda\in\R, \theta_0\in\R, \rho_0 \in \R)$, and of the map $g$ allow us to conclude the proof of the claim.\\
	We now take $s=-\tau/2 < 0$. We work on the event of probability one on which the claim and the assertions in Proposition \ref{Prop:Spec} hold. We deduce that for any $\lambda \in \R$: $\lambda$ is an eigenvalue of $\cS_\tau$ if and only if $Y_\lambda(-\tau/2,\tau/2)b_- \parallel b_+$. Recall that in polar coordinates $b_-$ corresponds to the angle $\beta_-$ modulo $\pi$, and by the (Symmetry) property in Proposition \ref{Prop:FlowThetaRho} $\theta_\lambda(\beta_-+k\pi,-\tau/2,\tau/2) = \theta_\lambda(\beta_-,-\tau/2,\tau/2) + k\pi$. Consequently $Y_\lambda(-\tau/2,\tau/2)b_- \parallel b_+$ is equivalent to $\theta_\lambda(\beta_-,-\tau/2,\tau/2) = \beta_+$ modulo $\pi$.
\end{proof}

\subsection{Density of states}\label{sec:DOS}

\begin{lemma}\label{Lemma:AngleAsymp}
	Fix $\lambda \in \R$. As $t\to\infty$, $\E[|\theta_\lambda(\beta_-,0,t) - \lambda t|] / t \to 0$.
\end{lemma}
\begin{proof}
	By \eqref{Eq:SDEPolars}, the process $(\theta_\lambda(\beta_-,0,t) - \lambda t, t\ge 0)$ is a Brownian motion with variance $\kappa_0^2+\kappa^2$ so that, once divided by $t$, it converges to $0$ a.s.~and in $L^p$ (for any $p\ge 1$) as $t\to\infty$.
\end{proof}

For $x\in \R$, we let $\{x\}_\pi$ be the unique real value in $[0,\pi)$ such that $x-\{x\}_\pi$ lies in $\pi \Z$, and we set $[x]_\pi := (x-\{x\}_\pi)/ \pi \in \Z$.

\begin{proof}[Proof of Proposition \ref{Th:IDS}]
	To alleviate the notation, we set $\theta_\lambda(t) := \theta_\lambda(\beta_-,-\tau/2,t)$ for all $t\ge -\tau/2$. Recall the definition of $N_\tau$ from \eqref{Eq:NtauStates}. From Proposition \ref{Prop:SpecAngle} and using the fact that $\lambda\mapsto \theta_\lambda(\tau/2)$ is continuous increasing, we deduce that $N_\tau(a,b)$ a.s.~coincides with the number of times the process $(a,b] \ni \lambda\mapsto \theta_\lambda(\tau/2)$ hits $\beta_++\pi \Z$ and thus a.s.
	\begin{equation}\label{Eq:NtauNice}
		N_\tau(a,b) = [\theta_b(\tau/2) - \beta_+]_\pi - [\theta_a(\tau/2) - \beta_+]_\pi\;.
	\end{equation}
	Therefore
	\begin{align*}
		\Big|\frac{N_\tau(a,b)}{\tau} - \frac{b-a}{\pi}\Big| &\le \Big| \frac{[\theta_b(\tau/2) - \beta_+]_\pi}{\tau} - \frac{b}{\pi}\Big| + \Big| \frac{[\theta_a(\tau/2) - \beta_+]_\pi}{\tau} - \frac{a}{\pi}\Big|\\
%		&\le \frac{|\theta_b(\tau/2)-\beta_+-b\tau|}{\pi \tau} + \frac{\{\theta_b-\beta_+\}_\pi}{\pi \tau} + \frac{|\theta_a(\tau/2)-\beta_+ - a\tau|}{\pi \tau} + \frac{\{\theta_a-\beta_+\}_\pi}{\pi \tau}\\
		&\le \frac2{\tau} + \frac{|\theta_b(\tau/2)-\beta_+-b\tau|}{\pi \tau}  + \frac{|\theta_a(\tau/2)-\beta_+ - a\tau|}{\pi \tau}\;.
	\end{align*}
	Since $\theta_\lambda(\tau/2)$ has the same law as $\theta_\lambda(\beta_-,0,\tau)$, Lemma \ref{Lemma:AngleAsymp} entails that the expectation of the l.h.s.~goes to $0$ as $\tau\to\infty$.
\end{proof}

\subsection{GMP formula}\label{Subsec:GMP}

Fix $\tau>0$ and some b.c.~$b_\pm \in \R^2$. For convenience, we consider the corresponding b.c.~in polar coordinates $\beta_\pm \in[0,\pi)$. We let $(\lambda_{j,\tau},\Psi_{j,\tau})_{j \in \Z}$ be the eigenvalues and normalized eigenfunctions of the operator $\cS_\tau$ (we have seen that it has compact resolvents). The result below establishes a very tractable identity on the intensity measure of the point process
$$ \mu_{\tau} := \sum_{j\in \Z} | \Psi_{j,\tau}(0)|^2 \delta_{(\lambda_{j,\tau},\Psi_{j,\tau})} \;.$$
It has the same flavor as a formula due to Goldsheid, Molchanov and Pastur~\cite{GMP} on a random Schr\"odinger operator, so we dub it \emph{GMP formula}. To spell it out, we need to introduce some notations. We give ourselves a two-sided Brownian motion $W_1$ and an independent process $Z_\tau$ distributed as a Brownian bridge from $(-\tau/2,0)$ to $(\tau/2,0)$, that is, a centered Gaussian process with covariance function
$$\E[Z_\tau(s)Z_\tau(t)] =  (\min(s,t) + \tau/2)\Big(1 - \frac{\max(s,t) + \tau/2}{\tau}\Big)\;,\quad s,t \in [-\tau/2,\tau/2]\;.$$
Then we set
$$ S(t) := \exp\Big(\kappa W_1(t) - \kappa^2|t|\Big)\;,\quad t\in\R\;,$$
and, for any $x\in\R$,
$$ \beta^{(x)}_\tau(t) := \beta_-+ \frac{x - \beta_-}{\tau} \big(t + \frac{\tau}{2}\big) + \sqrt{\kappa_0^2+\kappa^2} Z_\tau(t)\;,\quad t\in [-\tau/2,\tau/2]\;.$$
Finally we let $g_{a,b}$ be the density function of the Gaussian law of mean $a$ and variance $b$.

\begin{proposition}[GMP formula -- finite volume]\label{Prop:GMP2}
	Take $\kappa_0, \kappa \geq 0$. Fix $\tau > 0$ and let $\cS_{\tau}$ be endowed with the b.c.~$\beta_\pm$ in polar coordinates. Let $H:\R\times \cC([-\tau/2,\tau/2],\R^2) \to \R$ be a non-negative measurable s.t.~$H(\gl,\Psi) = H(\gl,-\Psi)$. Then
	\begin{align}
		&\E\Big[\sum_{j\in \Z} | \Psi_{j,\tau}(0)|^2 H(\lambda_{j,\tau},\Psi_{j,\tau})\Big] \label{e:Htau}\\
		&= \int_{\R}\sum_{n\in\Z}  \E\bigg[ H\bigg(\lambda , \frac{S(\cdot)}{\| S\|_\tau} \begin{pmatrix} \sin \beta^{(\beta_+  +n\pi)}_\tau (\cdot) \\\cos \beta^{(\beta_+  +n\pi)}_\tau (\cdot) \end{pmatrix}\bigg)\bigg] g_{\lambda \tau,(\kappa_0^2+\kappa^2)\tau}(\beta_+-\beta_-+n\pi)d\lambda\;, \notag
\end{align}
where $\| S\|_\tau$ denotes the $L^2([-\tau/2,\tau/2],dt)$ norm of $S$.
\end{proposition}
\begin{remark}
	In the particular case where $\kappa_0=\kappa=0$, the two processes $S$ and $\beta^{(x)}_\tau$ are deterministic and $g_{\lambda \tau,(\kappa_0^2+\kappa^2)\tau}$ must be interpreted as a Dirac mass at $\lambda\tau$.
\end{remark}
\begin{proof}
	We work under the assumption that $\kappa_0$ and $\kappa$ are not simultaneously equal to $0$; however the adaptation of the proof to $\kappa_0=\kappa=0$ is straightforward. We will prove that for any bounded, continuous, non-negative map $G:\R\times \cC([-\tau/2,\tau/2],\R^2) \to \R$ compactly supported in its first argument, and s.t.~$G(\lambda,\Psi)=G(\lambda,-\Psi)$ the following identity holds
	\begin{align}
			&\E\Big[\sum_{j\in\Z} G(\lambda_{j,\tau},\Psi_{j,\tau})\Big] \label{e:Gtau} \\
			&= \int_{\R}\sum_{n\in\Z}  \E\bigg[ \| S\|_\tau^2 G\bigg(\lambda , \frac{S(\cdot)}{\| S\|_\tau} \begin{pmatrix} \sin \beta^{(\beta_+  +n\pi)}_\tau (\cdot) \\\cos \beta^{(\beta_+ +n\pi)}_\tau (\cdot) \end{pmatrix}\bigg)\bigg] g_{\lambda \tau,(\kappa_0^2+\kappa^2)\tau}(\beta_+-\beta_-+n\pi)d\lambda\;. \notag
	\end{align}
	Once this identity is established, by monotone convergence it remains true without the boundedness and compactness assumptions on the support of $G$, in particular we can pick
	$$ G(\lambda,\Psi) := |\Psi(0)|^2 H(\lambda,\Psi)\;,$$
	for any non-negative continuous map $H$ s.t.~$H(\lambda,\Psi)=H(\lambda,-\Psi)$ (notice that the factor $|\Psi_{j,\tau}(0)|^2$ on the LHS becomes $1/\|S\|^2_\tau$ on the RHS). By a Monotone Class argument, this remains true for any non-negative measurable map $H$ and we recover~\eqref{e:Htau}.\\
	We thus concentrate on proving~\eqref{e:Gtau}. Let $a>0$ be such that $G(\lambda,\cdot) = 0$ whenever $\lambda \notin [-a,a]$. For $t\in [-\tau/2,\tau/2]$, let us set
	$$\theta_\lambda(t) := \theta_\lambda(\beta_-,-\tau/2,t)\;,\quad \rho_\lambda(t) := \rho_\lambda(\beta_-,0,-\tau/2,t)\;,$$
	where the latter processes come from Proposition \ref{Prop:FlowThetaRho}. Using Proposition \ref{Prop:SpecAngle}, we have almost surely
	\begin{align*}
		\sum_{j\in \Z} G(\lambda_{j,\tau},\Psi_{j,\tau}) &= \sum_{\lambda\in \R: \theta_\lambda(\tau/2) \in \beta_+ + \pi \Z} G(\lambda,\frac{Y_\lambda(-\tau/2,\cdot)b_-}{\|Y_\lambda(-\tau/2,\cdot)b_-\|_\tau})\\
		&= \sum_{\lambda\in \R: \theta_\lambda(\tau/2) \in \beta_+ + \pi \Z} \tilde{G}(\lambda, \theta_\lambda,\rho_\lambda)\;,
	\end{align*}
	where $\tilde{G}$ is defined through
	\begin{equation}\label{Eq:GGtilde}
		\tilde{G}(\lambda,\theta,\rho) = G\Big(\lambda, \frac{e^\rho}{\| e^\rho\|_\tau} \begin{pmatrix} \sin \theta\\\cos \theta\end{pmatrix}\Big)\;,
	\end{equation}
	By Proposition \ref{Prop:SpecAngle}, the map $\Theta \,:\, \lambda \mapsto \theta_\lambda(\tau/2)$ is a.s.~a $\cC^1$-diffeomorphism: we denote by $\theta\mapsto \Lambda(\theta)$ its inverse. Then almost surely
	\begin{align*}
		\sum_{j\in \Z} G(\lambda_{j,\tau},\Psi_{j,\tau})= \sum_{n\in \Z}  \tilde{G}({\Lambda(\beta_+ +n\pi)}, \theta_{\Lambda(\beta_+ +n\pi)},\rho_{\Lambda(\beta_++n\pi)})\;.
	\end{align*}
	From the (Continuity) property of Proposition \ref{Prop:FlowThetaRho}, this last r.v.~is the a.s.~limit as $\epsilon\downarrow 0$ of
	$$ X_\epsilon := \frac1{2\epsilon} \int_{\theta \in [\beta_+  - \epsilon,\beta_+ + \epsilon]+\pi \Z} \tilde{G}({\Lambda(\theta)}, \theta_{\Lambda(\theta)},\rho_{\Lambda(\theta)}) d\theta\;.$$
	Provided $\epsilon < \pi/2$, we find the bound
	\begin{align*}
		 |X_\epsilon| &\le \frac{1}{2\epsilon}	\|\tilde{G}\|_\infty  \int_{\theta \in [\beta_+ -\epsilon,\beta_+ + \epsilon] + \pi \Z} 1_{\theta \in [\Theta(-a),\Theta(a)]} d\theta \\
		 &\le 	\|\tilde{G}\|_\infty \# \{n \in \Z,\;  [\beta_+ -\epsilon + n \pi,\beta_+ + \epsilon + n \pi] \cap[\Theta(-a),\Theta(a)] \neq \emptyset  \}\\
%		 &\le |\tilde{G}\|_\infty (\# \{k \in \pi \Z,\;  [\beta_+  + k \pi,\beta_+ + k \pi] \cap[\Theta(-a),\Theta(a)] \neq \emptyset  \} + 2)\\
	&\le	\|\tilde{G}\|_\infty \big(\#\{\lambda_{j,\tau} \in [-a,a]\} + 2\big)\;.
	\end{align*}
	By the same arguments as those presented in Subsection \ref{sec:DOS}, the r.v.~on the r.h.s.~has finite expectation. By the Dominated Convergence Theorem we thus deduce that
	\begin{align*}
		\E\bigg[\sum_{j\in\Z} G(\lambda_{j,\tau},\Psi_{j,\tau})\bigg]= \lim_{\varepsilon \downarrow 0} \frac1{2\varepsilon} \E\Big[\int_{\theta \in [\beta_+  - \epsilon,\beta_+  + \epsilon]+\pi \Z} \tilde{G}({\Lambda(\theta)}, \theta_{\Lambda(\theta)},\rho_{\Lambda(\theta)}) d\theta\Big]\;.
	\end{align*}
	Applying a change of variable and Fubini Theorem, we obtain
	\begin{align*}
		&\E\bigg[\sum_{j\in \Z} G(\lambda_{j,\tau},\Psi_{j,\tau})\bigg]\\
		&= \lim_{\varepsilon \downarrow 0} \frac1{2\varepsilon}\E\Big[\int_{\{\lambda \in \R: \theta_\lambda(\tau/2) \in [\beta_+ -\varepsilon,\beta_+ +\varepsilon] + \pi \Z\}} \tilde{G}(\lambda, \theta_\lambda,\rho_\lambda) \partial_\lambda \theta_\lambda(\tau/2) d\lambda\Big]\\
		&= \lim_{\varepsilon \downarrow 0} \frac1{2\varepsilon}\int_{\lambda \in \R} \sum_{n\in \Z} \E\Big[ \mathbf{1}_{\{\theta_\lambda(\tau/2) \in [\beta_+ +n\pi-\varepsilon,\beta_+ +n\pi +\varepsilon]\}}  \tilde{G}(\lambda, \theta_\lambda,\rho_\lambda) \partial_\lambda \theta_\lambda(\tau/2)\Big]\;.
	\end{align*}
	From~\eqref{Eq:SDEPartial}, we deduce that
	\begin{align*}
		\partial_\lambda \theta_\lambda(t) &= \int_{-\tau/2}^{t} e^{-2\rho_\lambda(t)+2\rho_\lambda(s)} ds\;,\quad t\in [-\tau/2,\tau/2]\;,
	\end{align*}
	so that $\partial_\lambda \theta_\lambda(\tau/2) = e^{- 2\rho_\lambda(\tau/2)}\| e^{\rho_\lambda} \|_\tau^2$.\\
	We note that for fixed $\lambda$, the process $\rho_\lambda$ is independent of $\theta_\lambda$. Furthermore, $\theta_\lambda(\tau/2)$ has a Gaussian law with mean $\beta_-+\lambda \tau$ and variance $(\kappa_0^2+\kappa^2)\tau$, and conditionally given $\{\theta_\lambda(\tau/2) = x\}$, $\theta_\lambda$ has the same law as $\beta^{(x)}_\tau$. Putting everything together, we deduce that
	\begin{align*}
		&\E\Big[ \mathbf{1}_{\{\theta_\lambda(\tau/2) \in [\beta_+  +n\pi-\varepsilon,\beta_+  +n\pi +\varepsilon]\}} \tilde{G}(\lambda, \theta_\lambda,\rho_\lambda) \partial_\lambda \theta_\lambda(\tau/2)\Big]\\
		&= \int_{x\in [\beta_+  +n\pi-\varepsilon,\beta_+ +n\pi +\varepsilon]} g_{\lambda \tau, (\kappa_0^2+\kappa^2)\tau}(x-\beta_-) \E\Big[e^{- 2\rho_\lambda(\tau/2)}\| e^{\rho_\lambda} \|_\tau^2 \tilde{G}(\lambda,\beta^{(x)}_\tau,\rho_\lambda)\Big] dx\;.
	\end{align*}
	Applying the Dominated Convergence Theorem on both the integral in $\lambda$ and the sum over $n$, we deduce that
	\begin{align*}
		&\E\bigg[\sum_{j\in\Z} G(\lambda_{j,\tau},\Psi_{j,\tau})\bigg]\\
		&=\int_{\lambda \in \R} \sum_{n\in \Z}\lim_{\varepsilon \downarrow 0} \frac1{2\varepsilon}\int_{x\in [\beta_+ +n\pi-\varepsilon,\beta_+ +n\pi +\varepsilon]} g_{\lambda \tau, (\kappa_0^2+\kappa^2)\tau}(x-\beta_-)\\
		&\qquad\qquad\qquad\qquad \qquad\qquad\times\E\Big[e^{- 2\rho_\lambda(\tau/2)}\| e^{\rho_\lambda} \|_\tau^2 \tilde{G}(\lambda,\beta^{(x)}_\tau,\rho_\lambda)\Big] dx\, d\lambda\\
		&=\int_{\lambda \in \R} \sum_{n\in \Z}g_{\lambda \tau, (\kappa_0^2+\kappa^2)\tau}(\beta_+ - \beta_- +n\pi)  \E\Big[e^{- 2\rho_\lambda(\tau/2)}\| e^{\rho_\lambda} \|_\tau^2 \tilde{G}(\lambda,\beta^{(\beta_+ + n\pi)}_\tau,\rho_\lambda)\Big]d\lambda\;.
	\end{align*}
	Let us set $Y(t) := \rho_\lambda(t)-\rho_\lambda(0)$ and observe that there exists a two-sided Brownian motion $W_1$ such that
	$$ Y(t) = \kappa^2t + \kappa W_1(t)\;.$$
	Applying Girsanov's Theorem as in~\cite[Lemma 2.4]{DL_Critical}, we deduce that for any bounded continuous map $Q:\cC([-\tau/2,\tau/2],\R) \to\R$
	$$ \E[e^{- 2 Y(\tau/2)} Q(Y)] = \E[Q(\hat{Y})]\;,$$
	where $\hat{Y}(t) := -\kappa^2|t|+ \kappa W_1(t)$. Note that $S=\exp(\hat{Y})$ and recall that $\rho_\lambda$ is independent from $\beta^{(\beta_+  +n\pi)}_\tau$. From~\eqref{Eq:GGtilde}, one can check that $\tilde{G}(\lambda,\beta^{(\beta_+  +n\pi)}_\tau,\rho_\lambda) = \tilde{G}(\lambda,\beta^{(\beta_+ +n\pi)}_\tau,\rho_\lambda-\rho_\lambda(0))$. Then a straightforward computation shows that
	\begin{align*}
		\E[e^{- 2\rho_\lambda(\tau/2)}\| e^{\rho_\lambda} \|_\tau^2 \tilde{G}(\lambda,\beta^{(\beta_+ +n\pi)}_\tau,\rho_\lambda)]&=\E[e^{-2Y(\tau/2)}\| e^{Y} \|_\tau^2 \tilde{G}(\lambda,\beta^{(\beta_+ +n\pi)}_\tau,Y)]\\
		&= \E\Big[ \| S\|_{\tau}^2 G\Big(\lambda , \frac{S(\cdot)}{\| S\|_{\tau}} \begin{pmatrix} \sin \beta^{(\beta_+ +n\pi)}_\tau (\cdot) \\\cos \beta^{(\beta_+ +n\pi)}_\tau (\cdot) \end{pmatrix}\Big)\Big]\;,
	\end{align*}
	and the proof is complete.
\end{proof}

\section{Spectral properties of the operator}\label{Sec:Spectral}

\subsection{The matrix spectral measure}

Matrix spectral measures play an important rôle in the spectral study of canonical systems and Dirac operators. They allow one to build explicit spectral representations: for the sake of completeness we will recall the main arguments. Furthermore, we will see that they produce Radon measures from which one can sample eigenvalues, even in situations where the spectrum is dense.

We start by introducing spectral measures in finite volume. Fix $\tau > 0$ and some b.c.~$b_\pm\in\R^2$, or equivalently $\beta_\pm\in[0,\pi)$. We have seen that $\cS_\tau$ has discrete simple spectrum, we thus denote by $(\lambda_{j,\tau},\Psi_{j,\tau})_{j\in\Z}$ its eigenvalues and normalized eigenfunctions. We define the \emph{matrix spectral measure} by
\begin{align*}
	\cM_\tau  = \sum_j \Psi_{j,\tau}(0)   \Psi_{j,\tau}(0) ^\intercal \delta_{\lambda_{j,\tau}}\,.
	\end{align*}
It is a random positive semidefinite matrix-valued measure. Denote its entries by
$$ \cM_\tau := \begin{pmatrix} \xi_\tau & \eta_\tau\\
	\eta_\tau & \zeta_\tau
\end{pmatrix}\;,$$ 
where $\xi_\tau, \zeta_\tau$ are random non-negative Radon measures on $\R$ and $\eta_\tau$ a real random Radon measure that are equal to:
$$ \xi_\tau = \sum_{j} \Psi_{j,\tau}^{(1)}(0)^2 \delta_{\lambda_{j,\tau}}\;,\quad \zeta_\tau = \sum_{j} \Psi_{j,\tau}^{(2)}(0)^2 \delta_{\lambda_{j,\tau}}\;,\quad \eta_\tau = \sum_{j} \Psi_{j,\tau}^{(1)}(0) \Psi_{j,\tau}^{(2)}(0) \delta_{\lambda_{j,\tau}}\;,$$
where $\Psi_{j,\tau}^{(1)}$, resp.~$\Psi_{j,\tau}^{(2)}$, denotes the first, resp.~second, coordinate of the vector $\Psi_{j,\tau}$.\\

The matrix spectral measure $\cM_\tau$ produces an explicit spectral representation of the self-adjoint operator $\cS_\tau$, that we now recall. For any $f\in L^2_{I}((-\tau/2,\tau/2),\R^2)$ let 
$$ F^{N}_\tau(f,\lambda) := \int_{-\tau/2}^{\tau/2} (y_\lambda^{N}(t))^\intercal f(t) dt\;,\quad F^{D}_\tau(f,\lambda) := \int_{-\tau/2}^{\tau/2} (y_\lambda^{D}(t))^\intercal f(t) dt\;,$$
where $y^{N}$ and $y^D$ are the Neumann and Dirichlet solutions introduced in~\eqref{Eq:DirichletNeumann}.
\begin{lemma}[Parseval identity in finite-volume]\label{Lem:Parcevalfinite}
Almost surely,	the map $$F_\tau: f \mapsto (F^N_\tau(f,\cdot), F^D_\tau(f,\cdot))^\intercal$$ is unitary from $L^2_{I}((-\tau/2,\tau/2),\R^2)$ into $L^2_{\cM_\tau}(\R,\R^2)$, and for all $f \in L^2((-\tau/2,\tau/2),\R^2)$
	\begin{equation}\label{Eq:Parseval}
		\|f\|_{L^2((-\tau/2,\tau/2),\R^2)}^2 = \int_{\lambda \in \R} F_\tau(f,\lambda)^\intercal \cM_\tau(d\lambda)F_\tau(f,\lambda)\;.
	\end{equation}
\end{lemma}
\begin{proof}
	Since $\Psi_{j,\tau} = \Psi_{j,\tau}^{(1)}(0) y^N_{\lambda_{j,\tau}} + \Psi_{j,\tau}^{(2)}(0) y^D_{\lambda_{j,\tau}}$, we obtain by linearity
	\begin{align}\label{identityFouriercoeff}
		 \langle \Psi_{j,\tau} , f \rangle_{L^2_{I}((-\tau/2,\tau/2),\R^2)} = \Psi^{(1)}_{j,\tau}(0) F^N(f,\lambda_{j,\tau}) + \Psi^{(2)}_{j,\tau}(0) F^D(f,\lambda_{j,\tau})\;.
	\end{align}
	Therefore
	\begin{align*}
		\langle \Psi_{j,\tau} , f \rangle_{L^2_{I}((-\tau/2,\tau/2),\R^2)}^2 = (F_\tau f(\lambda_{j,\tau}))^\intercal \, \cM_\tau(\{\lambda_{j,\tau}\}) \, F_\tau f(\lambda_{j,\tau})\;.
	\end{align*}
	Summing over $j\in\Z$, the l.h.s.~equals the square of the $L^2_{I}((-\tau/2,\tau/2),\R^2)$-norm of $f$, while the r.h.s.~equals the square of the $L^2_{\cM_\tau}(\R,\R^2)$-norm of $F_\tau f$.
\end{proof}

One would like to construct the matrix spectral measure $\cM_\infty$ of $\cS_\infty$. However its definition is not as straightforward as in finite volume. Indeed the operator may not have genuine eigenvalues (i.e.~the operator might have a continuous spectrum), and even if it had pure point spectrum, this spectrum might not be discrete and it would then be unclear whether the above definitions would make sense in that case.

Fortunately the theory of canonical systems and Dirac operators provides the existence of $\cM_\infty$ as the almost sure vague limit of $\cM_\tau$ as $\tau\uparrow \infty$. We will briefly outline the arguments. However, we will need not only the almost sure but also the $L^p$ vague convergence, and we will provide a complete argument for it.

\begin{theorem}[Convergence of the spectral measures]\label{Th:M}
	Take $\kappa_0,\kappa\ge 0$. The matrix spectral measure $\cM_\tau$ converges vaguely, almost surely and in $L^p$ for any $p\ge 1$, to a matrix-valued measure $\cM_\infty = \begin{pmatrix} \xi_\infty & \eta_\infty\\
		\eta_\infty & \zeta_\infty
	\end{pmatrix}$ which is completely characterized by the Weyl functions $m^{\pm}_z$, $z\in\C\backslash\R$.
\end{theorem}
By convergence in $L^p$, we mean that for any continuous function $h:\R\to\R^2$ with compact support
$$ \E\Big[\Big| \int_\R h(\lambda)^\intercal \cM_\tau(d\lambda) h(\lambda) - \int_\R h(\lambda)^\intercal \cM_\infty(d\lambda) h(\lambda)\Big|^p\Big] \to 0\;.$$
This is equivalent to entrywise vague convergence in $L^p$ by polarization.

The key step for proving the $L^p$ convergence consists in the following probabilistic estimate (note that the restriction to all $\tau > 1$ is arbitrary).
\begin{proposition}\label{Prop:TightnessMtau}
	For any $A>0$ and any $p\ge 1$,
	$$\E\Big[\sup_{\tau > 1} \Big(\xi_\tau([-A,A])^p + \zeta_\tau([-A,A])^p\Big)\Big] < \infty\;.$$
\end{proposition}

\begin{proof}[Proof of Theorem \ref{Th:M}]
	The almost sure convergence of $\cM_\tau$ rests on classical, deterministic arguments from the theory of canonical systems~\cite[Sec.3.7]{Remling} that we briefly outline. For any $z\in\C\backslash\R$, we define
	\begin{equation}\label{Eq:Ntau}
		\cN_\tau(z) := \frac{-1}{m^{+,\tau}_z(\beta_+) + m^{-,\tau}_z(\beta_-)} \begin{pmatrix} 1 & \frac{m^{+,\tau}_z(\beta_+) - m^{-,\tau}_z(\beta_-)}{2} \\ \frac{m^{+,\tau}_z(\beta_+) - m^{-,\tau}_z(\beta_-)}{2} & - m^{+,\tau}_z(\beta_+) m^{-,\tau}_z(\beta_-)\end{pmatrix}\;.
	\end{equation}
	It belongs to the class of Herglotz matrix-valued functions, which display some rigid structure. In particular, $\cN_\tau$ is completely characterized by a triplet: two real symmetric matrices $A_\tau,B_\tau$ and some non-negative matrix-valued measure which happens to be $\cM_\tau$. The characterization takes the explicit form:
	\begin{equation}\label{Eq:Herglotztau}
		\cN_\tau(z) = A_\tau + B_\tau z + \int_{\R} \Big(\frac1{\lambda - z} - \frac{\lambda}{\lambda^2 + 1}\Big) \cM_\tau(d\lambda)\;,
	\end{equation}
On an event of probability $1$, we have the convergence of the Weyl functions $m^{\pm,\tau}_z(\beta_\pm)$ to $m^\pm_z$ as $\tau\uparrow\infty$. We then deduce the convergence of $\cN_\tau(z)$ to $\cN_\infty(z)$, where the latter is the counterpart of \eqref{Eq:Ntau}. One deduces then from this convergence that $\cM_\tau$ converges vaguely to some limit $\cM_\infty$ from the continuity theorem for Herglotz representations (see~\cite[Sec 3.5 and Th 7.3]{Remling} for the scalar version, the matrix version is obtained by projection), and the identity \eqref{Eq:Herglotztau} remains in force for $\tau=\infty$. Furthermore, the matrix spectral measure $\cM_\infty$ is ``explicitly'' determined by $\cN_\infty$ through the following convergence: $\cM_\infty(d\lambda)$ is the a.s.~vague limit as $y\downarrow 0$ of $(1/\pi)\Im(\cN_\infty(\lambda + iy)) d\lambda$.\\
	Proposition \ref{Prop:TightnessMtau} and the bound
	$$ 2 \, | \eta_\tau | \le \xi_\tau + \zeta_\tau\;,$$
	allow one to deduce that the convergence $\cM_\tau\to\cM_\infty$ holds not only almost surely but also in $L^p$, thus concluding the proof.
\end{proof}

The rest of this section is devoted to the proof of Proposition \ref{Prop:TightnessMtau}. Fix $\epsilon >0$ and recall the Dirichlet and Neumann solutions from~\eqref{Eq:DirichletNeumann}, and the vectors $e_1,e_2$ from~\eqref{Eq:e1e2}. For any $A>0$, let $q(A,\epsilon)$ be the largest value $h\in (0,1/2]$ such that for all $\lambda \in [-A,A]$ and all $t\in [0,h]$
$$ | y_\lambda^N(t) - e_1 | \le \epsilon \;,\mbox{ and } | y_\lambda^D(t) - e_2 | \le \epsilon\;.$$
\begin{lemma}\label{Lemma:Moment}
	For any $\epsilon >0$, for any $A>0$ and any $p\ge 1$, $\E[q(A,\epsilon)^{-p}] < \infty$.
\end{lemma}
\begin{proof}
	Recall that $Y_\lambda(0,t)x$, $t\ge 0$ satisfies the SDE \eqref{Eq:FwdSDE}. Recall also the estimates~\eqref{Eq:MomentSDE} and~\eqref{Eq:MomentSDE2}. By Kolmogorov Continuity Theorem~\cite[Theorem I.2.1]{RevuzYor}, there exists $b>0$ such that $\E[Z_A^q] < \infty$ for all $q\ge 1$ where
	$$ Z_A:= \sup_{(\lambda,t,x) \ne (\lambda',t',x')} \frac{|Y_\lambda(0,t)x - Y_{\lambda'}(0,t')x'|}{\big(|t-t'| + |x-x'| + |\lambda-\lambda'|\big)^{b}}\;,$$
	where the supremum is over all $(\lambda,t,x) \ne (\lambda',t',x')$ satisfying $-A \le \lambda, \lambda' \le A$, $0\le t, t' \le 1$ and $|x|,|x'| \le 1$.
	(Recall that the field is already continuous so the modification that Kolmogorov Continuity Theorem constructs necessarily coincides almost surely with the original field). Now observe that on the event $q(A,\epsilon)<1/2$ we have $Z_A \ge \epsilon / q(A,\epsilon)^{b}$ so that
	\begin{align*}
		q(A,\epsilon) &\ge (1/2) \mathbf{1}_{q(A,\epsilon)=1/2} + (\epsilon/Z_A)^{1/b} \mathbf{1}_{q(A,\epsilon)<1/2} \ge (1/2) \wedge (\epsilon/Z_A)^{1/b}\;,
	\end{align*}
	so that
	$$ q(A,\epsilon)^{-1} \le \frac1{(1/2) \wedge (\epsilon / Z_A)^{1/b}} \le 2 + (Z_A / \epsilon)^{1/b}\;,$$
	and the desired result follows.
\end{proof}

\begin{proof}[Proof of Proposition \ref{Prop:TightnessMtau}]
Fix $\epsilon \in (0,1/8)$.	For any $h\in (0,1/2]$, set $g_h(t) := h^{-1} \mathbf{1}_{[0,h]}(t)$. Applying~\eqref{Eq:Parseval} to $f := (g_h, 0)^\intercal$ we obtain (as $\cM_\tau$ is a non-negative matrix-valued measure):
	\begin{align*}
		\| f\|_{L^2((-\tau/2,\tau/2),\R^2)}^2 &= \int_{\lambda \in \R} \Big(F^N_\tau(\lambda)^2 \xi_\tau(d\lambda) + F^D_\tau(\lambda)^2 \zeta_\tau(d\lambda) + 2 F^N_\tau(\lambda)F^D_\tau(\lambda) \eta_\tau(d\lambda)\Big)\\
		&\ge \int_{\lambda \in [-A,A]} \Big(F^N_\tau(\lambda)^2 \xi_\tau(d\lambda) + F^D_\tau(\lambda)^2 \zeta_\tau(d\lambda) + 2 F^N_\tau(\lambda)F^D_\tau(\lambda) \eta_\tau(d\lambda)\Big)\;.
	\end{align*}
	The l.h.s.~equals $h^{-1}$. Note that
	$$ F^N_\tau(\lambda) = \frac1{h} \int_0^h y_\lambda^{N,(1)}(t) dt \;, \quad F^D_\tau(\lambda) = \frac1{h} \int_0^h y_\lambda^{D,(1)}(t) dt\;.$$
	We now take $h=q(A,\epsilon)$ and find
	$$ |F^N_\tau(\lambda) - 1 | \le \epsilon \;,\quad |F^D_\tau(\lambda) | \le \epsilon\;.$$
	Since in addition $2 \, | \eta_\tau| \le \xi_\tau + \zeta_\tau$, we find
	\begin{align*}
		q(A,\epsilon)^{-1} &\ge \int_{\lambda \in [-A,A]} \Big((1-\epsilon)^2 \xi_\tau(d\lambda) - (1+\epsilon)\epsilon |2\eta_\tau(d\lambda)|\Big)\\
		&\ge \int_{\lambda \in [-A,A]} \Big((1-3\epsilon) \xi_\tau(d\lambda) - (1+\epsilon)\epsilon\; \zeta_\tau(d\lambda)\Big)\;.
	\end{align*}
	Applying the same arguments with $f = (0,g_h)^\intercal$, we arrive at
	\begin{align*}
		q(A,\epsilon)^{-1} &\ge \int_{\lambda \in [-A,A]} \Big(-(1+\epsilon)\epsilon\; \xi_\tau(d\lambda) + (1-3\epsilon) \zeta_\tau(d\lambda)\Big)\;.
	\end{align*}
	Summing up these two inequalities, we obtain
	$$ 2q(A,\epsilon)^{-1} \ge (1-4\epsilon-\epsilon^2) \xi_\tau([-A,A]) + (1-4\epsilon-\epsilon^2) \zeta_\tau([-A,A])\;.$$
As $\epsilon$ is small enough, we conclude using Lemma \ref{Lemma:Moment} noticing that $q(A,\epsilon)$ does not depend on $\tau$.
\end{proof}

We now prove that indeed $\cM_\infty$ is a spectral matrix measure for the infinite-volume operator $\cS_{\infty}$.
\begin{proposition}[Spectral representation of $\cS_\infty$]\label{Propo:Spectralmatrix}
	Let us define for compactly supported functions $f \in L^2_{I}(\R,\R^2)$,
	\begin{align*}
		(F_\infty f)(\gl) := \bin{\int_{\R} y_\gl^N(t)^\intercal f(t) dt}{\int_{\R} y_\gl^D(t)^\intercal f(t) dt}
	\end{align*}
	Then for all $f \in L^2_{I}(\R,\R^2)$, the following limit exists in $L^2_{\cM_\infty}(\R,\R^2)$:
	\begin{align*}
		F_\infty f = \lim_{L \to \infty} F_\infty (\un_{(-L,L)}f )\,.
	\end{align*}
	It gives a unitary map from $L^2_{I}(\R,\R^2)$ to $L^2_{\cM_\infty}(\R,\R^2)$. The operator $F_\infty \cS_\infty F_\infty^{-1}$ is a multiplication operator by $\gl$ in the Hilbert space $L^2_{\cM_\infty}(\R,\R^2)$.
\end{proposition}

\begin{proof}
%	Recall that $\cS_\infty$ is obtained from the canonical system by the unitary conjugation
%	\begin{align*}
%		\cS_\infty = Y_0(0,\cdot) \cT_\infty Y_0(0,\cdot)^{-1} \,.
%	\end{align*}
	By the spectral representation theorem for the whole-line canonical systems in the limit-point case, see~\cite[Th 3.21]{Remling}, the matrix measure $\cM_\infty$ is a spectral measure for $\cT_\infty$. More precisely, recall the definition of $U_\gl$ in \eqref{Eq:RJconstant}. Then for all $g \in L^2_K(\R,\R^2)$ compactly supported, define:
	\begin{align*}
		g \mapsto \bin{\int_{\R} (U_\gl(0,t)e_1)^\intercal K(t)g(t) dt}{\int_{\R}(U_\gl(0,t)e_2)^\intercal K(t)g(t) dt} \,.
		\end{align*}
For all $g \in  L^2_K(\R,\R^2)$, the following limit exists in $L^2_{\cM_\infty}(\R,\R^2)$, for $i = 1,2$:
\begin{align*}
\lim_{L \to \infty}	\int_{-L}^L (U_\gl(0,t)e_i)^\intercal K(t)g(t) dt
\end{align*}
and is a unitary map from $ L^2_K(\R,\R^2)$ to $L^2_{\cM_\infty}(\R)$. Moreover, it is proved in~\cite[Th 3.21]{Remling} that this map unitarily transforms $\cT_\infty$ into the multiplication operator by $\gl$ in the space $L^2_{\cM_\infty}(\R)$ (note that Remling uses the terminology of relations in his book, but here the relations are operators as $K(t)$ is positive definite for every $t \in \R$).

It suffices to transfer those facts to $\cS_\infty$ using the unitary map $Y_0(0,\cdot)$ from $L^2_{K}(\R,\R^2)$ to $L^2_{I}(\R,\R^2)$. Indeed, for $f \in  L^2_{I}(\R,\R^2)$, we have for $g := Y_0(0,\cdot)^{-1} f$:
\begin{align*}
	(U_\gl(t)e_1)^\intercal K(t)g(t) &= e_1^\intercal U_\gl(0,t)^\intercal Y_0(0,\cdot)^{\intercal} f(t)  \\
	&= e_1^\intercal Y_\gl(0,t)^\intercal  f(t)\\
	&= y_\gl^N(t)^\intercal  f(t)
\end{align*}
where we used the identity of Lemma \ref{Lemma:Flowvz} for the second line. Similarly,
\begin{align*}
		(U_\gl(t)e_2)^\intercal K(t)g(t)  = y_\gl^D(t)^\intercal  f(t) \,,
\end{align*}
which implies the result.
\end{proof}

\subsection{Anderson localization}

In this subsection, we always assume that $\kappa >0$. Our proof of Theorem \ref{Th:Localization} relies on the following crucial estimate, where $\sigma_\infty := \xi_\infty + \zeta_\infty$ is the trace of the matrix spectral measure $\cM_\infty$.

\begin{proposition}\label{Prop:LocEstimate}
	Take $\kappa >0$. For any bounded interval $E \subset \R$, any $\epsilon > 0$ and any $p \in (0, 2\epsilon/\kappa^2)$, we have
	\begin{equation}
		\E\bigg[\int_E \inf_{\{x \in \R^2:|x|=1\}} \sup_{t\in\R} \Big(|Y_\lambda(0,t)x| e^{(\kappa^2-\epsilon)|t|} + \frac1{|Y_\lambda(0,t)x|} e^{-(\kappa^2+\epsilon)|t|} \Big)^p \sigma_\infty(d\lambda)\bigg] < \infty\;.
	\end{equation}
\end{proposition}
\begin{proof}
	The proof follows the same line as~\cite[Section 3.3]{CarmonaDuke} and~\cite[Section 4.4]{DL_Loc}. We can assume w.l.o.g.~that $E$ is an open bounded interval. Set for any $A>1$
	$$ G_A(\lambda) := \inf_{\{x \in \R^2:|x|=1\}} \sup_{t\in (-A/2,A/2)} \Big(|Y_\lambda(0,t)x| e^{(\kappa^2-\epsilon)|t|} + \frac1{|Y_\lambda(0,t)x|} e^{-(\kappa^2+\epsilon)|t|} \Big)^p\;,\quad \lambda \in\R\;.$$
	By the Monotone Convergence Theorem it suffices to prove that
	$$ \sup_{A>0} \E\Big[\int_{\lambda \in E} G_A(\lambda) \sigma_\infty(d\lambda)\Big] < \infty\;.$$
	By Theorem \ref{Th:M}, $\sigma_\tau := \xi_\tau+\zeta_\tau$ converges a.s.~vaguely to $\sigma_\infty$. For convenience, we endow the operator $\cS_{\tau}$ with the b.c.~$b_\pm=e_2$, or equivalently $\beta_\pm = 0$. For any $A>1$, the continuity properties of the stochastic flow ensure that $\lambda\mapsto G_A(\lambda)$ is continuous and therefore $G_A \sigma_\tau$ converges a.s.~vaguely to $G_A \sigma_\infty$ as $\tau \to\infty$. Since $E$ is open and bounded, we deduce that
	$$ \int_E G_A(\lambda) \sigma_\infty(d\lambda) \le \liminf_{\tau\to\infty} \int_E G_A(\lambda) \sigma_\tau(d\lambda)\;.$$
	Applying Fatou's Lemma, we deduce that we only have to prove that
	$$ \sup_{A>1} \sup_{\tau/2 > A} \E\Big[\int_{\lambda \in E} G_A(\lambda) \sigma_\tau(d\lambda)\Big] < \infty\;.$$
	Now observe that
	\begin{align*}
		\int_{\lambda \in E} G_A(\lambda) \sigma_\tau(d\lambda) &= \sum_{j\in \Z} \Big| \Psi_{j,\tau}(0)\Big|^2 \mathbf{1}_{\lambda_{j,\tau} \in E} \,G_A(\lambda_{j,\tau})\\
		&\le \sum_{j \in \Z} \Big| \Psi_{j,\tau}(0)\Big|^2 H_A(\lambda_{j,\tau}, \big| \Psi_{j,\tau} \big|)\;,
	\end{align*}
	where for any function $f:\R\to (0,\infty)$
	$$ H_A(\lambda,f) := \mathbf{1}_{E}(\lambda) \sup_{t\in (-A/2,A/2)} \Big(\frac{f(t)}{f(0)} e^{(\kappa^2-\epsilon)|t|} + \frac{f(0)}{f(t)} e^{-(\kappa^2+\epsilon)|t|} \Big)^p\;.$$
	The GMP formula stated in Proposition \ref{Prop:GMP2} yields
	\begin{align*}
		\E\bigg[\sum_{j\in\Z} \Big| \Psi_{j,\tau}(0)\Big|^2 H_A(\lambda_{j,\tau}, \big| \Psi_{j,\tau} \big|) \bigg] = \int_{\R}\sum_{n\in\Z}  \E\bigg[ H_A\bigg(\lambda , \frac{S(\cdot)}{\| S\|_{\tau}} \bigg)\bigg] g_{\lambda \tau,(\kappa_0^2+\kappa^2)\tau}(n\pi)d\lambda\;,
	\end{align*}
	where $S(t) = \exp(\kappa W_1(t) - \kappa^2|t|)$. At this point, we note that
	\begin{align*}
		H_A\bigg(\lambda , \frac{S(\cdot)}{\| S\|_{\tau}}\Big) &= \mathbf{1}_{\lambda \in E} \sup_{t\in (-A/2,A/2)} \Big( S(t)e^{(\kappa^2-\epsilon)|t|} + \frac1{S(t)} e^{-(\kappa^2+\epsilon)|t|} \Big)^p\\
		&= \mathbf{1}_{\lambda \in E} \sup_{t\in (-A/2,A/2)} \Big( e^{\kappa W_1(t) -\epsilon|t|} + e^{- \kappa W_1(t) -\epsilon|t|} \Big)^p\;.
	\end{align*}
	Basic properties of Brownian motion ensure the existence of a constant $C_{\epsilon,p} >0$ such that,
	$$ \sup_{A>1} \sup_{\tau/2 > A} \E\bigg[ H_A\bigg(\lambda , \frac{S(\cdot)}{\| S\|_{\tau}} \bigg)\bigg] < C_{\epsilon,p} \mathbf{1}_{E}(\lambda)\;,\quad \forall \lambda \in \R\;,$$
	if $p \in (0,2\epsilon/\kappa^2)$
	and therefore
	$$ \sup_{A>1} \sup_{\tau/2 > A} \E\bigg[\sum_{j\in\Z} \Big| \Psi_{j,\tau}(0)\Big|^2 H_A(\lambda_{j,\tau}, \Psi_{j,\tau}) \bigg]  \le C_{\epsilon,p} \sup_{\tau > 2} \int_E \sum_{n\in\Z}g_{\lambda \tau,(\kappa_0^2+\kappa^2)\tau}(n\pi) d\lambda <\infty\;,$$
	using Lemma \ref{Lemma:n0} for the last bound.
\end{proof}
We can now proceed with the proof of Theorem \ref{Th:Localization}, except the identification of the spectrum of $\cS_{\infty}$ that will be carried out at the end of the article.
\begin{proof}[Proof of Theorem \ref{Th:Localization} - first part]
	Let $\epsilon \in (0,\kappa^2)$ and let $\cA_\epsilon$ be the set of $\lambda \in \R$ for which there exists $x\in \R^2$, with $|x|=1$, such that
	$$ \sup_{t\in\R} \Big(|Y_\lambda(0,t)x| e^{(\kappa^2-\epsilon)|t|} + \frac1{|Y_\lambda(0,t)x|} e^{-(\kappa^2+\epsilon)|t|} \Big) < \infty\;.$$
	By Proposition \ref{Prop:CaracEigenvalueTinfty}, any such $\lambda$ is an eigenvalue of $\cS_\infty$. As a consequence, almost surely $\cA_\epsilon$ is a countable set.\\
	Proposition \ref{Prop:LocEstimate} ensures that almost surely $\sigma(\R\backslash \cA_\epsilon) = 0$, and therefore $\sigma$ is almost surely an atomic measure. Using that for any Borel set $B$, $\sigma_\infty(B) = 0$ implies $\cM_\infty(B) = 0$, and the fact that $\cM_\infty$ is a spectral matrix measure by Proposition \ref{Propo:Spectralmatrix}, we deduce that $\cS_\infty$ is pure point. Note that the multiplicity of each eigenvalue is necessarily one (since otherwise there would exist two linearly independent integrable solutions of the eigenproblem and the operator would not be limit point). The bound stated in Proposition \ref{Prop:LocEstimate} ensures the exponential localization stated in \eqref{Eq:ExpoLoc}.
\end{proof}

\section{Convergence of the eigenfunctions and intensity}\label{Sec:CVeigenfunctions}

\subsection{Convergence of the eigenpair point measure}

We assume in this section that $\kappa > 0$. Recall the Banach space $\cF$ introduced in~\eqref{Eq:DefF}, and its quotient $\tilde{\cF}$ by the equivalence relation corresponding to sign change. By Theorem \ref{Th:Localization} - first part, a.s~every eigenfunction $\Psi_{j,\infty}$ belongs to $\cF$. On the other hand, a.s.~every eigenfunction $\Psi_{j,\tau}$ is a continuous function on $[-\tau/2,\tau/2]$ ; we extend them into elements of $\cF$ by setting their values to zero outside $[-\tau/2-1,\tau/2+1]$ and by linear interpolation in between. This choice of extension is arbitrary and has essentially no impact on the results to come. From now on, we view the eigenfunctions as elements of $\tilde{\cF}$.\\

For any $\tau > 0$, we define the random measures
$$ \mu_{\tau} := \sum_{j\in\Z}  |\Psi_{j,\tau}(0)|^2 \delta_{(\lambda_{j,\tau},\Psi_{j,\tau})}\;,\quad \mu_{\infty} := \sum_{j\in\Z}  |\Psi_{j,\infty}(0)|^2 \delta_{(\lambda_{j,\infty},\Psi_{j,\infty})}\;,$$
that take values in the space $\ccM$ of non-negative measures on $\R \times \tilde{\cF}$ which are finite on $[a,b]\times\tilde{\cF}$ for all $a<b$, endowed with the smallest topology that makes the maps
\begin{align*}
	\mu \mapsto \int g(\gl,\psi) d\mu(\gl,\psi)
\end{align*}
continuous for all bounded continuous $g:\R \times \tilde{\cF} \to \R$ compactly supported in the first coordinate. By a slight adaptation of standard proofs, see~\cite[Theorem A2.6.III]{DaleyVereJones}, $\ccM$ is a Polish space.

 We now collect two estimates on $\mu_\tau$.
\begin{lemma}\label{Lemma:Tight}
	For any $C>0$, there exists a relatively compact subset $\tilde{K}(C)$ of $\tilde{\cF}$ such that, for any $[a,b]\subset \R$
	\begin{equation}\label{Eq:TightKC}
		\limsup_{C\to\infty} \limsup_{\tau \to \infty} \E\Big[\mu_{\tau}\big([a,b]\times \tilde{\cF}\backslash \tilde{K}(C)\big)\Big] = 0\;.
	\end{equation}
\end{lemma}
The proof uses the GMP formula of Proposition \ref{Prop:GMP2}. As it is quite technical, it is deferred to Appendix \ref{Appendix:Compact}.

The next lemma shows that, uniformly over all $\tau > 1$, the measure $\mu_\tau$ ``mostly'' charges pairs $(\lambda,\Psi)$ such that $|\Psi(0)|$ is not small.

\begin{lemma}\label{Lemma:ZeroMass}
	For any $[a,b]\subset \R$, we have
	$$ \lim_{\epsilon\downarrow 0} \sup_{\tau \geq 1}\E\Big[\mu_{\tau}\big([a,b]\times \{\Psi: |\Psi(0)|^2 \le \epsilon\}\big)\Big] = 0\;.$$
\end{lemma}
\begin{proof}
	By Proposition \ref{Prop:GMP2}, the expectation of the statement equals
	\begin{align*}
		&\int_{\lambda\in [a,b]} \sum_{n\in\Z}  \P\Big( \frac{|S(0)|^2}{\|S\|^2_\tau} \le \epsilon\Big) g_{\lambda \tau,(\kappa_0^2+\kappa^2)\tau}(\beta_+-\beta_-+n\pi) d\lambda\;.
	\end{align*}
	Note that
	$$ \P\Big( \frac{|S(0)|^2}{\|S\|_\tau^2} \le \epsilon\Big) \le \P\Big( \frac{1}{\int_{\R} e^{2\kappa W_1(t) - 2\kappa^2 |t|} dt} \le \epsilon\Big)\;.$$
	This last quantity does not depend on $\tau$ and goes to $0$ as $\epsilon\downarrow 0$ (as $\int e^{2\kappa W_1(t) - 2\kappa^2 |t|} dt$ is finite a.s.). Applying Lemma \ref{Lemma:n0} we can conclude.
\end{proof}
	Before we proceed with the proof of the theorem, let us set $\sigma_\tau := \xi_\tau+\zeta_\tau$ and $\sigma_\infty := \xi_\infty+\zeta_\infty$. These two measures are the marginals of $\mu_\tau$ and $\mu_\infty$
	\begin{equation}\label{Eq:sigmamu}
		\sigma_\tau(d\lambda) = \int_{\Psi} \mu_{\tau}(d\lambda,d\Psi)\;,\quad \sigma_\infty(d\lambda) = \int_{\Psi} \mu_{\infty}(d\lambda,d\Psi)\;.
	\end{equation}
	By Theorem \ref{Th:M}, we already know that $\sigma_\tau$ converges a.s.~vaguely to $\sigma_\infty$ as $\tau\to\infty$. 
\begin{proof}[Proof of Theorem \ref{Th:CVrandommeas}]
	In the sequel, we view $Y := (Y_z(s,t);(s,t) \in \R^2,z\in\C)$ as a r.v.~with values in $\cC(\R^2\times\C,\C^4)$.\\
	Our strategy of proof is slightly unusual: it consists in proving that the triplet $X_\tau:=(\mu_\tau,\mu_\infty,Y)$ converges \textit{in law} to $X_\infty:=(\mu_\infty,\mu_\infty,Y)$. From the convergence of the first two marginals, we deduce that $d_{\ccM}(\mu_\tau,\mu_\infty)$ converges in law to $0$, which implies the convergence in probability.\\

	\textit{Step 1: tightness.} Fix some sequence $(\tau_k)_k$ converging to $\infty$. We claim that the collection of r.v.~$(X_{\tau_k})_k := (\mu_{\tau_k},\mu_\infty,Y)_{k}$ is tight in the Polish space $\ccM\times\ccM\times \cC(\R^2\times\C,\C^4)$. Since the second and third coordinates do not depend on $\tau_k$, we only have to prove that $(\mu_{\tau_k})_{k}$ is tight in $\ccM$. To that end, Prohorov’s Theorem states that it suffices to show that for any $[a,b]\subset \R$ and any $\delta > 0$, there exists a compact set $\tilde{K}\subset \tilde{\cF}$ such that
	\begin{equation}\label{Eq:mutauTight}
		\limsup_{k\to\infty} \E[\mu_{\tau_k}([a,b]\times (\tilde{\cF}\backslash \tilde{K}))] \le \delta\;.
	\end{equation}
and
	\begin{align*}
		\limsup_{k\to\infty} \E[\mu_{\tau_k}( [a,b]\times \tilde{\cF})] < \infty\;.
	\end{align*}
	Lemma \ref{Lemma:Tight} provides the existence of such a compact $\tilde{K}$ for the first inequality, while Proposition \ref{Prop:TightnessMtau} combined with~\eqref{Eq:sigmamu} implies the second statement.\\
	Let $\hat{X}$ be the limit in law of a converging subsequence, for notational convenience we still denote by $(X_{\tau_k})_{k \ge 1}$ this subsequence. We aim at showing that $\hat{X}$ has the same law as $X_\infty$. To ease the presentation, we use Skorokhod Representation Theorem: on some probability space, there exists a sequence $(\mu_{\tau_k}^{(k)},\mu_\infty^{(k)},Y^{(k)})_{k\ge 1}$ that converges almost surely to a triplet $(\hat{\mu},\mu_\infty^{(\infty)},Y^{(\infty)})$, where each $(\mu_{\tau_k}^{(k)},\mu_\infty^{(k)},Y^{(k)})$ equals in law $X_{\tau_k}$ and $(\hat{\mu},\mu_\infty^{(\infty)},Y^{(\infty)})$ equals in law $\hat{X}$.\\
	Observe that $(\sigma_{\tau_k}^{(k)},\sigma_\infty^{(k)})$ converges almost surely to $(\hat{\sigma},\sigma_\infty^{(\infty)})$ where $\hat{\sigma}(d\lambda) := \int_\Psi \hat{\mu}(d\lambda,d\Psi)$. Since we already know that $(\sigma_{\tau_k}^{(k)},\sigma_\infty^{(k)})$ converges in law to $(\sigma_\infty,\sigma_\infty)$, we deduce that $\hat{\sigma} = \sigma_\infty^{(\infty)}$ almost surely.\\
	
	\textit{Step 2: the limit is atomic.} Let us prove that $\hat{\mu}$ is atomic. To that end, we will first disregard the atoms of $\mu_{\tau_k}^{(k)}$ which are ``small'' and then, we will use Lemma \ref{Lemma:ZeroMass} to conclude. Let $\chi : \R_+ \to [0,1]$ be an arbitrary non-decreasing continuous function that equals $0$ on $[0,1/2)$ and $1$ on $[1,\infty)$. For any $\epsilon > 0$, we let $\chi_\epsilon(x) = \chi(x/\epsilon)$ such that it equals $0$ on $[0,\epsilon/2)$ and $1$ on $[\epsilon,\infty)$. Notice that $\chi_\epsilon \uparrow 1$ when $\epsilon \downarrow 0$. 
	We set
	$$ \hat{\mu}_{\epsilon}(d\lambda,d\Psi) := \chi_\epsilon(|\Psi(0)|^2) \hat{\mu}(d\lambda,d\Psi)\;,$$
	as well as
	$$ \mu_{\tau_k,\epsilon}^{(k)}(d\lambda,d\Psi) := \chi_\epsilon(|\Psi(0)|^2) \mu_{\tau_k}^{(k)}(d\lambda,d\Psi)\;.$$	
	Since $\Psi\mapsto \chi_\epsilon(|\Psi(0)|^2)$ is a bounded continuous function on $\tilde{\cF}$, we deduce that almost surely $\mu_{\tau_k,\epsilon}^{(k)} \to \hat{\mu}_{\epsilon}$ as $k\to\infty$.\\
	We start by proving that $\hat{\mu}_{\epsilon}$ is atomic almost surely. Recall that this is the case for every measure $\mu_{\tau_k,\epsilon}^{(k)}$ by definition. Fix $[a,b] \subset \R$. The number of atoms of $\mu_{\tau_k,\epsilon}^{(k)}$ in $[a,b]\times \tilde{\cF}$ is bounded by the number of atoms of the measure
	$$ \mathbf{1}_{[\epsilon/2,\infty)}(|\Psi(0)|^2) \mu_{\tau_k}^{(k)}(d\lambda,d\Psi)$$
	in $[a,b]\times\tilde{\cF}$, which is itself bounded by $(2/\epsilon) \sigma_{\tau_k}^{(k)}([a,b])$. The almost sure convergence of $(\sigma_{\tau_k}^{(k)})_{k\ge 1}$ ensures that $m:= \sup_{k\ge 1} \sigma_{\tau_k}^{(k)}([a,b])$ is a.s.~finite. Consequently the number of atoms of $\mu_{\tau_k,\epsilon}^{(k)}$ in $[a,b]\times \tilde{\cF}$ is bounded uniformly over all $k \ge 1$. We thus deduce that the restriction to $(a,b)$ of $\hat{\mu}_{\epsilon}$ is atomic almost surely. Since $[a,b]$ was arbitrary, this suffices to deduce that $\hat{\mu}_{\epsilon}$ is atomic almost surely.\\

	Now we observe that $\hat{\mu}_{\epsilon}$ is non-decreasing as $\epsilon \downarrow 0$, and is bounded from above by $\hat{\mu}$. 
Let us fix a bounded interval $[a,b]$. Let $f_{a,b} :\R\to [0,1]$ be a continuous function supported in $[a-1,b+1]$ that equals one on $[a,b]$. We have
	\begin{align*}
		\int_{\gl, \Psi} f_{a,b}(\gl) d(\hat \mu(\gl, \Psi) - \hat \mu_\epsilon(\gl, \Psi)) = \lim_{k \to \infty}	\int_{\gl, \Psi} f_{a,b}(\gl) d( \mu^{(k)}_{\tau_k}(\gl, \Psi) -  \mu^{(k)}_{\tau_k,\epsilon}(\gl, \Psi)) \,.
	\end{align*}
	Therefore, using Fatou's Lemma,
		\begin{align*}
			\E\Big[ 	\int_{\gl, \Psi} f_{a,b}(\gl) d(\hat \mu(\gl, \Psi) - \hat \mu_\epsilon(\gl, \Psi))\Big] &\leq \liminf_{k \to \infty}\E\Big[ \int_{\gl, \Psi} f_{a,b}(\gl) d( \mu^{(k)}_{\tau_k}(\gl, \Psi) -  \mu^{(k)}_{\tau_k,\epsilon}(\gl, \Psi)) \Big] \\
			&\leq   \sup_{k \geq 1} 	\E\Big[ \int_{\gl, \Psi} \un_{[a-1,b+1]}(\gl) d( \mu^{(k)}_{\tau_k}(\gl, \Psi) -  \mu^{(k)}_{\tau_k,\epsilon}(\gl, \Psi)) \Big]\,.
		\end{align*}
	But thanks to Lemma \ref{Lemma:ZeroMass}, this last quantity tends to $0$ when $\epsilon \downarrow 0$. It implies 
		\begin{align*}
		\E\Big[ 	\hat \mu([a,b] \times \tilde{\cF}) - \hat \mu_\epsilon([a,b] \times \tilde{\cF})\Big] \to 0\,.
	\end{align*}
	As $\hat \mu([a,b] \times \tilde{\cF})- \hat \mu_\epsilon([a,b] \times \tilde{\cF}) \geq 0$ is non-increasing when $\epsilon \downarrow 0$, by the Monotone Convergence Theorem, we obtain that almost surely,
	\begin{align}\label{CVtildemu}
			\lim_{\epsilon \downarrow 0}  \hat \mu_\epsilon([a,b] \times \tilde{\cF}) = \hat \mu([a,b] \times \tilde{\cF})\,.
	\end{align}
	Let us take a countable sequence $\eps_j \downarrow 0$ and introduce the countable set $$A := \bigcup_j \{\mbox{Atoms of } \hat{\mu}_{\epsilon_j}|_{[a,b] \times \tilde{\cF}} \} \,.$$ We have
	\begin{align*}
	0 \leq 	\hat \mu(([a,b] \times \tilde{\cF}) \backslash A) - 	\hat \mu_{\epsilon_j}(([a,b] \times \tilde{\cF})\backslash A) \leq \hat \mu([a,b] \times \tilde{\cF}) -\hat \mu_{\epsilon_j}([a,b] \times \tilde{\cF}) \to 0
	\end{align*}
as $\hat \mu_{\epsilon_j}(([a,b] \times \tilde{\cF})\backslash A) = 0$, we deduce $	\hat \mu(([a,b] \times \tilde{\cF}) \backslash A)  =0$.	Since $[a,b]$ is arbitrary, this completes the proof of the atomicity of $\hat{\mu}$.\\

	\textit{Step 3: identification of the atoms.} Let $(\lambda,\Psi)$ be an atom of $\hat{\mu}$. We aim at showing that: 
	\begin{enumerate}
		\item[(1)] $\hat{\mu}(\{(\lambda,\Psi)\}) = |\Psi(0)|^2$ ;
		\item[(2)] $(\lambda,\Psi)$ is an atom of $\mu_\infty^{(\infty)}$.
	\end{enumerate}
From \eqref{CVtildemu}, we have $\hat \mu(\{(\gl,\Psi) \,:\, \Psi(0) = 0\}) = 0$. Therefore $|\Psi(0)|^2 > 0$. There must exist a sequence $(\lambda_k,\Psi_k)_{k\ge 1}$ of atoms of $\mu_{\tau_k}^{(k)}(d\lambda,d\Psi)$ that converges in $\R\times\tilde{\cF}$ towards $(\lambda,\Psi)$. The topology of $\tilde{\cF}$ ensures that $|\Psi_k(0)|^2 > 0$ for all $k$ large enough. Consequently for all $k$ large enough $(\lambda_k,\Psi_k)$ is an atom of mass $|\Psi_k(0)|^2$ of $\mu_{\tau_k}^{(k)}(d\lambda,d\Psi)$. Denoting by $\bar{B}((x,f),r)$ the closed ball of center $(x,f) \in \R\times \tilde{\cF}$ with radius $r>0$, the a.s.~convergence of $\mu_{\tau_k}^{(k)}$ to $\hat{\mu}$ ensures that
\begin{align*}
	\hat{\mu}\big(\bar{B}((\lambda,\Psi),r)\big) &\ge \limsup_{k\to\infty} \mu_{\tau_k}^{(k)}\big(\bar{B}((\lambda,\Psi),r)\big)\\
	&\ge \limsup_{k\to\infty}\mu^{(k)}_{\tau_k}(\{(\lambda_k,\Psi_k)\}) = \limsup_{k\to\infty} |\Psi_k(0)|^2 = |\Psi(0)|^2\;.
\end{align*}
Passing to the limit $r\downarrow 0$, we deduce that
	\begin{equation}\label{Eq:LBmass}
		\hat{\mu}(\{(\lambda,\Psi)\}) \ge |\Psi(0)|^2\;.
	\end{equation}
	Before proving the converse inequality, let us prove that $(\lambda,\Psi)$ is an atom of $\mu_\infty^{(\infty)}$. Since $(\lambda_k,\Psi_k)$ is an eigenpair of $\cS_{\tau_k}^{(k)}$ we have
	$$ \Psi_k(t) = Y_{\lambda_k}^{(k)}(0,t) \Psi_k(0)\;,\quad t\in [-\tau_k/2,\tau_k/2]\;.$$
	Since $(\lambda_k,\Psi_k)$ converges to $(\lambda,\Psi)$ in $\R\times\tilde{\cF}$ we deduce that
	$$ \Psi(t) = Y_\lambda^{(\infty)}(0,t) \Psi(0)\;,\quad t\in \R\;.$$
	Recall that $\Psi_k$ is the extension to $\R$ of a normalized eigenfunction on $[-\tau_k/2,\tau_k/2]$. A simple computation yields
	$$ \|\Psi_k\|_{L^2}^2 = 1 + \frac13\Big(|\Psi_k(-\tau_k/2)|^2 + |\Psi_k(\tau_k/2)|^2 \Big)\;.$$
	The topology on $\tilde{\cF}$ ensures that $\|\Psi\|_{L^2} = \lim_{k\to\infty} \|\Psi_k\|_{L^2} = 1$. In turn, Proposition \ref{Prop:CaracEigenvalueTinfty} ensures that $(\lambda,\Psi)$ is an atom of $\mu_\infty^{(\infty)}$.\\
	Now, since $\hat{\sigma} = \sigma_\infty^{(\infty)}$, we obtain
	$$ |\Psi(0)|^2 = \sigma_\infty^{(\infty)}(\{\lambda\}) = \hat{\sigma}(\{\lambda\}) \ge \hat{\mu}(\{(\lambda,\Psi)\})\;,$$
	and, combined with~\eqref{Eq:LBmass}, we deduce that $\hat{\mu}(\{(\lambda,\Psi)\}) = |\Psi(0)|^2$.\\
	
	\textit{Step 4: identification of the limit.} We have thus proven that almost surely
	$$ \mu_\infty^{(\infty)} \ge \hat{\mu}\;,\quad \sigma_\infty^{(\infty)} = \hat{\sigma}\;.$$
	Given the relationships~\eqref{Eq:sigmamu} between the $\mu$'s and the $\sigma$'s, necessarily $\mu_\infty^{(\infty)} = \hat{\mu}$ and this completes the proof of the convergence in probability.
	Let us now prove the $L^p$ convergence for any given $p \geq 1$. Since we already have convergence in probability, it suffices to prove uniform integrability of the $p$-th moment. Fix $q>p$. For all $F \;:\; [a,b] \times \tilde{\cF} \to \R$ bounded continuous, we have
	\begin{align*}
	\Big| \int F d \mu_\tau -  \int F d \mu_\infty \Big|^q &\leq \|F\|^q_\infty (\sigma_\tau([a,b]) + \sigma_\infty([a,b]))^q \\
		&\leq C_q \|F\|^q_\infty (\sigma_\tau([a,b])^q + \sigma_\infty([a,b])^q )\,.
	\end{align*}
	Thanks to Proposition \ref{Prop:TightnessMtau}, it implies
	\begin{align*}
		\sup_\tau \E\Big[	\Big| \int F d \mu_\tau -  \int F d \mu_\infty \Big|^q \Big] < \infty\,.
	\end{align*}
\end{proof}
\begin{remark}
The preceding proof gives an alternative route to prove pure point
	spectrum. Indeed, one could consider instead of the triplet $(\mu_\tau, \mu_\infty,Y)$ only the first coordinate $\mu_\tau$ (the other coordinates are used only to identify the atoms and the proof of its atomicity did not use them). 
	Since the projection of $\mu_\tau$ onto the spectral coordinate is
	$\sigma_\tau$, and since $\sigma_\tau$ converges vaguely to
	$\sigma_\infty$, one obtains $\hat\mu (\cdot \times \tilde{\cF})\overset{(d)}{=}
	\sigma_\infty(\cdot)$. Hence the atomicity of $\hat\mu$ implies the
one of $\sigma_\infty$. 
\end{remark}

\subsection{Convergence of the eigenpairs}\label{Subsec:CVeigenfunction}

We now turn to the proof of Corollary \ref{Cor:CVeigenfunctions}. As explained in the introduction, since the limiting eigenvalues are a.s.~dense in $\R$ there is no obvious way to order them. We will rely on the weight function
\begin{align*}
	w(\gl,\Psi) := \frac{|\Psi(0)|^2}{1+ |\gl|^2}\;,
\end{align*}
which ``favors'' those eigenpairs $(\lambda,\Psi)$ s.t.~$\lambda$ is close to $0$ and $|\Psi(0)|$ is large. From the first spectral properties we established above, $w$ is almost surely integrable w.r.t. $\sum_j \delta_{\gl_{j,\infty}, \Psi_{j,\infty}}$. For all $\epsilon >0$, this implies that the number of atoms $(\gl,\Psi)$ such that $w(\gl,\Psi) \geq \epsilon$ is finite. We can therefore order the atoms in decreasing order of the weights $w$. We then prove that the convergence in probability of $\mu_\tau$ towards $\mu_\infty$ implies the convergence of individual atoms after this choice of enumeration.

\begin{proof}[Proof of Corollary \ref{Cor:CVeigenfunctions}]
Recall the characterization of $\cM_\infty$ in \eqref{Eq:Herglotztau}, valid for $\tau = \infty$ (see the discussion below \eqref{Eq:Herglotztau}). Taking $z=i$ and the imaginary part, we get
\begin{align*}
	\Im \cN_\infty(i) = B_\infty + \int_{\R} \frac{1}{1 + \gl^2} \cM_\infty(d\gl)\,.
\end{align*}	
By applying the trace, we obtain
	$$ \sum_{j\in\Z} \frac1{1+\lambda_{j,\infty}^2}|\Psi_{j,\infty}(0)|^2 < \infty\;.$$
	We now introduce a convenient enumeration of the eigenvalues/eigenfunctions of $\cS_{\infty}$. 
	We denote by $(\nu_{m,\infty}, \Phi_{m,\infty})_{m\ge 1}$ the eigenvalues and eigenfunctions $(\lambda_{j,\infty}, \Psi_{j,\infty})_{j\in\Z}$ arranged in such a way that for any $m\ge 1$:\begin{enumerate}
		\item either $\frac1{1+\nu_{m,\infty}^2}|\Phi_{m,\infty}(0)|^2 > \frac1{1+\nu_{m+1,\infty}^2}|\Phi_{m+1,\infty}(0)|^2$
		\item or $\frac1{1+\nu_{m,\infty}^2}|\Phi_{m,\infty}(0)|^2 = \frac1{1+\nu_{m+1,\infty}^2}|\Phi_{m+1,\infty}(0)|^2$, and $\nu_{m,\infty} < \nu_{m+1,\infty}$.
	\end{enumerate}
	Now consider the distance
	$$ d((\lambda,\Psi);(\nu,\Phi)) := |\lambda-\nu| + \tilde{d}(\Psi , \Phi)\;.$$
	
	For any $\tau > 0$, we define $(\nu_{m,\tau},\Phi_{m,\tau})_{m\ge 1}$ according to the following algorithm. For any $m\ge 1$, the pair $(\nu_{m,\tau},\Phi_{m,\tau})$ minimizes the distance $d$ between $(\nu_{m,\infty}, \Phi_{m,\infty})$ and the finite-volume eigenpairs which have not yet been selected $$(\lambda_{m,\tau},\Psi_{m,\tau})_{m\in \Z} \backslash (\nu_{i,\tau},\Phi_{i,\tau})_{1 \le i \le m-1}\;.$$
	If this minimum is achieved by several pairs, pick the pair with the smallest eigenvalue.
	
	We stop the algorithm at $j = N_\tau := \lfloor \tau \rfloor$, and then enumerate in an arbitrary order all the remaining pairs of eigenvalues / eigenfunctions of $\cS_\tau$. This gives a complete enumeration of the finite-volume eigenpairs.

	Let us prove by contradiction that for all $M \geq 1$, $(\nu_{m,\tau}, \Phi_{m,\tau})_{1 \le m\le M}$ converges in probability to $(\nu_{m,\infty}, \Phi_{m,\infty})_{1 \le m\le M}$. We work by induction. Assume that there exists $j_0\ge 1$ such that $(\nu_{i,\tau}, \Phi_{i,\tau})_{1\le i \le j_0-1}$ converges in probability to $(\nu_{i,\infty}, \Phi_{i,\infty})_{1\le i \le j_0-1}$ but $(\nu_{j_0,\tau}, \Phi_{j_0,\tau})$ does not converge in probability to $(\nu_{j_0,\infty}, \Phi_{j_0,\infty})$. Contradicting the convergence in probability, we deduce that there exist $\delta_0>0$, $\rho_0 > 0$ and a sequence $\tau_k$ such that for all $k$, $\P\big( A_k\big) > \delta_0$
	where $$A_k := \Big\{d((\nu_{j_0,\tau_k}, \Phi_{j_0,\tau_k}) ; (\nu_{j_0,\infty}, \Phi_{j_0,\infty})) > \rho_0 \Big\} \,.$$
	Furthermore, there exists $\epsilon> 0$ such that the first $j_0$ points $(\nu_{i,\infty}, \Phi_{i,\infty})_{1\le i \le j_0}$ are at a distance at least $\epsilon$ from each other and $|\Phi_{j_0,\infty}(0)|^2 > \epsilon $ 
	with probability at least $1-\delta_0/2$. Let us denote by $\Omega_0$ this set.
	Up to extracting a further subsequence from $\tau_k$ (that we still denote by $\tau_k$), we can assume that $\mu_{\tau_k}$ and $(\nu_{i,\tau_k},\Phi_{i,\tau_k})_{1 \leq i <j_0}$ converge almost surely to $\mu_{\infty}$ and $(\nu_{i,\infty},\Phi_{i,\infty})_{1 \leq i <j_0}$.  We denote by $\Omega_1$ this event of probability $1$.
	
		Then, we let $\Omega_2 := \limsup_k A_k$. We have $\P(\Omega_2) \geq \delta_0$ and therefore $\P(\Omega_0 \cap \Omega_1  \cap \Omega_2) \geq \delta_0/2$. On $\Omega_2$, we have a (random) subsequence $\tau_{k_\ell}$ from $\tau_k$ on which $A_{k_\ell}$ holds --  that is
	$$d((\nu_{j_0,\tau_{k_\ell}}, \Phi_{j_0,\tau_{k_\ell}}) ; (\nu_{j_0,\infty}, \Phi_{j_0,\infty})) > \rho_0\,.$$
	Let $F:\R\times\tilde{\cF}\to[0,1]$ be a continuous map satisfying $F(\nu_{j_0,\infty},\Phi_{j_0,\infty}) =1$ and
	$$ d((\lambda,\Psi);(\nu_{j_0,\infty},\Phi_{j_0,\infty})) > \rho_0\wedge \epsilon/2 \Rightarrow F((\lambda,\Psi)) = 0 \;.$$
	Note that $F$ is compactly supported in its first coordinate. Necessarily on $\Omega_1$, as $k\to\infty$
	$$\int F(\lambda,\Psi) \mu_{\tau_k}(d\lambda,d\Psi) \to \int F(\lambda,\Psi) \mu_{\infty}(d\lambda,d\Psi)\;.$$
On the one hand, on $\Omega_0$,
\begin{align}\label{ineqFgeqeps}
	 \int F(\lambda,\Psi) \mu_{\infty}(d\lambda,d\Psi) \ge |\Phi_{j_0,\infty}(0)|^2 F(\nu_{j_0,\infty},\Phi_{j_0,\infty}) \geq \epsilon \;.
\end{align}
	On the other hand, on the event $\Omega_2$ by construction we know that for all $j\ge j_0$ and all $\ell \ge 1$
	$$ d((\nu_{j,\tau_{k_\ell}},\Phi_{j,\tau_{k_\ell}});(\nu_{j_0,\infty},\Phi_{j_0,\infty})) > \rho_0\;,$$
	and therefore 
	$$ \int F(\lambda,\Psi) \mu_{\tau_{k_\ell}}(d\lambda,d\Psi) = \sum_{i=1}^{j_0-1} F((\nu_{i,\tau_{k_\ell}},\Phi_{i,\tau_{k_\ell}})) |\Phi_{i,\tau_{k_\ell}}(0)|^2\;.$$
	But using the continuity of the evaluation map at $0$ in the topology of $\tilde{\cF}$ and the separation of the first $j_0$ atoms, if $\omega \in \Omega_0 \cap \Omega_1 \cap \Omega_2$, this quantity converges to $0$. As it should also be at least $\epsilon$ from \eqref{ineqFgeqeps}, we get a contradiction on an event of positive probability.
\end{proof}

\subsection{Intensity measure}\label{subsec:intensity}

In this section, we prove Theorem \ref{Th:Intensity} and thus work under the assumption $\kappa >0$. The main step consists in establishing the identity of the statement for functions $G(\lambda,\Psi)$ that depend only on the restriction of $\Psi$ to some bounded interval $[-A,A]$: to that end, we will rely on the convergence of $\mu_\tau$ to $\mu_\infty$ and on the GMP formula established in Proposition \ref{Prop:GMP2}.

\begin{proof}[Proof of Theorem \ref{Th:Intensity}]
	We recall the identity we aim for
	\begin{equation}\label{Eq:IdIntensity}
		\E\Big[\sum_{j\in \Z} | \Psi_{j,\infty}(0)|^2 G\Big(\lambda_{j,\infty},\Psi_{j,\infty}\Big)\Big]= \int \E\bigg[G\bigg(\gl, \frac{S(\cdot)}{\| S\|_{L^2(\R)}} \begin{pmatrix} \sin \gamma^{(\lambda)} (\cdot) \\\cos \gamma^{(\lambda)} (\cdot) \end{pmatrix}\bigg)\bigg] \frac{d\gl}{\pi}\;.
	\end{equation}
	By approximation, it suffices to prove it for any bounded continuous function $G:\R\times\tilde{\cF}\to\R$ compactly supported in its first coordinate. Fix such a function $G$. For any $A>0$, we introduce a cutoff function $\chi_A$ by setting $\chi_A(\cdot) := \chi(\cdot/A)$ where $\chi:\R\to [0,1]$ is a smooth function that vanishes outside $[-1,1]$ and equals one on $[-1/2,1/2]$. Then we set $G_A(\lambda,\Psi) := G(\lambda,\chi_A \Psi)$ and we observe that $G_A$ is a bounded continuous function from $\R\times\tilde{\cF}$ to $\R$.\\
	\textit{Step 1: from $G_A$ to $G$.} Assume that~\eqref{Eq:IdIntensity} holds for $G_A$. Let us show that it also holds for $G$. The exponential decay of $S$ ensures that a.s.
	$$ \sup_{t\in \R} e^{\frac{\kappa^2 |t|}{2}}(1-\chi_A(t)) \frac{|S(t)|}{\| S\|_{L^2(\R)}} \to 0\;,\quad A\to\infty\;.$$
	By passing to the limit $A\to\infty$, it is straightforward to check that the r.h.s.of~\eqref{Eq:IdIntensity}~converges to the desired quantity. Let us now concentrate on the l.h.s. Since $G$ is compactly supported in its first coordinate, there exist $a>0$ such that
	$$ \sum_{j\in \Z} | \Psi_{j,\infty}(0)|^2 \bigg| G\Big(\lambda_{j,\infty},\Psi_{j,\infty}\Big)-G_A\Big(\lambda_{j,\infty},\Psi_{j,\infty}\Big)\bigg|\le 2\|G\|_\infty \sum_{j\in \Z} | \Psi_{j,\infty}(0)|^2 \mathbf{1}_{[-a,a]}(\lambda_{j,\infty})\;,$$
	and, by Theorem \ref{Th:M}, the r.v.~on the r.h.s.~is bounded in $L^1$. Furthermore, by Theorem \ref{Th:Localization} the r.v.~on the l.h.s.~converges almost surely to $0$ as $A\to\infty$. By the Dominated Convergence Theorem the convergence also holds in $L^1$ and this suffices to extend the identity of the statement to the function $G$.\\
	\textit{Step 2: identity on $G_A$.} We now establish the identity for $G_A$. The convergence in $L^p$ stated in Theorem \ref{Th:CVrandommeas} ensures that as $\tau\uparrow\infty$
	\begin{equation}\label{Eq:CVIntensity}
		\E\Big[\sum_{j\in\Z} | \Psi_{j,\tau}(0)|^2 G_A(\lambda_{j,\tau},\Psi_{j,\tau})\Big] \to \E\Big[\sum_{j\in\Z} | \Psi_{j,\infty}(0)|^2 G_A(\lambda_{j,\infty},\Psi_{j,\infty})\Big]\;.
	\end{equation}
	This convergence holds for any b.c.~imposed on the operator $\cS_{\tau}$. For convenience, we pick~$b_\pm=e_2$, or equivalently $\beta_\pm = 0$. We set
	\begin{align*}
		U_{\tau,n} = \bin{\sin \beta^{(n\pi)}_\tau (\cdot) }{\cos \beta^{(n\pi)}_\tau (\cdot) }, \qquad U = \bin{\sin \gamma^{(\lambda)}}{\cos \gamma^{(\lambda)}}\,.
	\end{align*}
	Given the definition of $\tilde{\cF}$, $G_A:\R\times\tilde{\cF}\to\R$ can be seen as a bounded continuous map from $\R\times\cF$ into $\R$ that satisfies $G_A(\lambda,\Psi)=G_A(\lambda,-\Psi)$. The quantity on the l.h.s.~of~\eqref{Eq:CVIntensity} was computed in Proposition \ref{Prop:GMP2}:
	\begin{align*}
		\E\Big[\sum_{j\in\Z} | \Psi_{j,\tau}(0)|^2 G_A(\lambda_{j,\tau},\Psi_{j,\tau})\Big]= \int_{\lambda \in\R}\sum_{n\in\Z}  \E\Big[ G_A\Big(\lambda , \frac{S(\cdot)}{\| S\|_{\tau}} U_{\tau,n}\Big)\Big] g_{n\pi,(\kappa_0^2+\kappa^2)\tau}(\lambda \tau)d\lambda\;.
	\end{align*}
	Given $\lambda \in\R$, the Gaussian factor on the r.h.s.~is maximized for $n$ close to $n_0(\lambda,\tau) := \lfloor \frac{\lambda \tau}{\pi} \rfloor$. We will make use of the Gaussian estimates collected in Lemma \ref{Lemma:n0}: in particular, we will neglect the contributions coming from those $n$ that lie outside $I_{n_0,\tau} := [n_0- \tau^{2/3},n_0+\tau^{2/3}] \cap \Z$. Thus it suffices to pass to the limit on
	$$ \int_{\lambda\in\R}\sum_{n\in I_{n_0,\tau}}  \E\bigg[ G_A\bigg(\lambda , \frac{S(\cdot)}{\| S\|_{\tau}} U_{\tau,n}\bigg)\bigg] g_{n\pi,(\kappa_0^2+\kappa^2)\tau}(\lambda \tau)d\lambda\;.$$
	We wish to apply the Dominated Convergence Theorem on the integral in $\lambda$. From the assumptions made on $G$, there exist $K,a >0$ such that $|G(\lambda,\Psi)| \leq K \mathbf{1}_{[-a,a]}(\lambda)$. Therefore
	\begin{align*}
		&\sum_{n\in I_{n_0,\tau}}  \bigg| \E\bigg[ G_A\bigg(\lambda , \frac{S_\tau(\cdot)}{\| S_\tau\|_{L^2}} U_{\tau,n}\bigg)\bigg] \bigg| g_{n\pi,(\kappa_0^2+\kappa^2)\tau}(\lambda \tau)\le  K \mathbf{1}_{[-a,a]}(\lambda) \sum_{n\in \Z} g_{n\pi,(\kappa_0^2+\kappa^2)\tau}(\lambda \tau)\;.
	\end{align*}
	For all $\tau$ large enough, given Lemma \ref{Lemma:n0}, this last term is bounded by
	$$ \frac{2K}{\pi} \mathbf{1}_{[-a,a]}(\lambda)\;,$$
	which is integrable in $\lambda$. This will serve as a domination estimate.\\
	Using Lemmas \ref{Lemma:n0} and \ref{Lemma:phase}, for any $\lambda \in\R$ we have as $\tau\to\infty$
	$$ \sum_{n\in I_{n_0,\tau}}   \E\bigg[ G_A\bigg(\lambda , \frac{S_\tau(\cdot)}{\| S_\tau\|_{L^2}} U_{\tau,n}\bigg)\bigg]  g_{n\pi,(\kappa_0^2+\kappa^2)\tau}(\lambda \tau) \to \frac1{\pi} \E\bigg[ G_A\bigg(\lambda , \frac{S(\cdot)}{\| S\|_{L^2}} U\bigg)\bigg]\;.$$
	The Dominated Convergence Theorem, combined with \eqref{Eq:CVIntensity}, then yields the desired identity.
\end{proof}
We now prove the convergence in law of the r.v.~appearing in the expectation
\begin{lemma}\label{Lemma:phase}
	Fix $\lambda \in \R$. Let $\gamma_0$ be a uniform r.v.~over $[0,2\pi]$, and $W_2$ be a two-sided Brownian motion. Assume that $\gamma_0$, $W_1$ and $W_2$ are independent, and define
	$$ \gamma^{(\lambda)}(t) := \gamma_0 + \sqrt{\kappa_0^2+\kappa^2}\,W_2(t) + \lambda t\;,\quad t\in\R\;.$$
	Then, as $\tau\to\infty$,
	\begin{equation}\label{Eq:ToProveH}
		\sup_{n\in I_{n_0,\tau}} \bigg| \E\bigg[ G_A\bigg(\lambda , \frac{S(\cdot)}{\| S\|_{\tau}} U_{\tau,n}\bigg)\bigg]- \E\bigg[ G_A\bigg(\lambda , \frac{S(\cdot)}{\| S\|_{L^2}} U\bigg)\bigg] \bigg| \to 0\;.
	\end{equation}
\end{lemma}
\begin{proof}
	We fix $\lambda$ in the whole proof. As $\tau\to\infty$, $S/\| S\|_{\tau}$ converges a.s.~to $S/\| S\|_{L^2}$ for the topology of uniform convergence (not only on compact sets, although this would have sufficed for our purpose). Therefore it suffices to prove~\eqref{Eq:ToProveH} with $\frac{S(\cdot)}{\| S\|_{\tau}}$ replaced by $\frac{S(\cdot)}{\| S\|_{L^2}}$.\\
	Suppose that for any given bounded continuous function $F$ from $\cC([-A,A],\R^2)$ into $\R$ the following holds
	\begin{align}\label{CVphase}
		\sup_{n\in I_{n_0,\tau}} &\bigg| \E\big[ F\big(U_{\tau,n}|_{[-A,A]}\big)\big]- \E\big[ F\big(U|_{[-A,A]}\big)\big] \bigg| \to 0\;.
	\end{align}
	Since $S$ is independent of both $U_{\tau,n}$ and $U$, this is precisely the desired convergence provided one takes $F=F_S$ where
	$$ F_S(u) := \E\Big[G_A\Big(\lambda,\frac{S(\cdot)}{\| S\|_{L^2}} u(\cdot)\Big)\Big]\;.$$
	It remains to prove \eqref{CVphase}.
	
	As far as this convergence is concerned, we only care about $\beta^{(n\pi)}_\tau$ \emph{modulo} $2\pi$. We thus decompose $\beta^{(n\pi)}$ into the sum of three terms
	$$ \beta^{(n\pi)}(t) = \frac{n\pi t}{\tau} + \sqrt{\kappa_0^2+\kappa^2}\Big( Z_\tau(t) -  Z_\tau(0)\Big) + \Big( \frac{n\pi}{2} +  \sqrt{\kappa_0^2+\kappa^2} Z_\tau(0) \Big) \;.$$
	We are going to show that:\begin{enumerate}[label=(\roman*)]
		\item the first term converges to $\lambda t$ as $\tau\to\infty$, uniformly over all $n\in I_{n_0,\tau}$ and over all $t \in [-A,A]$,
		\item $Z_\tau(\cdot) -  Z_\tau(0)$ converges in law to $W_2$ in $\cC([-A,A], \R)$,
		\item the third term, modulo $2\pi$, converges in law to $\gamma_0$,
		\item the pair formed by $Z_\tau(\cdot) -  Z_\tau(0)$ and the third term, modulo $2\pi$, converges in law to the pair formed by $W_2$ and $\gamma_0$.
	\end{enumerate}
	Item (i) is straightforward. We turn to (ii). Tightness is straightforward from a moment estimate, and any limit point is necessarily a centered Gaussian process whose covariance can be computed explicitly and equals that of $W_2$. We now deal with (iii). We let $\{x\}_{2\pi}$ denote the unique real value in $[0,2\pi)$ such that $x-\{x\}_{2\pi}$ is a multiple of $2\pi$ for any given $x\in\R$. Fix $\alpha \in \R$. Let us define $\tilde{Z} := \frac{n\pi}{2} + \sqrt{\kappa_0^2+\kappa^2} Z_\tau(0)$, which has law $\cN(n\pi/2,(\kappa_0^2+\kappa^2)\tau/4)$. Let us abbreviate by $g$ the density function $g_{n\pi/2,(\kappa_0^2+\kappa^2)\tau/4}$. Then
	\begin{align*}
		\E\Big[e^{i\alpha \{ \tilde{Z}\}_{2\pi}} \Big] &= \sum_{k\in\Z} \E\bigg[e^{i\alpha ( \tilde{Z} - 2k\pi)} \mathbf{1}_{\{ \tilde{Z} \in [2k\pi, 2(k+1) \pi)\}}\bigg]= \sum_{k\in\Z} \int_0^{2\pi} e^{i\alpha x} g(2k\pi + x) dx\\
		&= \frac1{2\pi}\int_0^{2\pi} e^{i\alpha x} dx + R(\alpha,\tau)\;,
	\end{align*}
	the first term is the characteristic function of a uniform r.v.~over $[0,2\pi]$, and the second is a remainder that we need to control:
	\begin{align*}
		R(\alpha,\tau) &:= \int_0^{2\pi} e^{i\alpha x} \Big(\sum_{k\in\Z} g(2k\pi + x) - \frac1{2\pi}\Big) dx\\
		&=\int_0^{2\pi} e^{i\alpha x} \sum_{k\in\Z} \Big(g(2k\pi + x) - \frac1{2\pi} \int_{y\in[0,2\pi]} g(2k\pi + y) dy\Big) dx\;,
	\end{align*}
	where we have used at the second line the fact that $g$ is a density function. We have
	\begin{align*}
		\Big| g(2k\pi + x) - \frac1{2\pi} \int_{y\in[0,2\pi]} g(2k\pi + y) dy \Big|&\le \frac1{2\pi} \int_{y\in[0,2\pi]} \Big| g(2k\pi + x) - g(2k\pi + y)\Big| dy\\
		&\lesssim \sup_{z\in [0,2\pi]} |g'(2k\pi + z)|\;,
	\end{align*}
	uniformly over all $k\in\Z$ and all $x\in[0,2\pi)$. It is not hard to check that $|R(\alpha,\tau)|$ is bounded by a term of order $1/\sqrt{\tau}$, and is therefore negligible in the limit $\tau\to\infty$. Item (iii) is thus proved.\\
	We turn to (iv). It suffices to show asymptotic independence of $\{\frac{n\pi}{2} +  \sqrt{\kappa_0^2+\kappa^2}Z_\tau(0) \}_{2\pi}$ with the finite-dimensional marginals of the process $(Z_\tau(t) - Z_\tau(0), t\in\R)$. We provide the details for the one-dimensional marginals, the general case follows from the same arguments. Fix $\alpha, \beta \in \R$ and $t\in\R$. Set $X = Z_\tau(0)$ and $Y=Z_\tau(t) - Z_\tau(0)$. We aim at showing that as $\tau\to\infty$
	$$ \E[e^{i\alpha \{\frac{n\pi}{2}  + \sqrt{\kappa_0^2+\kappa^2} X\}_{2\pi} + i \beta Y}] \to \E[e^{i\alpha \gamma_0}] e^{-\frac{\beta^2 |t|}{2}}\;.$$
	The covariance structure of $Z_\tau$ ensures that there exists a r.v.~$V \sim\cN(0,1)$ independent of $X$ such that $Y = -(2|t|/\tau)X + \sqrt{|t| - (2t^2/\tau)} V$. Consequently
	$$  \E\Big[e^{i\alpha \{\frac{n\pi}{2}  + \sqrt{\kappa_0^2+\kappa^2}X\}_{2\pi} + i \beta Y}\Big] = \E\Big[e^{i\alpha \{\frac{n\pi}{2}  + \sqrt{\kappa_0^2+\kappa^2}X\}_{2\pi} - i (2|t| \beta/\tau) X}\Big] \E\Big[e^{i \beta \sqrt{|t| - (2t^2/\tau)} V}\Big]\;.$$
	The second term on the r.h.s.~satisfies
	$$ \E\Big[e^{i \beta \sqrt{|t| - (2t^2/\tau)} V}\Big] = \exp(-\frac{\beta^2(|t| - (2t^2/\tau))}{2}) \to \exp(-\frac{\beta^2 |t|}{2})\;,$$
	as $\tau\to\infty$. Regarding the first term, since the variance of $X$ equals $\tau/4$ we have
	$$ \Big|\E\Big[e^{i\alpha \{\frac{n\pi}{2} + \sqrt{\kappa_0^2+\kappa^2}X\}_{2\pi} - i (2|t| \beta/\tau) X} \mathbf{1}_{\{|X| > \tau^{2/3}\}}\Big]\Big| \le \P(|X| > \tau^{2/3}) \to 0\;,\quad \tau\to\infty\;.$$
	On the other hand,
	\begin{align*}
		\E\Big[e^{i\alpha \{\frac{n\pi}{2} + \sqrt{\kappa_0^2+\kappa^2}X\}_{2\pi} - i (2|t| \beta/\tau) X} \mathbf{1}_{\{|X| \le \tau^{2/3}\}}\Big] &= \E\Big[e^{i\alpha \{\frac{n\pi}{2} + \sqrt{\kappa_0^2+\kappa^2}X\}_{2\pi}} \mathbf{1}_{\{|X| \le \tau^{2/3}\}}\Big] + o(1)\\
		&=\E[e^{i\alpha \{\frac{n\pi}{2} + \sqrt{\kappa_0^2+\kappa^2}X\}_{2\pi}}] + o(1)\;,
	\end{align*}
	and we already saw that this quantity goes to $\E[e^{i\alpha \gamma_0}]$ as $\tau\to\infty$.
\end{proof}

\subsection{Some consequences of Theorem \ref{Th:Intensity}}

In this paragraph, we use Theorem \ref{Th:Intensity} to obtain two results: the density-density correlator expression of Corollary \ref{Cor:correlator} and the fact that the spectrum of $\cS_{\infty}$ is $\R$. Let us first prove Corollary \ref{Cor:correlator}. 
\begin{proof}
	We apply Theorem \ref{Th:Intensity} with $G(\gl,\Psi) := 1_{A}(E - \gl) |\Psi(r)|^2$ and any Borel set $A \subset \R$
	\begin{align*}
		\E\Big[\sum_{j \in \Z} |\Psi_{j,\infty}(0)|^2 |\Psi_{j,\infty}(r)|^2 1_{E - \gl_{j,\infty} \in A} \Big] &= \frac{1}{\pi}\int 1_{E-\gl \in A} \E\Big[\frac{S^2(r)}{\|S\|^2_{L^2(\R)}}\Big] d\gl \\
		&=\frac{1}{\pi} |A| \times \E\Big[\frac{S^2(r)}{\|S\|^2_{L^2(\R)}}\Big] \,.
	\end{align*}
	This implies Corollary \ref{Cor:correlator} in the special case $t=0$. Its extension to any $t \in \R$ follows from spatial invariance given by Proposition \ref{Prop:SymmetrySpatialShift}: under the shift $T_t$ each eigenfunction $\Psi_{j,\infty}(\cdot)$ is mapped to $\Psi_{j,\infty}(t + \cdot)$ while the eigenvalues remain unchanged.
\end{proof}

It is natural to conjecture the following formula, which matches the one found by Kolokolov \cite{Kolokolov1993} for the high energy of the 1d-Schr\"odinger operator with white noise potential (see~\cite[Formula (70)]{Kolokolov1993}: the expression below corresponds to the large $E$ limit of $2\pi \sqrt{E} \; p_{E}(x,x')$ using the notations therein).
\begin{conjecture}\label{Conj:Kolokolov}
	We have the identity
	\begin{align*}
		\E\bigg[
		\frac{S(r)^2}{\|S\|_{L^2(\R)}^2}
		\bigg] =
		4\pi^2\kappa^2
		e^{-\frac{\kappa^2}{2}|r|}
		\int_0^\infty
		\frac{t\sinh(\pi t)}{\cosh^3(\pi t)}
		\left(t^2+\frac14\right)^2
		e^{-2\kappa^2t^2|r|}
		\, dt \,.
	\end{align*} 
\end{conjecture}

Taking $r=0$, we recover the identity
\begin{align*}
		\E\bigg[
	\frac{1}{\|S\|_{L^2(\R)}^2}
	\bigg] =
	4\pi^2\kappa^2
	\int_0^\infty
	\frac{t\sinh(\pi t)}{\cosh^3(\pi t)}
	\left(t^2+\frac14\right)^2
	\, dt \,.
\end{align*}
which can be established using~\cite[Formula (6.2)]{MatsumotoYor2005}.
%which is equivalent to
%\begin{align*}
%	\E\bigg[
%	\frac{1}{\int_{-\infty}^{\infty} \exp(-2|t| + 2 W_1(t)) dt}
%	\bigg] =
%	4\pi^2
%	\int_0^\infty
%	\frac{t\sinh(\pi t)}{\cosh^3(\pi t)}
%	\left(t^2+\frac14\right)^2
%	\, dt \,.
%\end{align*}
%
%Theorem 6.2 in Matsumoto Yor. It is proved that if $\mu >0$ and $A_{\mu} := \int_0^{\infty} \exp(2 B(t) - 2 \mu t) dt$ then $(2\gamma_\mu)^{-1}$, where $\gamma_\mu$
%denotes a gamma random variable with parameter $\mu$.
%
%If $Y \sim \int_{0}^{\infty} \exp(-2 \, t + 2 W_1(t)) dt$ and $Y'$ is an independent copy of $Y$, then
%\begin{align*}
%	\E[\frac{1}{Y + Y'}] &= \int_0^\infty \E[\exp(- t (Y+Y'))] dt\\
%	&=\int_0^\infty \E[\exp(- t Y)]^2 dt\\
%	&= \int_0^\infty \E[\exp(- t A_{\infty}^{(-1)})]^2 dt\\
%	&= \int_0^\infty (\sqrt{2t} K_1(\sqrt{2t}))^2 dt \\
%	&=  2\int_0^\infty t (K_1(\sqrt{2t}))^2 dt 
%\end{align*}
%using formula (6.2) of Matsumoto Yor.

\medskip
We now prove that the spectrum of $\cS_\infty$ equals $\R$, thus completing the proof of Theorem \ref{Th:Localization}.
\begin{proof}[Proof of Theorem \ref{Th:Localization} - second part]
 In the proof of Proposition \ref{Prop:SymmetrySpatialShift}, we showed that the operator $T_{t_0} \cS_\infty  T_{t_0}^{-1}$ coincides with the operator associated with the shifted noise $(\kappa_0 \dot{B}_0(\cdot+t_0), \kappa \dot{B}_1(\cdot+t_0), \kappa \dot{B}_2(\cdot+t_0))$. Since white noise is ergodic, we deduce from~\cite[Lemma V.2.1]{CarLac} that the rank of any spectral projector is either a.s.~$0$ or a.s.~$\infty$. Fix some non-empty interval $[a,b] \subset \R$. By~\eqref{Eq:IntensityLambda} of Theorem \ref{Th:Intensity}, we deduce that $$\P( \#\{j \in \Z \,:\, \lambda_{j,\infty} \in [a,b]\} = 0) < 1 \,,$$ so the rank of the spectral projector on $[a,b]$ is not a.s.~$0$, consequently it equals $+\infty$ a.s. Therefore a.s.~the spectrum contains all non-empty intervals $[a,b]$ with $a,b\in\Q$ so it equals $\R$.
\end{proof}

\subsection{The case $\kappa=0$}\label{Subsec:k=0}

\begin{proof}[Proof of Proposition \ref{Th:Spectrumk=0}]
	We take $\kappa =0$. By Proposition \ref{Prop:Variousk}, the operators $\cS_{\infty,\kappa_0,0}$ and $\cS_{\infty,0,0}$ are unitary equivalent. Therefore, it suffices to prove that the latter is absolutely continuous. For any $h\in L^2(\R,\C)$, we define its Fourier transform
	$$ \hat{h}(\xi) := \frac1{\sqrt{2\pi}} \int_{\R} e^{it\xi} h(t) dt\;,\quad \xi\in\R\;.$$
	We define the unitary operator $A$ on $L^2(\R,\C^2)$
	$$ A\begin{pmatrix} f_1 \\ f_2\end{pmatrix} = \frac1{\sqrt{2}} \begin{pmatrix} \hat{f}_1(\xi) + i\hat{f}_2(\xi)  \\ \hat{f}_2(-\xi) + i\hat{f}_1(-\xi)\end{pmatrix}\;.$$
	Then the operator $\tilde{\cS} := A \cS_{\infty,0,0} A^{-1}$ is a self-adjoint operator on $L^2(\R,\C^2)$ that satisfies
	$$ \tilde{\cS} \begin{pmatrix} g_1(\xi) \\ g_2(\xi)\end{pmatrix} = \xi \begin{pmatrix} g_1(\xi) \\ g_2(\xi)\end{pmatrix}\;.$$
	We can thus compute the Stieltjes transform of the spectral measure associated to any $f\in L^2$: for any $z\in\C\backslash\R$
	$$ \langle f, (\cS_{\infty,0,0}-z)^{-1} f\rangle = \langle A f, (\tilde{\cS}-z)^{-1} A f\rangle = \int_{\R} |A f(\xi)|^2 (\xi-z)^{-1} d\xi\;.$$
	We deduce that the spectral measure associated to $f$ is absolutely continuous, with density equal to $|A f(\xi)|^2$. This proves that the spectrum of $\cS_{\infty,0,0}$ is absolutely continuous.
\end{proof}

\appendix

\section{Some technical estimates}\label{Appendix}
\subsection{Gaussian estimates}\label{Appendix:proofGaussian}

We prove here some useful estimates on the Gaussian kernel that appears in the GMP formula of Proposition \ref{Prop:GMP2}. Recall that $g_{a,b}$ is the density function of the Gaussian law of mean $a$ and variance $b$. To alleviate the notations, we will deal with $g_{\lambda \tau,\sigma^2\tau}(q+n\pi)$ for some arbitrary parameters $\sigma > 0$, $q,\lambda\in\R$. Given $\lambda,q \in \R$, this Gaussian factor is maximized for $n$ close to $n_0(\lambda,\tau) := \lfloor \frac{\lambda \tau - q}{\pi} \rfloor$. The lemma below states several estimates on its sum over $n\in\Z$.

\begin{lemma}\label{Lemma:n0}
	Fix $\sigma > 0$. There exist $c,C > 0$ such that for all $\lambda \in \R$, all $\tau > 1$ and all $q\in\R$
	$$ \sum_{n\in\Z} g_{\lambda \tau,\sigma^2\tau}(q+n\pi) \leq c\;,$$
	and
	$$ \sum_{n\in\Z: n \notin [n_0- \tau^{2/3},n_0+\tau^{2/3}]} g_{\lambda \tau,\sigma^2\tau}(q+n\pi) \leq c e^{-C\tau^{1/3}}\;.$$
	Furthermore, as $\tau \to \infty$
	$$ \sum_{n\in\Z} g_{\lambda \tau,\sigma^2\tau}(q+n\pi) \to \frac1{\pi}\;,$$
	uniformly over all $\lambda \in \R$ and all $q\in\R$. 
\end{lemma}
\begin{proof}
	We set $r_0 := \frac{\lambda \tau- q}{\pi} - n_0 \in [0,1)$. Using the identity
	$$g_{a,b}(x) = \frac1{\sqrt{b}} \,g_{0,1}\Big(\frac{x-a}{\sqrt{b}}\Big)\;,$$
	and the change of variable $m =n_0-n$ we obtain
	\begin{align*}
		\sum_{n\in\Z} g_{\lambda \tau,\sigma^2\tau}(q+n\pi) &= \sum_{n\in \Z} \frac{1}{\sigma\sqrt{\tau}} g_{0,1}\Big(\frac{q + n\pi - \lambda \tau}{\sigma\sqrt{\tau}}\Big)=\sum_{m\in \Z} \frac{1}{\sigma\sqrt{\tau}} g_{0,1}\Big(\frac{(m+r_0)\pi}{\sigma\sqrt{\tau}}\Big)\\
		&=A(\tau) + R(\tau,r_0)\;,
	\end{align*}
	where
	$$ A(\tau) := \sum_{m\in \Z} \frac{1}{\sigma\sqrt{\tau}} g_{0,1}\Big(\frac{m\pi}{\sigma\sqrt{\tau}}\Big)\;,$$
	and
	$$ R(\tau,r_0) := \sum_{m\in \Z} \frac{1}{\sigma\sqrt{\tau}} \bigg(g_{0,1}\Big(\frac{(m+r_0)\pi}{\sigma\sqrt{\tau}}\Big) - g_{0,1}\Big(\frac{m\pi}{\sigma\sqrt{\tau}}\Big)\bigg)\;.$$
	The factor $A(\tau)$ is independent of $\lambda$ and $q$, it converges as $\tau\to\infty$ to
	$$ \int_{\R} g_{0,1}(x\pi) dx = \frac1{\pi}\;,$$
	and it is bounded uniformly over all $\tau > 1$. On the other hand, $R(\tau,r_0)$ is a remainder term which is bounded by a term of order $1/\sqrt{\tau}$ uniformly over all $\lambda,q\in \R$ and $\tau >1$, and thus converges to $0$ as $\tau\to\infty$.\\
	%	Indeed
	%	\begin{align*}
		%		\sum_{n\in\Z} g_{a,b}(n\pi) &\le \sum_{n\in\Z:|n\pi - a| > \pi} g_{a,b}(n\pi) + \sum_{n\in\Z:|n\pi - a| \le \pi} g_{a,b}(n\pi)\\
		%		&\le \sum_{n\in\Z:n\pi > a+ \pi} \int_{[(n-1)\pi,n\pi]}g_{a,b}(x)dx + \sum_{n\in\Z:n\pi < a - \pi} \int_{[n\pi,(n+1)\pi]}g_{a,b}(x)dx + 3 g_{a,b}(a)\\
		%		&\le \int_\R g_{a,b}(x) dx + \frac{3}{\sqrt{2\pi b}}\;,
		%	\end{align*}
	Following the same computation as above we deduce that for all $\tau > 1$ and all $\lambda\in\R$ and all $q\in\R$
	\begin{align*}
		\sum_{n\in\Z: n \notin [n_0- \tau^{2/3},n_0+\tau^{2/3}]} g_{\lambda \tau,\sigma^2\tau}(q+n\pi) &\le \sum_{m\in \Z: |m| > \tau^{2/3}} \frac{1}{\sigma\sqrt{\tau}} g_{0,1}\Big(\frac{(m+r_0)\pi}{\sigma\sqrt{\tau}}\Big)\\
		&\le \sum_{m\in \Z: |m| > \frac12 \tau^{2/3}} \frac{1}{\sigma\sqrt{\tau}} g_{0,1}\Big(\frac{m\pi}{\sigma\sqrt{\tau}}\Big)\\
		&\le \frac1{\pi}\int_{x\in \R: |x| > \frac{\tau^{1/6}}{4\sigma}} g_{0,1}(x \pi) dx\;,
	\end{align*}
	and the asserted bound follows from standard bounds on the Gaussian density function.
\end{proof}

\subsection{Estimates on the exponential of Brownian motion}\label{Appendix:BMestimates}

We start with a lemma on the moments of the exponential of a Brownian motion with drift.
\begin{lemma}\label{Lem:momentexpBM}
	Recall that $W_1$ is a two-sided Brownian motion. Fix $\kappa,\alpha,p > 0$. We have
\begin{align*}
	\E\Big[\Big(\sup_{t \in \R} \exp(\kappa W_1(t) - \alpha |t|) \Big)^p\Big] < \infty \quad \mbox{iff} \quad p < 2\frac{\alpha}{\kappa^2}\,.
\end{align*}
\end{lemma}
\begin{proof}
	By symmetry, we restrict to $t\ge 0$. We have
	\begin{align*}
		\sup_{t \geq 0} (\kappa W_1(t) - \alpha t) \stackrel{(d)}{=} 	\sup_{t \geq 0} (W_1(t) - \frac{\alpha}{\kappa^2} t) 
	\end{align*}
	By~\cite[Formula 1.1.4 (1) in Chapter II.2 p.251]{Handbook}, for any $u\ge 0$
	\begin{align*}
		\P\Big(	\sup_{t \geq 0} (\kappa W_1(t) - \alpha t) \geq u\Big) = \exp(- 2\frac{\alpha}{\kappa^2} u)\,.
	\end{align*}
	Therefore
	\begin{align*}
		\P \Big(	\sup_{t \geq 0} \exp(\kappa W_1(t) - \alpha t) \geq u \Big) = u^{- 2 \frac{\alpha}{\kappa^2} }\,,
	\end{align*}
	which immediately gives the result.
\end{proof}
Let $Y_\tau$ be either the two-sided Brownian motion $W_1$ on $\R$ (in which case $Y_\tau$ does not depend on $\tau$) or the Brownian bridge $Z_\tau$ on $[-\tau/2,\tau/2]$ introduced in Subsection \ref{Subsec:GMP}. In the next lemma, we establish a moment bound on the increments of $Y_\tau$.
\begin{lemma}\label{Lem:momentincrementsexpBM}
	For any $\alpha \in (0,1/2)$, any $\eta > 0$ and any $p > 0$
	\begin{align*}
		&\sup_{\tau > 1}\E \left[\sup_{t,t+s \in [-\tau/2,\tau/2], \; 0 < s \leq 1} \Big(e^{-\eta |t|} \frac{|Y_\tau(t+s) - Y_\tau(t)|}{s^{\alpha}}\Big)^p\right] < \infty\,.
	\end{align*}
	\end{lemma}
\begin{proof}
Decompose $\R$ into subintervals $I_j := [j,j+2]$ for $j \in \Z$ (so that every pair $(t,t+s)$ is at least in one $I_j$). Then we have
\begin{align*}
	\sup_{\substack{t,t+s \in [-\tau/2,\tau/2] \\0 < s \leq 1}} e^{-\eta |t|} \frac{|Y_\tau(t+s) - Y_\tau(t)|}{s^{\alpha}} &=\sup_{j\in\Z} \sup_{\substack{t,t+s \in [-\tau/2,\tau/2] \cap I_j \\ 0 < s \leq 1}} e^{-\eta |t|} \frac{|Y_\tau(t+s) - Y_\tau(t)|}{s^{\alpha}}  \\
	&\leq e^{2\eta} \sup_{j\in\Z}  e^{-\eta |j|}\sup_{\substack{t,t+s \in [-\tau/2,\tau/2] \cap I_j \\ 0 < s \leq 1}} \frac{|Y_\tau(t+s) - Y_\tau(t)|}{s^{\alpha}}  \,.
\end{align*}
Therefore
\begin{align*}
	&\E \left[\sup_{\substack{t,t+s \in [-\tau/2,\tau/2] \\ 0 < s \leq 1}} \Big(e^{-\eta |t|} \frac{|Y_\tau(t+s) - Y_\tau(t)|}{s^{\alpha}}\Big)^p\right] \\
	&\leq e^{2p \eta} \sum_{j \in \Z} e^{-\eta p|j|}\E\left[\sup_{\substack{t,t+s \in [-\tau/2,\tau/2]\cap I_j \\  0 < s \leq 1}} \Big(\frac{|Y_\tau(t+s) - Y_\tau(t)|}{s^{\alpha}}\Big)^p \right]\;.
\end{align*}
It can be checked that $\E[|Y_\tau(t+s) - Y_\tau(t)|^2] \leq s$ (for both $Y_\tau=W_1$ and $Y_\tau=Z_\tau$). Since this r.v.~is gaussian, we further get that for any $m \geq 1$, there exists a universal constant $C_m>0$ such that $\E[|Y_\tau(t+s) - Y_\tau(t)|^m] \leq C_m |s|^{m/2}$. By Kolmogorov Continuity Theorem~\cite[Theorem I.2.1]{RevuzYor}, for any $p\ge 1$ there exists a constant $c_p>0$ such that for all $\tau>1$ and all $j\in\Z$
\begin{align*}
	\E \left[ \Big(\sup_{t,t+s \in [-\tau/2,\tau/2]\cap I_j, \; 0 < s \leq 1} \frac{|Y_\tau(t+s) - Y_\tau(t)|}{s^{\alpha}}\Big)^p \right] < c_p\;,
\end{align*}
 and consequently
\begin{align*}
	&\sup_{\tau > 1} \E \left[\sup_{t,t+s \in [-\tau/2,\tau/2], \; 0 < s \leq 1} \Big(e^{-\eta |t|} \frac{|Y_\tau(t+s) - Y_\tau(t)|}{s^{\alpha}}\Big)^p\right] < \infty\,.
\end{align*}
This bound remains true for $p\in (0,1)$ thanks to the trivial inequality $x^p \le 1 + x$ for all $x\ge 0$.
\end{proof}

\subsection{Compact sets in $\cF$}\label{Appendix:Compact}

We first identify convenient relatively compact subsets of $\cF$.
\begin{lemma}\label{Lemma:Compact}
	For any $C>0$, we let $K(C)$ be the set of all $f\in\mathcal{F}$ for which there exists $T\in (1,+\infty]$ such that $f(t) = 0$ for all $t\notin [-T/2-1,T/2+1]$, $f$ is affine on $[-T/2-1,-T/2]$ and $[T/2,T/2+1]$, and
	$$ \sup_{t\in\R} |f(t)| e^{\frac{2|t| \kappa^2}{3}} + \sup_{t\in\R, s\in (0,1]: -T/2 \le t < t+s \le T/2} \frac{|f(t+s)e^{\frac{|t+s| \kappa^2}{2}}-f(t)e^{\frac{|t|\kappa^2}{2}}|}{s^{1/3}} \le C\;.$$
	Then $K(C)$ is relatively compact in the Banach space $\cF$.
\end{lemma}
In practice we will apply this material to the eigenfunctions in finite volume and we will take $T=\tau$. The exponents appearing in this set are arbitrary up to some extent: the first exponent $2/3$ in front of $\kappa^2$ must be taken in $(1/2,1)$ and the exponent $1/3$ in the denominator must be smaller than $1/2$ in line with the H\"older regularity of Brownian motion.
\begin{proof}
	Let $\cC_0(\R,\R^2)$ be the space of continuous functions from $\R$ to $\R^2$ vanishing at infinity and endowed with the supremum norm. By the locally compact version of the Arzela-Ascoli Theorem, a subset $J \subset \cC_0(\R,\R^2)$ is relatively compact if:\begin{enumerate}
		\item (Boundedness) $\sup_{g \in J} \| g\|_\infty <\infty$,
		\item (Uniform equicontinuity) for any $x\in \R$, the functions $g\in J$ are uniformly equicontinuous at $x$,
		\item (Uniform vanishing at infinity) $\sup_{g\in J} \sup_{t\notin [-A,A]} |g(t)| \to 0$ as $A\to\infty$.
	\end{enumerate}
	Consider the continuous bijection $\iota:g\mapsto f$ from $\cC_0(\R,\R^2)$ to $\cF$ defined through
	$$ f(t) := e^{-\frac{\kappa^2 |t|}{2}} g(t)\;,\quad t\in\R\;.$$
	Define $J(C):= \iota^{-1}(K(C))$. If $J(C)$ is relatively compact in $\cC_0(\R,\R^2)$, then $K(C)$ is relatively compact in $\cF$. It thus remains to prove that $J(C)$ is relatively compact. Note that $J(C)$ is the set of all $g\in \cC_0(\R,\R^2)$ for which there exists $T\in (1,+\infty]$ such that $g(t) = 0$ for all $t\notin [-T/2-1,T/2+1]$, $e^{-\frac{\kappa^2 |t|}{2}} g(t)$ is affine on $[-T/2-1,-T/2]$ and $[T/2,T/2+1]$, and
	$$ \sup_{t\in\R} |g(t)| e^{\frac{|t| \kappa^2}{6}} + \sup_{t\in\R, s\in (0,1]: -T/2 \le t < t+s \le T/2} \frac{|g(t+s)-g(t)|}{s^{1/3}} \le C\;.$$
	The boundedness and uniform vanishing at infinity are satisfied by $J(C)$. Let us check the uniform equicontinuity condition. Fix $g\in J(C)$ and let $f := \iota(g)$. For any $t,t+s \in [T/2,T/2+1]$
	\begin{align*}
		|g(t+s)-g(t)| &= |f(t+s) e^{\frac{\kappa^2 (t+s)}{2}} - f(t) e^{\frac{\kappa^2 t}{2}}|\\
		&\le e^{\frac{\kappa^2 (t+s)}{2}} |f(t+s)  - f(t)| + |f(t)| (e^{\frac{\kappa^2 (t+s)}{2}} - e^{\frac{\kappa^2 t}{2}})\;.
	\end{align*}
	Recall that $f$ is linear on $[T/2,T/2+1]$. Using an elementary bound on the increment of the exponential term, we deduce that
	$$ \frac{|g(t+s)-g(t)|}{s}  \le e^{\frac{\kappa^2 (T/2+1)}{2}} |f(T/2)| \Big(1+ \frac{\kappa^2}{2}\Big) \le C e^{-\frac{\kappa^2 T}{12}} e^{\frac{\kappa^2}{2}} \Big(1+ \frac{\kappa^2}{2}\Big)\;.$$
	A similar bound holds on $[-T/2-1,-T/2]$. Putting everything together, we deduce easily that there exists a constant $C'>0$ such that for all $g\in J(C)$
	$$ \sup_{t\in\R, s\in (0,1]} \frac{|g(t+s)-g(t)|}{s^{1/3}} \le C'\;,$$
	and the uniform equicontinuity holds.
\end{proof}
We proceed with the proof of the delicate tightness estimate stated in Lemma \ref{Lemma:Tight}.
\begin{proof}[Proof of Lemma \ref{Lemma:Tight}]
	In the previous lemma we identified relatively compact sets $K(C)$ in $\cF$. Since the map $\pi$ that associates to a function in $\cF$ its equivalence class in $\tilde{\cF}$ is $1$-Lipschitz, the set $\tilde{K}(C) := \pi(K(C))$ is relatively compact in $\tilde{\cF}$.\\
	For any $n\in\Z$ and any $\tau > 0$ we let $R_{\tau,n}:\R\to\R^2$ be the continuous function that vanishes outside $[-\tau/2-1,\tau/2+1]$, that satisfies
	$$ R_{\tau,n}(t) = \frac{S(t)}{{\|S\|}_{\tau}}U_{\tau,n}(t)\;,\quad t\in [-\tau/2,\tau/2]\;,$$
	and that is affine on $[\tau/2,\tau/2+1]$ and $[-\tau/2-1,-\tau/2]$. By Proposition \ref{Prop:GMP2}, the expectation in the statement equals
	$$ I_\tau := \int_{\lambda\in[a,b]}\sum_{n\in\Z}  A_{\tau,n}(C) g_{\lambda \tau,(\kappa_0^2+\kappa^2)\tau}(\beta_+-\beta_-+n\pi)d\lambda\;,$$
	where
	$$A_{\tau,n}(C):=\P\Big( R_{\tau,n} \in \cF\backslash K(C)\Big)\;,\quad U_{\tau,n}(\cdot) := \begin{pmatrix} \sin \beta^{(\beta_++n\pi)}_\tau (\cdot) \\\cos \beta^{(\beta_++n\pi)}_\tau (\cdot) \end{pmatrix}\;.$$
	Set $n_0(\lambda,\tau) := \lfloor \frac{\lambda \tau - \beta_+ + \beta_-}{\pi} \rfloor$ and write $I_\tau \le I_\tau^{(1)} + I_\tau^{(2)}$ where
	\begin{align*}
		I_\tau^{(1)} &:= \int_{\lambda\in[a,b]}\sum_{n\notin [n_0- \tau^{2/3},n_0+\tau^{2/3}]}  g_{\lambda \tau,(\kappa_0^2+\kappa^2)\tau}(\beta_+-\beta_-+n\pi)d\lambda\;,\\
		I_\tau^{(2)} &:= \int_{\lambda\in[a,b]}\sum_{n\in [n_0- \tau^{2/3},n_0+\tau^{2/3}]}  A_{\tau,n}(C) g_{\lambda \tau,(\kappa_0^2+\kappa^2)\tau}(\beta_+-\beta_-+n\pi)d\lambda\;.
	\end{align*}
	By Lemma \ref{Lemma:n0}, $\limsup_{\tau \to \infty} I_\tau^{(1)} = 0$. If we prove that as $C\to\infty$
	\begin{equation}\label{Eq:TightToProve}
		\sup_{n\in [n_0- \tau^{2/3},n_0+\tau^{2/3}]\cap \Z}\sup_{\tau > 1} A_{\tau,n}(C) \to 0\;,
	\end{equation}
	then Lemma \ref{Lemma:n0} will ensure that $\limsup_{C\to\infty} \limsup_{\tau \to \infty} I_\tau^{(2)} = 0$, and the proof will be complete.
	We now prove~\eqref{Eq:TightToProve}. Let $V_\tau:\R\to\R$ be the continuous function that coincides with $S$ on $[-\tau/2,\tau/2]$, vanishes outside $[-\tau/2-1,\tau/2+1]$ and is affine in between. Given the definition of $K(C)$, it is elementary to check that
	\begin{equation}
		\left\{R_{\tau,n} \in \cF\backslash K(C) \right\} \subset E_{\tau}^{(1)} \cup E_{\tau}^{(2)} \cup E_{\tau,n}^{(3)}\;,
	\end{equation}
	where
	\begin{align*}
		E_{\tau}^{(1)} &:= \left\{\sup_{t\in \R} V_\tau(t) e^{2\kappa^2 |t|/3} \geq \sqrt{C/2}\right\}\;,\\
		E_{\tau}^{(2)} &:= \left\{{\|S\|}_{\tau} \leq \sqrt{2/C}\right\}\;,\\
		E_{\tau,n}^{(3)} &:= \left\{\sup_{\stackrel{t,t+s\in [-\tau/2,\tau/2]}{s\in (0,1]}} \frac{\big|e^{\kappa W_1(t+s)-\frac{\kappa^2 |t+s|}{2}}U_{\tau,n}(t+s) - e^{\kappa W_1(t)-\frac{\kappa^2 |t|}{2}}U_{\tau,n}(t) \big|}{s^{1/3}} \geq \sqrt{C/2} \right\}\;.
	\end{align*}
	Observe that for all $t\in [-\tau/2,\tau/2]$
	$$ V_\tau(t) e^{2\kappa^2 |t|/3} = e^{\kappa W_1(t)-\kappa^2 |t|/3}\;.$$
	As the random variable $\sup_{t\in \R} e^{\kappa W_1(t)-\kappa^2 |t|/3}$ is almost surely finite, $\sup_{\tau > 1} \P(E_{\tau}^{(1)})$ goes to $0$ as $C\to\infty$. Regarding $E_{\tau}^{(2)}$, since for any $\tau > 1$
	$$ {\|S\|}_{\tau}^2 \ge \int_{[-1/2,1/2]} e^{2\kappa W_1(t)-2\kappa^2 |t|} dt\;,$$
	and since the r.v.~on the r.h.s.~is positive almost surely, we easily deduce that $\sup_{\tau > 1} \P(E_{\tau}^{(2)})$ goes to $0$ as $C\to\infty$. We now concentrate on $E_{\tau,n}^{(3)}$. Denote $X(t) := \kappa W_1(t)-\frac{\kappa^2 |t|}{2}$ and observe that
	\begin{align*}
		\big|e^{X(t+s)}U_{\tau,n}(t+s) - e^{X(t)}U_{\tau,n}(t) \big| \le e^{X(t+s)}\big|U_{\tau,n}(t+s) - U_{\tau,n}(t) \big| + \big|e^{X(t+s)} - e^{X(t)} \big|\;.
	\end{align*}
	
	Recall the definition of $U_{\tau,n}$ and observe that the functions $\sin, \cos$ are $1$-Lipschitz. To prove that 
	$$\sup_{n\in [n_0- \tau^{2/3},n_0+\tau^{2/3}]\cap \Z} \sup_{\tau > 1} \; \P(E_{\tau,n}^{(3)}) \to 0\;,\quad \mbox{ as } C\to\infty\;,$$
	it thus suffices to show that for some $r >0$,
	\begin{align}
		\sup_{n\in [n_0- \tau^{2/3},n_0+\tau^{2/3}]\cap \Z} \sup_{\tau > 1} \; \E\Big[\sup_{\stackrel{t,t+s\in [-\tau/2,\tau/2]}{s\in (0,1]}} e^{r X(t+s)}\frac{\big|\beta_{\tau}^{(\beta_+ + n \pi)}(t+s) - \beta_{\tau}^{(\beta_+ + n \pi)}(t) \big|^r}{s^{r/3}} \Big] < \infty \label{Boundphase}
	\end{align}
	and
	\begin{align}
		\sup_{n\in [n_0- \tau^{2/3},n_0+\tau^{2/3}]\cap \Z} \sup_{\tau > 1}\;  \E\Big[\sup_{\stackrel{t,t+s\in [-\tau/2,\tau/2]}{s\in (0,1]}} \frac{\big|e^{X(t+s)} - e^{X(t)} \big|^r}{s^{r/3}}\Big] < \infty\,. \label{Boundmodulus}
	\end{align}
	For the first term \eqref{Boundphase}, we introduce $\eta \in (0,\kappa^2/2)$ and use H\"older inequality: 
	\begin{align*}
		&\E\Big[\sup_{t,t+s\in [-\tau/2,\tau/2], s\in (0,1]} 	e^{r X(t+s)}\frac{\big|\beta_{\tau}^{(\beta_+ + n \pi)}(t+s) - \beta_{\tau}^{(\beta_+ + n \pi)}(t) \big|^r}{s^{r/3}} \Big]\\
		&\leq  \E\Big[\sup_{t\in \R} 	e^{2 r X(t) +2 \eta r  |t|}\Big]^{1/2} \E\Big[\sup_{\stackrel{t,t+s\in [-\tau/2,\tau/2]}{s \in (0,1]}} e^{-2 \eta r  |t + s|}\frac{\big|\beta_{\tau}^{(\beta_+ + n \pi)}(t+s) - \beta_{\tau}^{(\beta_+ + n \pi)}(t) \big|^{2r}}{s^{2r/3}} \Big]^{1/2}
	\end{align*}
	By Lemma \ref{Lem:momentexpBM}, the first factor in the RHS is finite as soon as $2 r  < 1- 2\eta/\kappa^2$.
	Regarding the second factor, we decompose $\beta_{\tau}^{(\beta_+ + n \pi)}$ into its drift part and its Brownian bridge part. The drift is uniformly Lipschitz for all $\tau >1$ and $n \in [n_0 - \tau^{2/3},n_0 + \tau^{2/3}]$. For the Brownian bridge part using Lemma \ref{Lem:momentincrementsexpBM}
	we conclude to the finiteness of \eqref{Boundphase}.
	
	For the second term \eqref{Boundmodulus}, using the simple inequality,
	\begin{align*}
		|e^{X(t+s)} - e^{X(t)} | \leq e^{X(t + s) \vee X(t)} 	|X(t+s) -X(t) |
	\end{align*}
	and both Lemma \ref{Lem:momentexpBM} and Lemma \ref{Lem:momentincrementsexpBM}, we can conclude.
\end{proof}

\bibliographystyle{Martin}
\bibliography{library}

\end{document}